# Multi-Level Stochastic Gradient Methods for Nested Composition Optimization

Shuoguang Yang* Mengdi Wang† Ethan X. Fang‡

**Abstract**

Stochastic gradient methods are scalable for solving large-scale optimization problems that involve empirical expectations of loss functions. Existing results mainly apply to optimization problems where the objectives are one- or two-level expectations. In this paper, we consider the multi-level compositional optimization problem that involves compositions of multi-level component functions and nested expectations over a random path. It finds applications in risk-averse optimization and sequential planning. We propose a class of multi-level stochastic gradient methods that are motivated from the method of multi-timescale stochastic approximation. First we propose a basic $T$-level stochastic compositional gradient algorithm, establish its almost sure convergence and obtain an $n$-iteration error bound $\mathcal{O}(n^{-1/2^T})$. Then we develop accelerated multi-level stochastic gradient methods by using an extrapolation-interpolation scheme to take advantage of the smoothness of individual component functions. When all component functions are smooth, we show that the convergence rate improves to $\mathcal{O}(n^{-4/(7+T)})$ for general objectives and $\mathcal{O}(n^{-4/(3+T)})$ for strongly convex objectives. We also provide almost sure convergence and rate of convergence results for nonconvex problems. The proposed methods and theoretical results are validated using numerical experiments.

## 1 Introduction

Over the past decade, stochastic gradient-type methods have drawn significant attention from various communities such as mathematical programming, signal processing and machine learning, mainly due to their practical efficiency in minimizing expected-value objective functions or empirical sums of a large number of loss functions [2, 5, 10, 11, 12, 14, 18, 22, 30]. They are particularly popular methods for tackling large-scale problems such as statistical estimation [6, 20], matrix and tensor factorization [8] and training deep neural networks [13, 29]. Stochastic gradient methods mainly apply to minimizing the expectation of a stochastic function, i.e.,

$$\min_x \mathbb{E}_\omega[f_\omega(x)],$$

where the expectation is taken over a random variable $\omega$. Note that this problem involves one level of expectation.

*Department of Industrial Engineering and Operations Research, Columbia University, New York City, NY, USA; e-mail: sy2614@columbia.edu

†Department of Operations Research and Financial Engineering, Princeton University, Princeton, NJ, USA; email: mengdiw@princeton.edu

‡Department of Statistics and Department of Industrial and Manufacturing Engineering, Pennsylvania State University, University Park, PA, USA; email: xxf13@psu.edu

In this paper, we propose to study a richer class of stochastic optimization problems, which involve nested expectations over a sequence of random variables. In particular, we consider the *$T$-level stochastic compositional optimization problem*, given by

$$\min_{x\in\mathcal{X}} F(x) = \mathbb{E}_{\omega_1}\left[f^{(1)}_{\omega_1}\left(\mathbb{E}_{\omega_2}\left[f^{(2)}_{\omega_2}\left(\cdots\left(\mathbb{E}_{\omega_T}[f^{(T)}_{\omega_T}(x)]\right)\cdots\right)\right]\right)\right], \tag{1.1}$$

where $f^{(j)}_{\omega_j}(\cdot) : \mathbb{R}^{d_j} \mapsto \mathbb{R}^{d_{j-1}}$ for $j = 1, \cdots, T$ are continuous mappings, $\mathcal{X}$ is a convex and closed set, and $d_0 = 1$, i.e., $F(x)$ is a real-valued function. The nested composition structure provides a rich modeling tool for data analysis and decision-making applications. For instance, online principal component analysis and policy evaluation in reinforcement learning can be formulated into two-level stochastic compositional optimization [15, 27]. For more applications, we illustrate two examples that arise from operations research in Section 4. One example is a mean-deviation risk-averse optimization problem which can be formulated into a 3-level compositional problem [1, 21]. The other example is related to smooth approximations of multi-stage stochastic programming [24].

In problem (1.1), for each $f^{(j)}_{\omega_j}$, we use the subscript $\omega_j$ to denote a random variable and use the superscript $(j)$ to denote its level. We focus on situations whether there exist deterministic functions $f^{(1)}, \ldots, f^{(T)}$ such that

$$f^{(j)}(x_j) = \mathbb{E}[f^{(j)}_{\omega_j}(x_j)|\omega_1, \cdots, \omega_{j-1}],$$

*for all $j = 1, \ldots, T$ with probability 1.* We refer to $f^1, \ldots, f^T$ as *component functions.* However, these component functions are not explicitly known to us. Note that the multi-level random variables $\omega_1, \ldots, \omega_T$ are not necessarily independent of one another. When we sample from their joint distribution, we may generate a sample path $(\omega_1, \ldots, \omega_T)$ sequentially by sampling each $\omega_j$ conditioned on realizations at the previous level's $(\omega_1, \ldots, \omega_{j-1})$. Throughout this paper, we assume that the component functions $f^1, \ldots, f^T$ are continuous and that there exists at least one optimal solution $x^*$ to problem (1.1). In some part of our analysis, we require the overall objective function $F(x)$ be convex, but we *never* require that any individual component function $f^{(j)}_{\omega_j}(\cdot)$ be convex, linear or monotone. We say that a function $f$ is "smooth" if it has Lipschitz continuous gradients, and say that it is "non-smooth" otherwise.

Our goal is to solve the $T$-level stochastic compositional optimization problem (1.1) by sampling multiple paths of $(\omega_1, \ldots, \omega_T)$. We are interested in scenarios where we do not have the explicit knowledge of the expected-value component functions $f^{(j)}$'s. This often occurs when evaluating $f^{(j)}$ requires making expensive passes over large data sets. This also occurs in online learning applications where $f^{(j)}$ can not be accurately calculated using finitely many samples. Instead of knowing $f^{(j)}$'s, we suppose that there is a Sample Oracle ($\mathcal{SO}$) such that:

- Upon each query ($x \in \mathcal{X}, y_1 \in \mathbb{R}^{d_1}, \ldots, y_T \in \mathbb{R}^{d_T}$), the $\mathcal{SO}$ generates a sample path $(\omega_1, \ldots, \omega_T)$ independently from the query.
- The $\mathcal{SO}$ returns a vector $f^{(T)}_{\omega_T}(x) \in \mathbb{R}^{d_{T-1}}$ and a noisy gradient/subgradient $\widetilde{\nabla} f^{(T)}_{\omega_T}(x) \in \mathbb{R}^{d_T \times d_{T-1}}$.
- The $\mathcal{SO}$ returns a vector $f^{(j)}_{\omega_j}(y_j) \in \mathbb{R}^{d_j}$ and a noisy gradient $\nabla f^{(j)}_{\omega_j}(y_j) \in \mathbb{R}^{d_j \times d_{j-1}}$.
- The $\mathcal{SO}$ returns a noisy gradient $\nabla f^{(1)}_{\omega_1}(y_1) \in \mathbb{R}^{d_1}$.

In the above, we denote by $\widetilde{\nabla} f^{(T)}_{\omega_T}(x)$ a noisy gradient/subgradient, which is to be specified in the context. Let us emphasize that this $\mathcal{SO}$ does *not* return unbiased first-order information regarding the overall objective function. The $\mathcal{SO}$ can be viewed as a *component-wise stochastic first-order oracle* that returns noisy first-order information for individual component functions $f^{(j)}$'s.Detailed assumptions on the $\mathcal{SO}$ will be specified later.

One might attempt to apply the sample average approximation (SAA) method to attack the multi-level expectation problem (1.1). However, replacing the nested expectations with empirical averages will not solve the optimization problem. It will reduce one problem with expectations to another one with empirical expectations. However, the two problems share similar structures and the latter one is not necessarily easier to solve. What we need is an implementable algorithm that computes the optimal solution by iteratively querying the $\mathcal{SO}$ and making efficient updates.

Another attempt would be to use some version of gradient method or stochastic gradient method. Stochastic gradient method will not work automatically. The main challenge is that we do not have access to the unbiased sample gradient of $F$ due to the multi-level nested expectations. To see this, let us consider the case where each $f^{(j)}$ is differentiable and apply the chain rule to get

$$\nabla F(x) = \nabla f^{(T)}(x) \nabla f^{(T-1)}(f^T(x)) \cdots \nabla f^{(1)} \left( f^{(2)} \circ \cdots \circ f^{(T)}(x) \right).$$

For a given $x \in \mathcal{S}$ and a given sample path $(\omega_1, \ldots, \omega_T)$, one may formulate an unbiased estimate of $\nabla F(x)$ as

$$\nabla f^{(T)}_{\omega_T}(x) \nabla f^{(T-1)}_{\omega_{T-1}}(f^{(T)}(x)) \cdots \nabla f^{(1)}_{\omega_1} \left( f^{(2)} \circ \cdots \circ f^{(T)}(x) \right),$$

which unfortunately cannot be calculated by calling the $\mathcal{SO}$ once (or even finitely many times). This is because that computing the preceding unbiased gradient sample requires quering the $\mathcal{SO}$ at values $f^{(T)}(x), f^{(T-1)} \circ f^{(T)}(x), \ldots, f^{(2)} \circ \cdots \circ f^{(T)}(x)$, which are unfortunately not known. As a result, the nested composition structure induces substantial bias in the sample gradients for $F$ as long as $T \geq 2$. In contrast, when $T = 1$, the objective function is linear in the distribution of the random variable $\omega$. For problems with $T \geq 2$, the nonlinear composition between expectations and component functions creates an objective function that is highly nonlinear with respect to the joint probability distribution of $\omega_1, \ldots, \omega_T$. A graphical illustration of the level of difficulty for dealing with multi-level composition optimization is given in Figure 1. We can view the optimization problem (1.1) under the $\mathcal{SO}$ as a form of estimation problem, in which we want to estimate the optimal solution $x^*$ by taking independent sample paths. We can see that the nonlinear composition makes this estimation/optimization problem fundamentally challenging.

Existing work on stochastic compositional optimization traces back to [7] which considered the two-level problem. In Section 6.7 of [7], a two-timescale stochastic approximation scheme was proposed and its almost sure convergence was established without rate analysis. Recently, [26] developed a general class of stochastic compositional gradient descent (SCGD) method for two-level problems and established convergence rate results under various assumptions. [28] developed an accelerated stochastic compositional proximal gradient (ASC-PG) method for the two-level problem and proved faster convergence in some cases. [16] considered a special case of the two-level problem where each expectation takes the form of a finite sum of loss functions and developed variance-reduced versions of the compositional gradient methods. However, to the best of our knowledge, all existing results only apply to the case where $T = 1, 2$. Multi-level stochastic compositional optimization remains largely open.

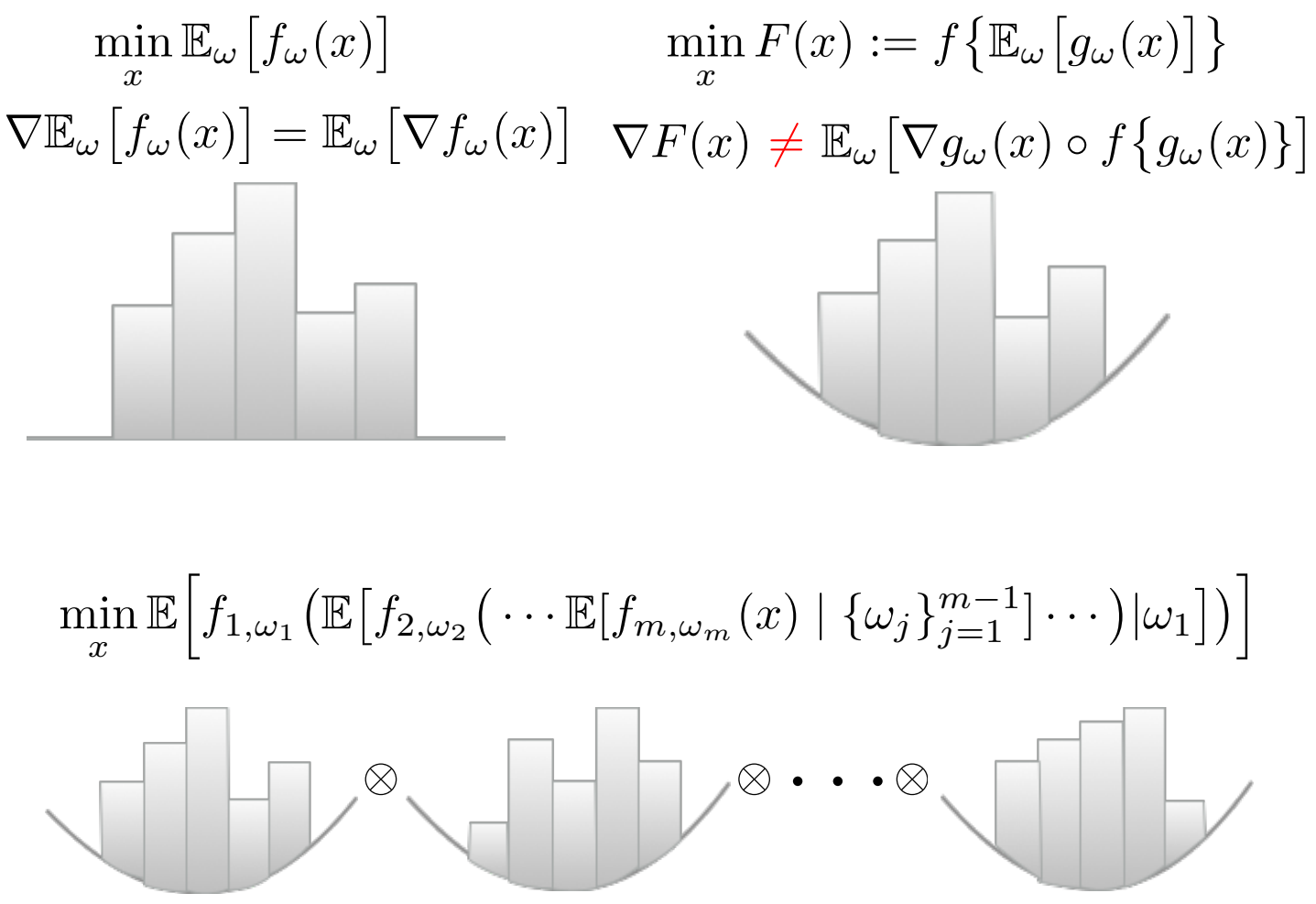

Figure 1: In one-level stochastic optimization, the objective function is linear in the probability distribution of $\omega$. In multi-level stochastic compositional optimization, the objective is no longer linear in the joint probability distribution of the random variables $(\omega_1, \ldots, \omega_m)$, making the problem fundamentally harder.

In this paper, we develop sampling-based algorithms and complexity theory for the $T$-level stochastic compositional problem (1.1). We draw motivation from the optimality conditions of problem (1.1). In particular, we expand the first-order condition into a system of variational equalities and inequalities by introducing auxiliary variables that correspond to a sequence of *value functions* at the optimal solution, i.e., tail compositions of the component functions. Our first attempt is a basic multi-timescale stochastic approximation iteration to solve this system. We establish its almost sure convergence using a $T$-element super-martingale argument for both convex and convex problems. We also show that it converges to the optimal solution at a rate of $\mathcal{O}(n^{-1/2^T})$ where $n$ is the number of iterations/oracle queries. This result suggests that the sample complexity for obtaining an approximate-optimal solution depends exponentially on the number of nested levels $T$. Such an exponential dependence is somewhat expected. It is consistent with the sample path complexity for solving multi-stage stochastic programming, although the optimization formulations and assumptions are slightly different.

Furthermore, we develop accelerated multi-level stochastic gradient methods. The accelerated algorithms apply to "smooth" composition problems and takes advantages of the smoothness of individual component functions $f^{(j)}$. An *extrapolation-interpolation* scheme is used to balance the bias-variance tradeoff in approximating each value function. The accelerated updates for the auxiliary variables can be viewed as first-oder running approximations of the true values, while the basic method without acceleration uses zeroth-order running approximations. As a result, the accelerated updates are more accurate and thus the overall convergence rate is improved. In the case when all component functions are smooth, we improve the convergence rate to $\mathcal{O}(n^{-4/(7+T)})$ for convex objective functions and $\mathcal{O}(n^{-4/(3+T)})$ for strongly convex ones. We have also obtained convergence and rate of convergence results for nonconvex problems. Table 1 summarizes our results and compare them with the best known ones for the single- and two-level stochastic compositional optimization problems [9, 19, 23, 26, 28]. We also provide numerical experiments with a risk-averse regression problem. The numerical results validate our theory.

| | | Non-convex | Convex | Strongly convex |
|---|---|---|---|---|
| 1-level | | $\mathcal{O}(n^{-1/2})$ [9] | $\mathcal{O}(n^{-1/2})$ [23] | $\mathcal{O}(n^{-1})$ [19] |
| 2-level | Smooth | $\mathcal{O}(n^{-4/9})$ [28] | $\mathcal{O}(n^{-4/9})$ [28] | $\mathcal{O}(n^{-4/5})$ [28] |
| | Non-smooth | $\mathcal{O}(n^{-1/4})$ [26] | $\mathcal{O}(n^{-1/4})$ [26] | NA |
| 3-level | Smooth | $\mathcal{O}(n^{-2/5})$ [*] | $\mathcal{O}(n^{-2/5})$ [*] | $\mathcal{O}(n^{-2/3})$ [*] |
| $T$-level | Smooth | $\mathcal{O}(n^{-4/(7+T)})$ [*] | $\mathcal{O}(n^{-4/(7+T)})$ [*] | $\mathcal{O}(n^{-4/(3+T)})$ [*] |

Table 1: Best-known $n$-sample error bound for solving multi-level stochastic compositional optimization. These bounds are achieved by stochastic gradient-type methods, so they are $n$-iteration error bounds at the same time. Note that we say the composition problem is "smooth" if *all* the component functions have Lipschitz continuous gradients. We use [*] to denote the current paper.

To the best of our knowledge, this paper proposes for the first time the multi-level stochastic gradient methods for the composition optimization problem (1.1), where we establish almost sure convergence results and obtain fast convergence rates. For the case where $T = 1$, our results match the best known sample complexity upper- and lower-bounds. For the case where $T = 2$, our results improve the convergence rate from $\mathcal{O}(n^{-2/9})$ of the a-SCGD in [26] to $\mathcal{O}(n^{-2/5})$. Besides, with additional assumption that the inner level function $f^{(T)}$ in (1.1) has Lipschitz continuous gradients, we obtain a convergence rate $\mathcal{O}(n^{-4/9})$ for two-level problems, which matches the state-of-art result achieved by ASC-PG in [28]. For the case where $T \geq 3$, our results fill the open gaps and provide the first few sample complexity benchmarks.

**Paper Organization.** Section 2 gives a basic algorithm based on multi-timescale stochastic approximation and establishes its convergence. Section 3 develops accelerated versions of the algorithm and shows that they achieve faster convergence for smooth problems. Section 4 illustrates two motivating applications in operations research and Section 5 gives numerical experiments.

**Notation and Definitions.** For $x \in \mathbb{R}^n$, we denote by $x'$ its transpose, and by $\|x\|$ its Euclidean norm (i.e.,$\|x\| = \sqrt{x'x}$). For two sequences $\{x_k\}$ and $\{y_k\}$, we write $x_k = \mathcal{O}(y_k)$ if there exists a constant $c > 0$ such that $\|x_k\| \leq c\|y_k\|$ for each $k$. We denote by $\mathbb{I}_{condition}^{value}$ the indicator function, which returns " value " if the " condition " is satisfied; otherwise 0. We denote by $F^*$ the optimal objective function value for (1.1), and denote by $\mathcal{X}^*$ the set of optimal solutions. For a set $\mathcal{X} \subset \mathbb{R}^n$ and a vector $y \in \mathbb{R}^n$, we denote by $\Pi_{\mathcal{X}}\{y\} = \mathrm{argmin}_{x\in\mathcal{X}} \|y - x\|^2$ the Euclidean projection of $y$ on $\mathcal{X}$, where the minimization is always uniquely attained if $\mathcal{X}$ is nonempty, convex and closed. For a function $f(x)$, we denote by $\nabla f(x)$ its gradient at $x$ if $f$ is differentiable, denote by $\partial f(x)$ its subdifferential at $x$, and denote by $\widetilde{\nabla} f(x)$ some noisy estimate of the gradient/subgradient of $f$ at $x$. We denote by " $\to$ " as "converges to", and denote by "w.p.1" as "with probability 1".

# 2 A Basic Algorithm Based On Multi-Timescale Stochastic Approximation

We start by writing down the optimality condition of problem (1.1) (assuming that the problem is convex):

$$\nabla F(x^*)'(x - x^*) \geq 0, \qquad \forall x \in \mathcal{X},$$

where

$$\nabla F(x) = \nabla f^{(T)}(x) \cdot \nabla f^{(T-1)}\big(f^T(x)\big) \cdots \nabla f^{(1)}\left(f^{(2)} \circ \cdots \circ f^{(T)}(x)\right).$$

However, this optimality condition is not easy to work with. As we have discussed in Section 1, the chain rule makes obtaining unbiased samples of $\nabla F(x)$ difficult. Let us rewrite the the optimality condition as follows

$$\begin{aligned}
&\left(\nabla f^{(T)}(x)\nabla f^{(T-1)}\big(y^{T-1}\big)\big) \cdots \nabla f^{(1)}\left(y^{(1)}\right)\right)'(x - x^*) \geq 0, \qquad \forall x \in \mathcal{X},\\
&y^{(T-1)} = f^{(T)}(x)\\
&y^{(T-2)} = f^{(T-1)}(y^{(T-2)}) = f^{(T-1)} \circ f^{(T)}(x)\\
&y^{(1)} = f^{(2)}(y^{(2)}) = f^{(2)} \circ \cdots \circ f^{(T)}(x).
\end{aligned}$$

We refer to $f^{(j)} \circ \cdots \circ f^{(T)}(x)$, $j = 1, \ldots, T-1$ as the *value functions*, i.e., tail compositions of multi-level component functions. By introducing the auxiliary variables $y^{(j)}$'s to represent the value functions, we can decouple the chain product. Now for a given $(x, y^{(1)}, \ldots, y^{(T-1)})$, our sampling oracle allows us to get unbiased estimates for all the quantities in the preceding system of optimality conditions.

## 2.1 A $T$-Level Stochastic Gradient Method

Motivated by the system of optimality conditions, we develop our first algorithm - a multi-timescale approximation iteration. It is also a generalization of the basic-SCGD in [26] which applies only to two-level problems. Our algorithm runs iteratively. Denote by $k$ the iteration counter. A key ingredient of our algorithm is to introduce auxiliary variables $y_k^{(j)}$'s, defined recursively, as running estimates for the value functions $\mathbb{E}_{\omega_{j,k}}[f^{(j)}_{\omega_{j,k}}(y_k^{(j+1)})|\omega_{1,k}, \cdots, \omega_{j-1,k}]$, where $j = 1, \cdots, T-1$, and $x_k = y_k^{(T)}$. At the $k$-th iteration, we update the current solution $x_k$ by using a quasi-stochastic gradient step given by

$$x_{k+1} = \Pi_{\mathcal{X}}\left\{x_k - \alpha_k \widetilde{\nabla} f^{(T)}_{\omega_{T,k}}(x_k)\nabla f^{(T-1)}_{\omega_{T-1,k}}(y_k^{(T-1)}) \cdots \nabla f^{(1)}_{\omega_{1,k}}(y_k^{(1)})\right\}.$$

**Algorithm 1** Basic Stochastic Compositional Gradient Descent (T-SCGD)

**Input :** $x_0 \in \mathbb{R}^{d_T}$, $y_0^{(j)} \in \mathbb{R}^{d_j}$, for $j = T-1, \cdots, 1$, $\mathcal{SO}$, $K$, stepsizes $\{\alpha_k\}_{k=0}^K$, $\{\beta_{j,k}\}_{k=0}^K$ for $j = T-1, ..., 1$.

**Output :** The sequence $\{x_k\}_{k=0}^{\infty}$.

**for** $k = 0, ...,$ **do**

Query the $\mathcal{SO}$ for the sample values of $f^{(T)}, \cdots, f^{(1)}$ at $(x_k, y_k^{(T-1)}, \cdots, y_k^{(1)})$, obtain the sample gradients/ subgradients $\widetilde{\nabla} f_{\omega_{T,k}}^{(T)}(x_k), \nabla f_{\omega_{T-1,k}}^{(T-1)}(y_k^{(T-1)}), \cdots, \nabla f_{\omega_{1,k}}^{(1)}(y_k^{(1)})$.

Update the main iterate by

$$x_{k+1} = \Pi_{\mathcal{X}}\Big\{x_k - \alpha_k \widetilde{\nabla} f_{\omega_{T,k}}^{(T)}(x_k) \nabla f_{\omega_{T-1,k}}^{(T-1)}(y_k^{(T-1)}) \cdots \nabla f_{\omega_{1,k}}^{(1)}(y_k^{(1)})\Big\}.$$

Query the $\mathcal{SO}$ for the sample value of $f^{(T)}(\cdot)$ at $x_k$, obtain $f_{\omega_{T,k+1}}^{(T)}(x_k)$.

Update $y_k^{(T-1)}$ by

$$y_{k+1}^{(T-1)} = (1 - \beta_{T-1,k}) y_k^{(T-1)} + \beta_{T-1,k} f_{\omega_{T,k+1}}^{(T)}(x_k).$$

**for** $j = T-2, \cdots, 1$ **do**

Query the $\mathcal{SO}$ for the sample value of $f^{(j)}$ at $y_{k+1}^{(j)}$, obtain $f_{\omega_{j,k+1}}^{(j)}(y_{k+1}^{(j)})$.

Update

$$y_{k+1}^{(j)} = (1 - \beta_{j,k}) y_k^{(j)} + \beta_{j,k} f_{\omega_{j+1,k+1}}^{(j+1)}(y_{k+1}^{(j+1)}).$$

**end for**

**end for**

Then, we update the auxiliary variables $y_k^{(j)}$'s by taking an weighted average between the previous values and the new samples returned by the $\mathcal{SO}$, i.e., for $j = T-1, T-2, \ldots, 1$,

$$y_{k+1}^{(j)} = (1 - \beta_{j,k}) y_k^{(j)} + \beta_{j,k} f_{\omega_{j+1,k+1}}^{(j+1)}(y_{k+1}^{(j+1)}), \tag{2.1}$$

where $\omega_{j,k}$ denotes the realization of $j$-th level random variable at the $k$-th iteration, $\beta_{j,k}$'s are pre-specified stepsizes. We refer to this update for $y_k^{(j)}$ as a *basic update step*. Letting $y_k^{(T)} = x_k$ and $\alpha_k = \beta_{T,k}$ to simplify the notation, we refer to the preceding iteration as the basic $T$-level Stochastic Compositional Gradient Descent ($T$-SCGD) method and summarize it in Algorithm 1. Note that we choose the stepsizes such that $\beta_{j+1,k}/\beta_{j,k} \to 0$ as $k \to \infty$ for all $j$'s, in order to control and balance the convergence speed for each auxiliary variables.

To analyze the convergence of the algorithm, we impose the following assumptions on the smoothness and bounded second-order moments for the stochastic component functions.

**Assumption 2.1.** Let $C_1, C_2, \cdots, C_T, V_1, \cdots, V_T, L_2, L_3, \cdots, L_T$ be positive scalars.

(i) The outer functions $f^{(T-1)}, f^{(T-2)}, \cdots, f^{(1)}$ are continuously differentiable, the inner function $f^{(T)}$ is continuous, the feasible set $\mathcal{X}$ is closed and convex, and there exists at least one optimal solution $x^*$ to problem (1.1).

(ii) The sample paths $(\omega_{1,0}, \omega_{2,0}, \cdots, \omega_{T,0})$, $(\omega_{1,1}, \omega_{2,1}, \cdots, \omega_{T,1})$,...,$(\omega_{1,k}, \omega_{2,k}, \cdots, \omega_{T,k})$ are independent across $k$ and satisfy with probability 1

$$\mathbb{E}[f^{(j)}_{\omega_{j,0}}(x_j)|\omega_{1,0}, \cdots, \omega_{j-1,0}] = f^{(j)}(x_j), \forall x_j \in \mathbb{R}^{d_j} \text{ for } j = 1, \cdots, T, \text{ and } \mathbb{E}[\widetilde{\nabla} F_{\omega_0}(x)] \in \partial F(x),$$

for all $x \in \mathcal{X}$, where $\widetilde{\nabla} F_{\omega_0}(x) \equiv \widetilde{\nabla} f^{(T)}_{\omega_{T,0}}(x) \nabla f^{(T-1)}_{\omega_{T-1,0}}\big(f^{(T)}(x)\big) \cdots \nabla f^{(1)}_{\omega_{1,0}}\Big(f^{(2)} \circ \cdots \circ f^{(T)}(x)\Big)$.

(iii) The function $f^{(T)}(\cdot)$ is Lipschitz continuous with parameter $C_T$, and the samples $f^{(T)}_{\omega_{T,0}}(\cdot)$, $\widetilde{\nabla} f^{(T)}_{\omega_{T,0}}(\cdot)$ have bounded second-order moments such that with probability 1

$$\mathbb{E}\big[\|\widetilde{\nabla} f^{(T)}_{\omega_{T,0}}(x)\|^2|\omega_{T-1,0}, \cdots, \omega_{1,0}\big] \le C_T, \mathbb{E}\big[\|f^{(T)}_{\omega_{T,0}}(x) - f^{(T)}(x)\|^2|\omega_{T-1,0}, \cdots, \omega_{1,0}\big] \le V_T,$$

for all $x \in \mathcal{X}$.

(iv) For $j = 1, \cdots, T-1$, the functions $f^{(j)}(\cdot)$'s and $f^{(j)}_{\omega_{j,0}}(\cdot)'s$ have $L_j$-Lipschitz continuous gradients such that with probability 1

$$\mathbb{E}\big[\|\nabla f^{(j)}_{\omega_{j,0}}(x_j)\|^2|\omega_{j-1,0}, \cdots, \omega_{1,0}\big] \le C_j, \mathbb{E}\big[\|f^{(j)}_{\omega_{j,0}}(x_j) - f^{(j)}(x_j)\|^2|\omega_{j-1,0}, \cdots, \omega_{1,0}\big] \le V_j,$$

$$\text{and } \|\nabla f^{(j)}_{\omega_{j,0}}(x_j) - \nabla f^{(j)}_{\omega_{j,0}}(\bar{x}_j)\| \le L_j\|x_j - \bar{x}_j\|,$$

for all $x_j, \bar{x}_j \in \mathbb{R}^{d_j}$.

In some part of the analysis, we also assume that the overall objective is sufficiently smooth as follows.

**Assumption 2.2.** The function $F(x)$ has Lipschitz continuous gradient, i.e.,there exists $L_F > 0$ such that

$$F(z) - F(x) \le \langle \nabla F(x), z - x\rangle + \frac{L_F}{2}\|z - x\|^2, \ \ \forall x, z.$$

Note that in Assumption 2.1, we require the functions $f^{(1)}(\cdot), \cdots, f^{(T-1)}(\cdot)$ to have Lipschitz continuous gradients, and we do not impose such assumptions on $f^{(T)}(\cdot)$. Hence, we cannot guarantee that $F(x)$ has a Lipschitz continuous gradient, which means Assumption 2.1 does not imply Assumption 2.2.

## 2.2 Almost Sure Convergence of $T$-SCGD

Theoretical analysis of Algorithm 1 is challenging due to the nested level of expections over a path of random variables. The multiple nested levels of expectations need to be carefully estimated and balanced to ensure convergence of the algorithm. We first prove the almost sure convergence of the algorithm as long as the step-sizes are properly chosen and diminishing. For convex problems, we show that the algorithm generates a sequence of solutions that converges to an optimal solution to problem (1.1) with probability 1. For nonconvex problems, we show that all limiting points of the sequence generated by this algorithm are stationary points with probability 1. In the rest of this subsection, we give a proof outline with all the lemmas. We defer the full proof to Appendix A.

**Theorem 2.1** (Almost sure convergence of $T$-SCGD)**.** Let Assumption 2.1 hold, and let the step-sizes $\{\alpha_{1,k}\}, \{\beta_{2,k}\}, \cdots, \{\beta_{T,k}\}$ be such that

$$\sum_{k=0}^{\infty} \alpha_k = \infty, \sum_{k=0}^{\infty} \beta_{j,k} = \infty, \text{ for all } j = T-1, ..., 1,$$

and

$$\sum_{k=0}^{\infty} \Big( \alpha_k^2 + \beta_{T-1,k}^2 + \cdots + \beta_{1,k}^2 + \frac{\alpha_k^2}{\beta_{2,k}} + \frac{\alpha_k^2}{\beta_{3,k}} + \cdots + \frac{\alpha_k^2}{\beta_{T-1,k}} + \frac{\beta_{T-1,k}^2}{\beta_{T-2,k}} + \cdots + \frac{\beta_{2,k}^2}{\beta_{1,k}} \Big) < \infty.$$

Let $\{(x_k, y_k^{(T-1)}, \cdots, y_k^{(1)})\}_{k=0}^{\infty}$ be the sequence generated by the $T$-SCGD Algorithm 1 starting with an arbitrary initial point $(x_0, y_0^{(T-1)}, \cdots, y_0^{(1)})$. Then:

(a) If $F$ is convex, $\{x_k\}$ converges almost surely to a random point in the set of optimal solutions to problem (1.1).

(b) Suppose in addition that Assumption 2.2 holds, $\mathcal{X} = \mathbb{R}^{d_T}$, and all samples generated by the $\mathcal{SO}$ are uniformly bounded. Then any limiting point of the sequence $\{x_k\}_{k=0}^{\infty}$ is a stationary point to problem (1.1) almost surely.

*Proof Outline.* We denote by $\mathbb{F}_k$ the collection of random variables up to the $k$-th iteration to help us better analyze the convergence properties:

$$\Big\{\{x_i\}_{i=0}^{k}, \{y_i^{(T-1)}\}_{i=0}^{k-1}, \cdots, \{y_i^{(1)}\}_{i=0}^{k-1}, \{\omega_{T,i}\}_{i=0}^{k-1}, \cdots, \{\omega_{1,i}\}_{i=0}^{k-1}\Big\}.$$

To derive the almost sure convergence of Algorithm 1, we construct two different $T$-element super-martingales for the convex and non-convex objectives, respectively.

Firstly, for problems with convex objective $F$, in the $k$-th iteration, we have the following lemma to analyze the improvement from $\|x_k - x^*\|$ to $\|x_{k+1} - x^*\|$ by $\|y_k^{(T-1)} - f^{(T)}(x_k)\|$, $\|y_k^{(T-2)} - f^{(T-1)}(y_k^{(T-1)})\|$, $\cdots$, and $\|y_k^{(1)} - f^{(2)}(y_k^{(2)})\|$.

**Lemma 2.1.** Let Assumption 2.1 hold, and let $F = f^{(1)} \circ f^{(2)} \circ \cdots \circ f^{(T)}$ be convex. Then Algorithm 1 generates a sequence $\{(x_k, y_k^{(T-1)}, \cdots, y_k^{(1)})\}_{k=0}^{\infty}$ such that there exists a constant $C_0 > 0$ and an optimal solution $x^* \in \mathcal{X}^*$, for all $k$, with probability 1,

$$\begin{aligned}
&\mathbb{E}[\|x_{k+1} - x^*\|^2 | \mathbb{F}_k] \\
\leq &\Big(1 + [\frac{\alpha_k^2}{\beta_{T-1,k}} + \cdots + \frac{\alpha_k^2}{\beta_{1,k}}] C_0\Big) \|x_k - x^*\|^2 + \alpha_k^2 C_1 C_2 \cdots C_T - 2\alpha_k \big(F(x_k) - F^*\big) \\
&+ (T-1)\beta_{T-1,k} \mathbb{E}[\|y_k^{(T-1)} - f^{(T)}(x_k)\|^2 | \mathbb{F}_k] + (T-2)\beta_{T-2,k} \mathbb{E}[\|y_k^{(T-2)} - f^{(T-1)}(y_k^{(T-1)})\|^2 | \mathbb{F}_k] \\
&+ \cdots + \beta_{1,k} \mathbb{E}[\|y_k^{(1)} - f^{(2)}(y_k^{(2)})\|^2 | \mathbb{F}_k].
\end{aligned} \tag{2.2}$$

Lemma 2.1 states that for $T$-level SCGD with convex objective function $F$, the optimality error $\|x_{k+1} - x^*\|$ can be bounded by $\|x_k - x^*\|$, $\|y_k^{(T-1)} - f^{(T)}(x_k)\|$, $\|y_k^{(T-2)} - f^{(T-1)}(y_k^{(T-1)})\|$, $\cdots$, and $\|y_k^{(1)} - f^{(2)}(y_k^{(2)})\|$ in a super-martingale form.

Next, we present a lemma used in the analysis in part (b).

**Lemma 2.2.** Suppose that Assumption 2.1 and 2.2 hold, and $\mathcal{X} = \mathbb{R}^{d_T}$. Let $F^* = \min_{x\in\mathcal{X}} F(x)$, then Algorithm 1 generates a sequence $\{(x_k, y_k^{(T-1)}, \cdots, y_k^{(1)})\}_{k=0}^{\infty}$ such that

$$\begin{aligned}
&\mathbb{E}[F(x_{k+1}) - F^*|\mathbb{F}_k]\\
&\leq F(x_k) - F^* - \frac{\alpha_k}{2}\|\nabla F(x_k)\|^2 + \frac{1}{2}\alpha_k^2 L_F C_1 C_2 \cdots C_T + (T-1)\beta_{T-1,k}\mathbb{E}[\|y_k^{(T-1)} - f^{(T)}(x_k)\|^2|\mathbb{F}_k]\\
&\quad + (T-2)\beta_{T-2,k}\mathbb{E}[\|y_k^{(T-2)} - f^{(T-1)}(y_k^{(T-1)})\|^2|\mathbb{F}_k] + \cdots + \beta_{1,k}\mathbb{E}[\|y_k^{(1)} - f^{(2)}(y_k^{(2)})\|^2|\mathbb{F}_k],
\end{aligned}$$

for $k$ sufficiently large, with probability 1.

This lemma tells us that for $T$-level SCGD with general nonconvex objective function $F$, $(F(x_{k+1}) - F^*)$ can be bounded by $(F(x_k) - F^*)$, $\|y_k^{(T-1)} - f^{(T)}(x_k)\|$, $\|y_k^{(T-2)} - f^{(T-1)}(y_k^{(T-1)})\|$, $\cdots$, and $\|y_k^{(1)} - f^{(2)}(y_k^{(2)})\|$ in a super-martingale form. Similar as in Lemma 2.1, we shall construct the super-martingales for $\|y_k^{(T-1)} - f^{(T)}(x_k)\|$ and $\|y_k^{(j)} - f^{(j+1)}(y_k^{(j+1)})\|$ for $j = T-2, \cdots, 1$ respectively, and then use Lemma 2.4 to show the almost sure convergence of $(F(x_k) - F^*)$ for a $T$-level SCGD with nonconvex objective $F$. With further analysis, we show that any limiting point of the sequence $\{x_k\}_{k=0}^{\infty}$ is a stationary point with probability 1, which proves part (b) of Theorem 2.1.

Next, we analyze the term $\|y_k^{(j)} - f^{(j+1)}(y_k^{(j+1)})\|$ for $j = T-1, \cdots, 1$ and construct the proper super-martingales for them respectively.

**Lemma 2.3.** Let Assumption 2.1 hold, and let $\{(x_k, y_k^{(T-1)}, \cdots, y_k^{(1)})\}_{k=0}^{\infty}$ be the sequence generated by Algorithm 1. For $j = T-1, \cdots, 1$, suppose $\mathbb{E}[\|y_{k+1}^{(j+1)} - y_k^{(j+1)}\|^2] \leq \mathcal{O}(\beta_{j+1,k}^2)$ for all $k$, then we have

(a) For all $k$, with probability 1,

$$\begin{aligned}
&\mathbb{E}[\|y_{k+1}^{(j)} - f^{(j+1)}(y_{k+1}^{(j+1)})\|^2|\mathbb{F}_{k+1}]\\
\leq&(1-\beta_{j,k})\|y_k^{(j)} - f^{(j+1)}(y_k^{(j+1)})\|^2 + \beta_j^{-1}C_{j+1}\mathbb{E}[\|y_{k+1}^{(j+1)} - y_k^{(j+1)}\|^2|\mathbb{F}_{k+1}] + 2V_{j+1}\beta_{j,k}^2.
\end{aligned} \tag{2.3}$$

(b) If $\sum_{k=1}^{\infty} \beta_{j+1,k}^2/\beta_{j,k} < \infty$, then

$$\sum_{k=1}^{\infty} \beta_{j,k}^{-1}\mathbb{E}[\|y_{k+1}^{(j)} - y_k^{(j)}\|^2|\mathbb{F}_{k+1}] < \infty, \quad w.p.1.$$

(c) There exists a constant $D_j \geq 0$ such that $\mathbb{E}[\|y_{k+1}^{(j)} - f^{(j+1)}(y_{k+1}^{(j+1)})\|^2] \leq D_j$ for all $k$.

(d) $\mathbb{E}[\|y_{k+1}^{(j)} - y_k^{(j)}\|^2] \leq \mathcal{O}(\beta_{j,k}^2)$ for all $k$.

Note that here we use $y_k^{(T)} = x_k$ and $\beta_{T,k} = \alpha_k$ for ease of notation. This lemma constructs super-martingales of $\{\|y_k^{(j)} - f^{(j+1)}(y_k^{(j+1)})\|\}_{k=1}^{\infty}$ for $j = T-1\cdots,1$ respectively, and it also shows that under proper assumptions, the tail part for the super-martingale, $\beta_j^{-1}C_{j+1}\mathbb{E}[\|y_{k+1}^{(j)} - y_k^{(j)}\|^2|\mathbb{F}_k] + 2V_{j+1}\beta_{j,k}^2$, converges almost surely.

Previous lemmas provide basic blocks for us to build a $T$-element super-martingale. We then provide the $T$-element super-martingale convergent lemma to establish the convergence property of $\{x_k - x^*\}$.

**Lemma 2.4** ($T$-element supermartingale convergence)**.** Let $\{X_k\}$,$\{Y_k^{(T-1)}\}$,$\cdots$,$\{Y_k^{(1)}\}$, $\{\eta_k\}$, and $\{u_k^{(j)}\}$,$\{\mu_k^{(j)}\}$,$\{\theta_k^{(j)}\}$ for $j = 1, \cdots, T$ be sequences of nonnegative random variables such that

$$\mathbb{E}[X_{k+1}|\mathbb{G}_k] \leq (1+\eta_k)X_k - u_k^{(T)} + \sum_{j=1}^{T-1} c_j \theta_k^{(j)} Y_k^{(j)} + \mu_k^{(T)},$$

and

$$\mathbb{E}[Y_{k+1}^{(T-1)}|\mathbb{G}_k] \leq (1-\theta_k^{(j)})Y_k^{(j)} - u_k^{(j)} + \mu_k^{(j)}, \text{for } j = T-1, ..., 1,$$

for all $k$, where $\mathbb{G}_k$ is the collection of random variables

$$\Big\{\{X_i\}_{i=0}^k, \{Y_i^{(T-1)}\}_{i=0}^k, \cdots, \{Y_i^{(1)}\}_{i=0}^k, \{\eta_i\}_{i=0}^k, \{u_i^{(j)}\}_{i=0}^k, \{\mu_i^{(j)}\}_{i=0}^k, \{\theta_i^{(j)}\}_{i=0}^k, \text{ for } j = 1, \cdots, T\Big\},$$

and $c_{T-1}, c_{T-2} \cdots, c_1$ are positive scalars. Assume that

$$\sum_{k=0}^{\infty} \eta_k < \infty, \sum_{k=0}^{\infty} \mu_k^{(j)} < \infty, \quad \text{for } j = 1, \cdots, T.$$

Then $\{X_k\}$, $\{Y_k^{(1)}\}$, $\{Y_k^{(2)}\}$, $\cdots$, $\{Y_k^{(T-1)}\}$ converge almost surely to $T$ nonnegative random variables respectively, and we have

$$\sum_{j=1}^{T}\sum_{k=0}^{\infty} u_k^{(j)} < \infty, \sum_{k=0}^{\infty}\sum_{j=1}^{T-1} c_j \theta_k^{(j)} Y_k^{(j)} < \infty \;\; w.p.1.$$

To prove Theorem 2.1 part (a), by Lemma 2.1 and Lemma 2.3, we construct a $T$-element super-martingale by letting

$$\begin{aligned}
&X_k = \|x_k - x^*\|^2, Y_k^{(T-1)} = \mathbb{E}[\|y_k^{(T-1)} - f^{(T)}(x_k)\|^2|\mathbb{F}_k],\\
&Y_k^{(T-2)} = \mathbb{E}[\|y_k^{(T-2)} - f^{(T-1)}(y_k^{(T-1)})\|^2|\mathbb{F}_k], \cdots, Y_k^{(1)} = \mathbb{E}[\|y_k^{(1)} - f^{(2)}(y_k^{(2)})\|^2|\mathbb{F}_k],\\
&\eta_k = [\frac{\alpha_k^2}{\beta_{T-1,k}} + \cdots + \frac{\alpha_k^2}{\beta_{1,k}}]C_0, u_k^{(T)} = 2\alpha_k(F(x_k) - F^*),\\
&u_k^{(1)} = u_k^{(2)} = \cdots = u_k^{(T-1)} = 0, c_1 = 1, \cdots, c_{T-1} = T-1,\\
&\mu_k^{(1)} = 2\beta_{1,k}^2 V_1 + \mathcal{O}(\frac{\mathbb{E}[\|y_{k+1}^{(2)} - y_k^{(2)}\|^2|\mathbb{F}_k]}{\beta_{1,k}}), \cdots,\\
&\mu_k^{(T-2)} = 2\beta_{T-2,k}^2 V_{T-1} + \mathcal{O}(\frac{\mathbb{E}[\|y_{k+1}^{(T-1)} - y_k^{(T-1)}\|^2|\mathbb{F}_k]}{\beta_{T-2,k}}),\\
&\mu_k^{(T-1)} = 2\beta_{T-1,k}^2 V_T + \mathcal{O}(\frac{\mathbb{E}[\|x_{k+1} - x_k\|^2|\mathbb{F}_k]}{\beta_{T-1,k}}),\\
&\mu_k^{(T)} = \alpha_k^2 C_1 C_2 \cdots C_T, \theta_j^{(1)} = \beta_{1,k}, \cdots, \theta_j^{(T-1)} = \beta_{T-1,k}.
\end{aligned}$$

Under the conditions in Theorem 2.1, we have that the $T$-element super-martingale converges almost surely to $T$ random variables by Lemma 2.4, thus $\|x_k - x^*\|$ converges almost surely, and

$$\sum_{k=0}^{\infty} \alpha_k(F(x_k) - F^*) < \infty, \qquad w.p.1,$$

which further implies that

$$\liminf_{k\to\infty} F(x_k) = F^*, \qquad w.p.1.$$

Finally, the following lemma shows the the sequence $\{x_k\}_{k=0}^{\infty}$ converges almost surely to an optimal solution to problem (1.1), which completes the proof of part (a).

**Lemma 2.5.** Let $\{(x_k, y_k^{(T-1)}, \cdots, y_k^{(1)})\}_{k=0}^{\infty}$ be the sequence generated by Algorithm 1. Let $F^* = F(x^*)$, where $x^*$ is an optimal solution to problem (1.1). Suppose

$$\liminf_{k\to\infty} F(x_k) = F^*, \qquad w.p.1,$$

then $\{x_k\}$ converges almost surely to a random point in the set of optimal solutions to problem (1.1).

For part (b), by Lemma 2.2 and Lemma 2.3, we construct the $T$-element super-martingale for general non-convex functions, and show that $\{F(x_k) - F^*\}$ converges almost surely by Lemma 2.4, which further implies $\sum_{k=0}^{\infty} \alpha_k \|\nabla F(x_k)\|^2 < \infty$ with probability 1. Then, we have the following lemma which shows that any limiting point of the sequence $\{x_k\}$ is a stationary point of $F(x)$ with probability 1.

**Lemma 2.6.** Let $\big\{(x_k, y_k^{(T-1)}, \cdots, y_k^{(1)})\big\}_{k=0}^{\infty}$ be the sequence generated by Algorithm 1. Suppose $\sum_{k=0}^{\infty} \alpha_k = \infty$ and $\sum_{k=0}^{\infty} \alpha_k \|\nabla F(x_k)\|^2 < \infty$ with probability 1, then any limiting point of the sequence $\{x_k\}$ is a stationary point of $F(x)$ with probability 1.

This concludes the proof for part (b). □

## 2.3 Convergence Rate of $T$-SCGD

In this subsection, we analyze the convergence rate of Algorithm 1. Specifically, we derive the rate through taking the averaged iterates $\widehat{x}_n = \frac{1}{N_n} \sum_{k=n-N_n+1}^{n} x_k$, where $N_n = \lceil n/2 \rceil$. Note that similar results still hold if we let $N_n = n/C$ for other constant $C > 0$.

Clearly, the rate of convergence is closely related to the stepsizes $\alpha_k$'s and $\beta_{j,k}$'s. We consider stepsizes of the form

$$\alpha_k = k^{-a} \text{ and } \beta_{j,k} = k^{-b_j} \text{ for all } j = T-1, ..., 1, \tag{2.4}$$

where $a$ and $b_j$'s are real numbers.

**Theorem 2.2** (Convergence rate of $T$-SCGD)**.** Suppose that Assumption 2.1 holds, and the objective function $F(\cdot)$ is convex. Let $D > 0$ be such that $\sup_{x_0 \in \mathcal{X}} \|x_0 - x^*\| < D$, and let the stepsizes be $\alpha_k = k^{-a}$, $\beta_{j,k} = k^{-b_j}$ for $j = T-1, \cdots, 1$, where $(a, b_{T-1}, b_{T-2}, \cdots, b_1) \in (0, 1)$. If we choose $a = 1 - \frac{1}{2^T}, b_{T-1} = 1 - \frac{1}{2^{T-1}}, \cdots, b_1 = 1 - \frac{1}{2}$, letting $\{(x_k, y_k^{(T-1)}, \cdots, y_k^{(1)})\}_{k=0}^{\infty}$ be the sequence generated by the $T$-SCGD Algorithm 1 starting with an arbitrary initial point $(x_0, y_0^{(T-1)}, \cdots, y_0^{(1)})$, we obtain

$$\mathbb{E}[F(\widehat{x}_n) - F^*] \le \mathcal{O}(n^{-1/2^T}).$$

*Proof Outline.* We present the outline of proof here and defer the details in Appendix B. We first derive the convergence rate of $\|y_{k+1}^{(j)} - f^{(j+1)}(y_{k+1}^{(j+1)})\|$ and $\|y_{k+1}^{(j)} - y_k^{(j)}\|$ for $j = T-1, \cdots, 1$. By Lemma 2.3 and Lemma B.1 in the Appendix, we have the following lemma characterizing the corresponding convergence rates:

**Lemma 2.7.** Let Assumption 2.1 hold, and let $\{(x_k, y_k^{(T-1)}, \cdots, y_k^{(1)})\}_{k=0}^{\infty}$ be the sequence generated by Algorithm 1. For any basic update step $j = T-1, \cdots, 1$, we have

$$\mathbb{E}[\|y_k^{(j)} - f^{(j+1)}(y_k^{(j+1)})\|^2] \le \mathcal{O}(k^{-2b_{j+1}+2b_j}) + \mathcal{O}(k^{-b_j}) \quad \text{for all } k.$$

Next, define the random variable

$$J_k \equiv \|x_k - x^*\|^2 + (T-1)\|y_k^{(T-1)} - f^{(T)}(x_k)\|^2 + (T-2)\|y_k^{(T-2)} + f^{(T-1)}(y_k^{(T-1)})\|^2 + \cdots + \|y_k^{(1)} - f^{(2)}(y_k^{(2)})\|^2,$$

so we have $\mathbb{E}[J_k] \le D_T + (T-1)D_{T-1} + (T-2)D_{T-2} + \cdots D_1 \equiv D_J$.

We multiply (2.3) by $j \times (1 + \beta_{j,k})$ for every $j$ from $T-1$ to 1, take their sums with (2.2), and take the expectation on both sides. By Lemma 2.7, we obtain

$$\begin{aligned}\mathbb{E}[J_{k+1}] \le & \Big(1 + [\frac{\alpha_k^2}{\beta_{T-1,k}} + \frac{\alpha_k^2}{\beta_{T-2,k}} + \cdots + \frac{\alpha_k^2}{\beta_{1,k}}]C_0\Big)\mathbb{E}[J_k] - 2\alpha_k\big(F(x_k - F^*)\big) \\ & + C_1C_2\cdots C_T\alpha_k^2 + 4(T-1)V_T\beta_{T-1,k}^2 + \frac{2(T-1)C_1C_2\cdots C_{T-1}C_T^2\alpha_k^2}{\beta_{T-1,k}} \\ & + \sum_{j=1}^{T-2} j \times \Big[4V_{j+1}\beta_{j,k}^2 + \frac{2C_{j+1}}{\beta_{j,k}}\mathcal{O}(\beta_{j+1,k}^2)\Big].\end{aligned}$$

Let $N > 0$, by reordering the terms in the preceding relation and taking its sum over $k-N, \cdots, k$, with basic algebra, we have

$$\begin{aligned}& 2\sum_{t=k-N}^{k} \mathbb{E}[F(x_t) - F^*] \\ \le\ & k^a D + \sum_{t=k-N}^{k} (k^{-a+b_{T-1}} + \cdots + k^{-a+b_1})C_0 D_J \\ & + \sum_{t=k-N}^{k} \mathcal{O}(k^{-a}) + \sum_{j=1}^{T-1}\sum_{t=k-N}^{k} \mathcal{O}(k^{-2b_j+a}) + \sum_{t=k-N}^{k} \mathcal{O}(k^{-a+b_{T-1}}) + \sum_{j=1}^{T-2}\sum_{t=k-N}^{k} \mathcal{O}(k^{-2b_{j+1}+a+b_j}).\end{aligned}$$

Finally, we optimize the convergence rate by choosing $a = 1 - \frac{1}{2^T}, b_{T-1} = 1 - \frac{1}{2^{T-1}}, \cdots, b_1 = 1 - \frac{1}{2}$ and obatin

$$\mathbb{E}[F(\widehat{x}_k) - F^*] \le \frac{1}{N_k}\sum_{t=k-N}^{k} \mathbb{E}[F(x_t) - F^*] \le \mathcal{O}(k^{-1/2^T}),$$

which completes the proof. □

This result provides a sample complexity upper bound for the multi-level stochastic compositional optimization problem. In the case where $T = 2$, this result matches the convergence of basic-SCGD given in [26].

# 3 Accelerated Multi-Level Stochastic Gradient Algorithm

In the previous section, we establish an $\mathcal{O}(n^{-1/2^T})$ rate of convergence for the $T$-level stochastic compositional optimization problem. A key question is whether and when we can better utilize noisy gradients of component functions and improve the overall convergence rate.

Throughout this section, in addition to Assumption 2.1, we impose the following assumption:

**Assumption 3.1.** Let $C_1, C_2, \cdots, C_T, V_1, \cdots, V_T$ be positive scalars.

(i) The samples $f^{(j)}_{\omega_{T,k}}(\cdot), \widetilde{\nabla} f^{(j)}_{\omega_{T,k}}(\cdot)$ have bounded fourth-order moments such that with probability 1,

$$\mathbb{E}[\|\widetilde{\nabla} f^{(T)}_{\omega_{T,0}}(x)\|^4 | \omega_{1,0}, \cdots, \omega_{T-1,0}] \le C_T^2,$$
$$\text{and } \mathbb{E}[\|f^{(T)}_{\omega_{T,0}}(x) - f^{(T)}(x)\|^4 | \omega_{1,0}, \cdots, \omega_{T-1,0}] \le V_T^2, \ \ \forall x \in \mathcal{X}.$$

(ii) The samples $f^{(j)}_{\omega_{j,k}}(\cdot)$'s and $\nabla f^{(j)}_{\omega_{j,k}}(\cdot)$'s have bounded fourth-order moments such that with probability 1,

$$\mathbb{E}[\|\nabla f^{(j)}_{\omega_{j,0}}(x_j)\|^4 | \omega_{1,0}, \cdots, \omega_{j-1,0}] \le C_j^2,$$
$$\text{and } \mathbb{E}[\|f^{(j)}_{\omega_{j,0}}(x_j) - f^{(j)}(x_j)\|^4 | \omega_{1,0}, \cdots, \omega_{j-1,0}] \le V_j^2, \ \ \forall x_j \in \mathbb{R}^{d_j}, \quad \text{and for } j = T-1, \cdots, 1.$$

We also consider the case when the first inner level function $f^{(T)}$ also has Lipschitz continuous gradients. In some part of our subsequent analysis, we make the following assumption.

**Assumption 3.2.** The function $f^{(T)}$ has Lipschitz continuous gradient such that

$$\|\nabla f^{(T)}(x) - \nabla f^{(T)}(\bar{x})\| \le L_T \|x - \bar{x}\|,$$

for all $x, \bar{x} \in \mathcal{X}$ .

In what follows, we propose an accelerated algorithm to better utilize those smoothness properties and achieve improved convergence rates. I

## 3.1 An Extrapolation-Interpolation Scheme For Acceleration

The basic idea of acceleration is to refine the running estimates of the value functions by using additional extrapolations. The same idea has been used for the case where $T = 2$. Specifically, in [26], with an additional bounded fourth moments assumption, the authors developed an accelerated SCGD (a-SCGD) algorithm and achieved faster convergence rate using an extra extrapolation step per iteration.

Now we develop a new accelerated algorithm for the multi-level problem that runs as follows: At the $k$-th iteration, we first update the main iterate solution $x_{k+1}$ by the chain rule,

$$x_{k+1} = \Pi_{\mathcal{X}} \Big\{ x_k - \alpha_k \widetilde{\nabla} f^{(T)}_{\omega_{T,k}}(x_k) \nabla f^{(T-1)}_{\omega_{T-1,k}}(y_k^{(T-1)}) \cdots \nabla f^{(1)}_{\omega_{1,k}}(y_k^{(1)}) \Big\}.$$

We then update the running estimate $y_k^{(T-1)}$ for $\mathbb{E}_{\omega_{T,k}}[f_{\omega_{T,k}}^{(T)}(x_k)|\omega_{1,k},\cdots,\omega_{T-1,k}]$ by taking weighted average between the new sample and the previous estimate. Specifically, we update $y_k^{(T-1)}$ by letting

$$y_{k+1}^{(T-1)} = (1-\beta_{T-1,k})y_k^{(T-1)} + \beta_{T-1,k} f_{\omega_{T,k+1}}^{(T)}(x_{k+1}).$$

Next, we conduct extrapolation steps for acceleration. The intuition is that we can use sample gradients of individual component functions more efficiently when these functions are smooth, which allows us to obtain better estimates of $f^{(j)}$'s. In particular, our accelerated updates for the auxiliary variables are performing first-oder running approximations of the true values. In comparison, the corresponding updates used in $T$-SCGD can be viewed as zeroth-order running approximations. Specifically, at the $k$-th iteration, we refine our estimate $y_{k+1}^{(j)}$ by taking an additional extrapolation step and obtaining a new auxiliary variable $\widehat{y}_{k+1}^{(j)}$:

$$\widehat{y}_{k+1}^{(j)} = (1-1/\beta_{j,k})y_k^{(j+1)} + y_{k+1}^{(j+1)}/\beta_{j,k}.$$

Then, when we update $y_{k+1}^{(j)}$, we plug in this auxiliary variable aiming for a better estimate that

$$y_{k+1}^{(j)} = (1-\beta_{j,k})y_k^{(j)} + \beta_{j,k}\cdot f_{\omega_{j+1,k+1}}^{(j+1)}(\widehat{y}_{k+1}^{(j)}).$$

We point out that this is essentially a weighted smoothing scheme, where $\widehat{y}_k^{(j)}$'s are obtained through extrapolation steps to further utilize the smoothness in order to improve the convergence rate. Roughly speaking, this further extrapolation step helps us achieve estimators $y_{k+1}^{(j)}$'s for $f^{(j+1)}(y_{k+1}^{(j+1)})$'s accurate up to the second order terms if we take Taylor expansions of $f^{(j)}$'s. In comparison, without the extrapolation, if we directly plug in $y_{k+1}^{(j+1)}$'s instead, the estimators are only accurate up to the first order terms. We call this an *accelerating update step*. Note that here we do not assume $f^{(T)}$ has Lipschitz continuous gradient as in some applications, $f^{(T)}$ includes some sparse-inducing regularization terms and is not continuously differentiable.

When Assumption 3.2 holds, we update the main iteration by the chain rule, and then apply extrapolation to this level to better utilize the smoothness. That is, we refine our estimate $y_{k+1}^{(T-1)}$ with an additional extrapolation step and an auxiliary variable $\widehat{y}_{k+1}^{(T-1)}$ as

$$\widehat{y}_{k+1}^{(T-1)} = (1-1/\beta_{T-1,k})x_k + x_{k+1}/\beta_{T-1,k}.$$

Next, we update $y_{k+1}^{(T-1)}$ by this auxiliary variable such that

$$y_{k+1}^{(T-1)} = (1-\beta_{T-1,k})y_k^{(T-1)} + \beta_{T-1,k} f_{\omega_{T,k+1}}^{(T)}(\widehat{y}_{k+1}^{(T-1)}).$$

For the remaining levels, we apply the same procedure as in the accelerating update steps previously described. We summarize those two slightly different accelerated algorithms in Algorithm 2.

In the remaining part of this section, we provide theoretical guarantees for this accelerated algorithm. We first provide the almost sure convergence result that almost surely, our algorithm converges to an optimal solution when the problem is convex, and any limiting point of the generated solution path is a stationary point. Next, we obtain an improved convergence rate for our algorithm for general nonconvex objective functions. Furthermore, we investigate the case when the objective function is strongly convex, and shows that one can achieve faster convergence. For all results, we provide outlines and key lemmas in the main text, and defer the detailed proofs in Appendix C, D and E.

**Algorithm 2** Accelerated $T$-Level Stochastic Compositional Gradient Descent ($a$-TSCGD)

**Input:** $x_0 \in \mathbb{R}^{d_T}, y_0^{(j)} \in \mathbb{R}^{d_j}$ for $j = T-1, ..., 1$, $\mathcal{SO}$, $K$, stepsizes $\{\alpha_k\}_{k=0}^K$,$\{\beta_{j,k}\}_{k=0}^K$ for $j = T-1, ..., 1$.

**Output:** The sequence $\{x_k\}_{k=0}^K$.

**for** $k = 0, ..., K$ **do**

Query the $\mathcal{SO}$ for the sample values of $f^{(T)}, \cdots, f^{(1)}$ at $x_k, y_k^{(T-1)}, \cdots, y_k^{(1)}$, obtain $\widetilde{\nabla} f_{\omega_{T,k}}^{(T)}(x_k)$, $\nabla f_{\omega_{T-1,k+1}}^{(T-1)}(y_k^{(T-1)}), \cdots, \nabla f_{\omega_{1,k}}^{(1)}(y_k^{(1)})$.

Update the main iterate by

$$x_{k+1} = \Pi_{\mathcal{X}}\Big\{x_k - \alpha_k \widetilde{\nabla} f_{\omega_{T,k}}^{(T)}(x_k) \nabla f_{\omega_{T-1,k}}^{(T-1)}(y_k^{(T-1)}) \cdots \nabla f_{\omega_{1,k}}^{(1)}(y_k^{(1)})\Big\}.$$

**if** Assumption 3.2 is known to hold **then**

Update the auxiliary variable $\widehat{y}_{k+1}^{(T-1)}$ by

$$\widehat{y}_{k+1}^{(T-1)} = (1 - 1/\beta_{T-1,k})x_k + x_{k+1}/\beta_{T-1,k}.$$

Query the $\mathcal{SO}$ for the sample value of $f^{(T)}$ at $\widehat{y}_{k+1}^{(T-1)}$, obtain $f_{\omega_{T,k+1}}^{(T)}(\widehat{y}_{k+1}^{(T-1)})$.

Update

$$y_{k+1}^{(T-1)} = (1 - \beta_{T-1,k})y_k^{(T-1)} + \beta_{T-1,k} f_{\omega_{T,k+1}}^{(T)}(\widehat{y}_{k+1}^{(T-1)}).$$

**else if** Assumption 3.2 is NOT known to hold **then**

Query the $\mathcal{SO}$ for the sample values of $f^{(T)}$ at $x_k$, obtain $f_{\omega_{T,k+1}}^{(T)}(x_k)$.

Update $y^{(T-1)}$ by

$$y_{k+1}^{(T-1)} = (1 - \beta_{T-1,k})y_k^{(T-1)} + \beta_{T-1,k} f_{\omega_{T,k+1}}^{(T)}(x_{k+1}).$$

**end if**

**for** $j = T-1, \cdots, 2$ **do**

Update the auxiliary variable $\widehat{y}_{k+1}^{(j-1)}$ by

$$\widehat{y}_{k+1}^{(j-1)} = (1 - \frac{1}{\beta_{j-1,k}})y_k^{(j)} + \frac{1}{\beta_{j-1,k}} y_{k+1}^{(j)}.$$

Query the $\mathcal{SO}$ for the sample value of $f^{(j)}$ at $z_{k+1}^{(j-1)}$, obtain $f_{\omega_{j,k+1}}^{(j)}(z_{k+1}^{(j-1)})$.

Update $y^{(j)}$ by

$$y_{k+1}^{(j-1)} = (1 - \beta_{j-1,k})y_k^{(j-1)} + \beta_{j-1,k} f_{\omega_{j,k+1}}^{(j)}(\widehat{y}_{k+1}^{(j-1)}).$$

**end for**

**end for**

### 3.2 Almost Sure Convergence Of $a$-TSCGD

We first investigate whether and under what condition the algorithm converges almost surely. In particular, we provide sufficient conditions of the stepsizes, such that when the problem is convex, the algorithm converges to an optimal solution almost surely, and when the problem is nonconvex, all limiting points of the solution path generated by the algorithm are stationary points almost surely when $F(x)$ has Lipschitz continuous gradient.

**Theorem 3.1** (Almost sure convergence for $a$-TSCGD)**.** Let Assumptions 2.1 and 3.1 hold, and let the stepsizes $\{\alpha_k\}, \{\beta_{T-1,k}\}, \cdots, \{\beta_{1,k}\}$ be such that

$$\sum_{k=0}^{\infty} \alpha_k = \infty, \sum_{k=0}^{\infty} \beta_{T,k} = \infty, \cdots, \sum_{k=0}^{\infty} \beta_{1,k} = \infty,$$

$$\sum_{k=0}^{\infty} \left( \alpha_k^2 + \beta_{T-1,k}^2 + \cdots + \beta_{1,k}^2 + \frac{\alpha_k^2}{\beta_{T-1,k}} + \cdots + \frac{\alpha_k^2}{\beta_{1,k}} \right) < \infty,$$

and

$$\sum_{k=0}^{\infty} \left( \frac{\beta_{T-1,k}^4}{\beta_{T-2,k}^3} + \cdots + \frac{\beta_{2,k}^4}{\beta_{1,k}^3} \right) < \infty.$$

Let $\left\{(x_k, y_k^{(T-1)}, \cdots, y_k^{(1)})\right\}_{k=0}^{\infty}$ be the sequence generated by Algorithm 2 starting with an arbitrary initial point $(x_0, y_0^{(T-1)}, \cdots, y_0^{(1)})$. Then:

(a) If $F$ is convex, the sequence $\{x_k\}_{k=0}^{\infty}$ converges almost surely to a random point in the set of optimal solutions to problem (1.1).

(b) Suppose in addition that Assumption 2.2 holds, $\mathcal{X} = \mathbb{R}^{d_T}$, and all samples generated by the $\mathcal{SO}$ are uniformly bounded. Then any limiting point of the sequence $\{x_k\}_{k=0}^{\infty}$ is a stationary point of problem (1.1) almost surely.

Furthermore, if Assumption 3.2 also holds, i.e., when $f^{(T)}$ has Lipschitz continuous gradient, then if the stepsizes also satisfy

$$\sum_{k=0}^{\infty} \frac{\alpha_k^4}{\beta_k^3} < \infty,$$

the assertions in (a) and (b) also hold.

*Proof Outline.* We provide the proof outline here for the case when the first inner level function $f^{(T)}$ is non-smooth. The analysis for problems with a smooth first inner level function could be derived from the non-smooth case, and we present the details for both cases in Appendix C.

Essentially, we construct a $T$-element super-martingale to derive the almost sure convergence of the algorithm. We denote by $\mathbb{F}_k$ the collection of random variables up to the $k$-th iteration, i.e.,

$$\mathbb{F}_k = \left\{ \{x_i\}_{i=0}^{k}, \{y_i^{(T-1)}\}_{i=0}^{k-1}, \cdots, \{y_i^{(1)}\}_{i=0}^{k-1}, \{\widehat{y}_i^{(T-2)}\}_{i=0}^{k-1}, \cdots, \{\widehat{y}_i^{(1)}\}_{i=0}^{k-1}, \{\omega_{T,i}\}_{i=1}^{k-1}, \cdots, \{\omega_{1,i}\}_{i=1}^{k-1} \right\}.$$

For the first inner level, since $f^{(T)}$ is non-smooth, we construct the super-martingale for this level by Lemma 2.3. With the additional finite fourth-moment Assumption 3.1, we can derive a stronger result in the following lemma.

**Lemma 3.1.** Let Assumptions 2.1 and 3.1 hold, and let $\{(x_k, y_k^{(T-1)}, \cdots, y_k^{(1)})\}_{k=0}^{\infty}$ be the sequence generated by Algorithm 2. Suppose $\mathbb{E}[\|x_{k+1}-x_k\|^4] \le \mathcal{O}(\alpha_k^4)$ for all $k$ and $\alpha_k/\beta_{T-1,k} \to 0$ as $k \to 0$, in addition to Lemma 2.3 (a) (b) and (c), we have:

(a) There exists a constant $S_{T-1} > 0$ such that $\mathbb{E}[\|y_k^{(T-1)} - f^{(T)}(x_k)\|^4] \le S_{T-1}$ for all $k$.

(b) $\mathbb{E}[\|y_{k+1}^{(T-1)} - y_k^{(T-1)}\|^4] \le \mathcal{O}(\beta_{T-1,k}^4)$ for all $k$.

Next, to construct the super-martingale for the accelerating update steps, we present the following lemma.

**Lemma 3.2.** Let Assumption 2.1 and 3.1 hold, and let $\{(x_k, y_k^{(T-1)}, \cdots, y_k^{(1)})\}_{k=1}^{\infty}$ be the sequence generated by Algorithm 2. For $j = T-2, \cdots, 1$, suppose $\mathbb{E}[\|y_{k+1}^{(j+1)} - y_k^{(j+1)}\|^4] \le \mathcal{O}(\beta_{j+1,k}^4)$ for all $k$ and $\beta_{j+1,k}/\beta_{j,k} \to 0$ as $k \to 0$, then there exists a random variable $e_k^{(j)} \in \mathbb{F}_{k+1}$ for all $k$ satisfying $\|y_k^{(j)} - f^{(j+1)}(y_k^{(j+1)})\| \le e_k^{(j)}$ such that:

(a) For all $k$, with probability 1,

$$\mathbb{E}[[e_{k+1}^{(j)}]^2|\mathbb{F}_{k+1}] \le (1 - \frac{\beta_{j,k}}{2})[e_k^{(j)}]^2 + 2\beta_{j,k}^2 V_{j+1} + \mathcal{O}\left(\frac{\mathbb{E}[\|y_{k+1}^{(j+1)} - y_k^{(j+1)}\|^4|\mathbb{F}_{k+1}]}{\beta_{j,k}^3}\right).$$

(b) If $\sum_{k=1}^{\infty} \beta_{j+1,k}^4/\beta_{j,k}^3 < \infty$, we have

$$\sum_{k=1}^{\infty} \frac{\mathbb{E}[\|y_{k+1}^{(j+1)} - y_k^{(j+1)}\|^4|\mathbb{F}_{k+1}]}{\beta_{j,k}^3} < \infty \quad w.p.1.$$

(c) There exists a constant $D_j \ge 0$ such that $\mathbb{E}[e_k^{(j)}]^2 \le D_j$ for all $k$.

(d) There exists a constant $S_j \ge 0$ such that $\mathbb{E}[\|y_k^{(j)} - f^{(j+1)}(y_k^{(j+1)})\|^4] \le S_j$ for all $k$.

(e) $\mathbb{E}[\|y_{k+1}^{(j)} - y_k^{(j)}\|^4] \le \mathcal{O}(\beta_{j,k}^4)$ for all $k$.

By Lemmas 2.1, 2.3, 3.1 and 3.2, we construct the $T$-element super-martingale and show its

convergence by letting

$$
\begin{aligned}
&X_k = \|x_k - x^*\|^2, Y_k^{(T-1)} = \mathbb{E}[\|y_k^{(T-1)} - f^{(T)}(x_k)\|^2|\mathbb{F}_k],\\
&Y_k^{(T-2)} = \mathbb{E}[[e_k^{(T-2)}]^2|\mathbb{F}_k], \cdots, Y_k^{(1)} = \mathbb{E}[[e_k^{(1)}]^2|\mathbb{F}_k],\\
&\eta_k = [\frac{\alpha_k^2}{\beta_{T-1,k}} + \cdots + \frac{\alpha_k^2}{\beta_{1,k}}]C_0, u_k^{(T)} = 2\alpha_k(F(x_k) - F^*),\\
&u_k^{(1)} = u_k^{(2)} = \cdots = u_k^{(T-1)} = 0, c_1 = 2, \cdots, c_{T-2} = 2(T-2), c_{T-1} = T-1,\\
&\mu_k^{(T-1)} = C_T\beta_{T-1,k}^{-1}\mathbb{E}[\|x_{k+1} - x_k\|^2|\mathbb{F}_k] + 2V_T\beta_{T-1,k}^2,\\
&\mu_k^{(T-2)} = 2\beta_{T-2,k}^2 V_{T-1} + \mathcal{O}(\frac{\mathbb{E}[\|y_{k+1}^{(T-1)} - y_k^{(T-1)}\|^4|\mathbb{F}_k]}{\beta_{T-2,k}^3}), \cdots,\\
&\mu_k^{(1)} = 2\beta_{1,k}^2 V_1 + \mathcal{O}(\frac{\mathbb{E}[\|y_{k+1}^{(2)} - y_k^{(2)}\|^4|\mathbb{F}_k]}{\beta_{1,k}^3}),\\
&\mu_k^{(T)} = \alpha_k^2 C_1 C_2 \cdots C_T,\\
&\theta_k^{(1)} = \beta_{1,k}/2, \cdots, \theta_k^{(T-2)} = \beta_{T-2,k}/2, \theta_k^{(T-1)} = \beta_{T-1,k}.
\end{aligned}
$$

The rest of the proof is similar to that of Theorem 2.1, and we defer the details in Appendix C to avoid repetition. Besides, to prove part (b), we construct a $T$-element super-martingale by applying Lemmas 2.2, 2.3, 3.1 and 3.2, and show any limiting point is a stationary point with probability 1 by Lemma 2.6, which is deferred in Appendix C as well.

□

### 3.3 Convergence Rate Results For $a$-TSCGD

In this subsection, we study the rate of convergence of the algorithm. We consider stepsizes of the form

$$\alpha_k = k^{-a}, \beta_{T-1,k} = k^{-b_{T-1}}, \text{ and } \beta_{j,k} = 2k^{-b_j} \text{ for all } j = T-2, \cdots, 1,$$

where $a$ and $b_j$'s are real numbers if the first inner level function $f^{(T)}$ is nonsmooth, and we choose the step-sizes to be

$$\alpha_k = k^{-a}, \text{ and } \beta_{j,k} = 2k^{-b_j} \text{ for all } j = T-1, \cdots, 1,$$

if $f^{(T)}$ is smooth. After optimizing the rate over all $a$ and $b_j$'s, we get the following result for both convex and nonconvex $F(x)$.

**Theorem 3.2** (Convergence rate of $a$-TSCGD)**.** Suppose that Assumptions 2.1, 2.2 and 3.1 hold and $\mathcal{X} = \mathbb{R}^{d_T}$. Let the stepsizes be $\alpha_k = k^{-a}$, $\beta_{T-1,k} = k^{-b_{T-1}}$ and $\beta_{j,k} = 2k^{-b_j}$ for $j = T-2, \cdots, 1$, where $a, b_{T-1}, ..., b_1 \in (0,1)$. If we choose the step-sizes as $a = \frac{4+T}{8+T}$ and $b_j = \frac{j+3}{8+T}$ for $j = T-2, ..., 1$, letting $\{(x_k, y_k^{(T-1)}, \cdots, y_k^{(1)})\}_{k=0}^{\infty}$ be the sequence generated by $a$-TSCGD Algorithm 2, we obtain

$$\frac{\sum_{k=1}^n \mathbb{E}[\|\nabla F(x_k)\|^2]}{n} \leq \mathcal{O}(n^{-4/(8+T)}).$$

Furthermore, if Assumption 3.2 also holds, Algorithm 2 achieves

$$\frac{\sum_{k=1}^{n}\mathbb{E}[\|\nabla F(x_k)\|^2]}{n} \le \mathcal{O}(n^{-4/(7+T)}),$$

with $\alpha_k = k^{-a}$ and $\beta_{j,k} = 2k^{-b_j}$, where $a = \frac{3+T}{7+T}$ and $b_j = \frac{j+3}{7+T}$ for $j = T-1, ..., 1$.

*Proof Outline.* We present the outline of proof here and defer the detailed analysis in Appendix D.

The analysis is similar to Theorem 2.2. For the accelerating update steps, by Lemma 3.2 and Lemma B.1 in the Appendix, we have the following result:

**Lemma 3.3.** Let Assumptions 2.1 and 3.1 hold, and let $\{(x_k, y_k^{(T-1)}, \cdots, y_k^{(1)})\}_{k=0}^{\infty}$ be the sequence generated by Algorithm 2. Then for any accelerated update step, we have for all $k$

$$\mathbb{E}[\|y_{k+1}^{(j)} - f^{(j+1)}(y_{k+1}^{(j+1)})\|^2] \le \mathcal{O}(k^{4(b_j - b_{j+1})}) + \mathcal{O}(k^{-b_j}), \quad j = T-2, \cdots, 1.$$

Under additional Assumption 2.2 that $F$ has Lipschitz gradient, we have the following result.

**Lemma 3.4.** Let Assumptions 2.1, 2.2 and 3.1 hold, and let $\{(x_k, y_k^{(T-1)}, \cdots, y_k^{(1)})\}_{k=0}^{\infty}$ be the sequence generated by Algorithm 2, then we have for all $k$

$$\begin{aligned}&\mathbb{E}[\|\nabla F(x_k)\|^2]\\ \le& 2\alpha_k^{-1}\mathbb{E}\Big[F(x_k)] - 2\alpha_k^{-1}\mathbb{E}[F(x_{k+1})\Big] + \mathcal{O}\left(\mathbb{E}[\|y_{k+1}^{(T-1)} - f^{(T)}(x_k)\|^2]\right) + \mathcal{O}\left(\mathbb{E}[\|y_{k+1}^{(T-2)} - f^{(T-1)}(y_{k+1}^{(T-1)})\|^2]\right)\\ &+ \cdots + \mathcal{O}\left(\mathbb{E}[\|y_{k+1}^{(1)} - f^{(2)}(y_{k+1}^{(2)})\|^2]\right) + \mathcal{O}(\alpha_k).\end{aligned}$$

Summing up the inequalities in the previous lemma from $k = 0$ to $n$, by Lemma 2.7 and Lemma 3.3, we obtain

$$\begin{aligned}\frac{\sum_{k=1}^{n}\mathbb{E}[\|\nabla F(x_k)\|^2]}{n} &\le \mathcal{O}\Big(n^{a-1} + n^{-2a+2b_{T-1}}\mathbb{I}_{2(a-b_{T-1})=1}^{\log n} + n^{-b_{T-1}} + n^{-a}\Big)\\ &\quad + \mathcal{O}\Big(\sum_{j=1}^{T-2}[n^{4(b_j-b_{j+1})}\mathbb{I}_{4(b_{j+1}-b_j)=1}^{\log n} + n^{-b_j}]\Big)\\ &\le \mathcal{O}(n^{-4/(8+T)}),\end{aligned}$$

by choosing $a = \frac{4+T}{8+T}$ and $b_j = \frac{3+j}{8+T}$ for $j = T-1, \cdots, 1$.

Furthermore, if Assumption 3.2 also holds, i.e., the first inner level function $f^{(T)}$ has Lipschitz continuous gradient, then the first inner level could also be updated by the accelerating update rule. By similar analysis as in Lemma 3.3, we have for all $k$,

$$\mathbb{E}[\|y_{k+1}^{(T-1)} - f^{(T)}(x_{k+1})\|^2] \le \mathcal{O}(k^{4(b_{T-1}-a)}) + \mathcal{O}(k^{-b_{T-1}}).$$

Combine this inequality with Lemmas 3.3 and 3.4, by choosing $a = \frac{3+T}{7+T}$ and $b_j = \frac{3+j}{7+T}$ for $j = T-1, \cdots, 1$, we obtain

$$\frac{\sum_{k=1}^{n}\mathbb{E}[\|\nabla F(x_k)\|^2]}{n} \le \mathcal{O}(n^{-4/(7+T)}),$$

which completes the proof. □

This result shows that one can solve the multi-level composition problem using few calls to the sampling oracle when individual component functions are smooth. In the special case where $T = 2$, when the first inner level is smooth, our result strictly improves the convergence rate of the a-SCGD in [26] from $\mathcal{O}(n^{-2/7})$ to $\mathcal{O}(n^{-4/9})$. In this case our result matches the convergence rate by ASC-PG in [28]. To the best of our knowledge, our results for the $T$-level problem strictly improve and generalize existing results which work for the case where $T = 2$.

Next we investigate the convergence rate of Algorithm 2 for *optimally strongly convex* objective functions. Denote by $\mathcal{X}^*$ the set of optimal solutions $x^*$ to problem (1.1). We say that the objective function $F$ is *optimally strongly convex* with parameter $\lambda > 0$ if

$$F(x) - F\big(\Pi_{\mathcal{X}^*}(x)\big) \geq \lambda \|x - \Pi_{\mathcal{X}^*}(x)\|^2, \quad \forall x \in \mathcal{X}. \tag{3.1}$$

Clearly, the class of optimally strongly convex functions strictly contains all strongly convex functions, and is thus more general.

In the next theorem, we prove that for optimally strongly convex objective, our algorithm converges faster. We defer the detailed proof to Appendix E.

**Theorem 3.3** (Convergence rate of $a$-TSCGD for strongly convex problems)**.** Let Assumptions 2.1, 2.2 and 3.1 hold. Suppose that the objective function $F(x)$ in (1.1) is optimally strongly convex with some parameter $\lambda > 0$ defined in (3.1) . Set $\alpha_k = \frac{1}{\lambda} k^{-a}$, $\beta_{T-1,k} = k^{-b_{T-1}}$ and $\beta_{j,k} = 2k^{-b_j}$ for $j = T-2, \cdots, 1$. Let $\{(x_k, y_k^{(T-1)}, \cdots, y_k^{(1)})\}_{k=0}^{\infty}$ be the sequence generated by $a$-TSCGD Algorithm 2, then

$$\mathbb{E}[\|x_n - \Pi_{\mathcal{X}^*}(x_n)\|^2] \leq \mathcal{O}\Big(n^{-a} + n^{-2(a-b_{T-1})} + n^{-b_{T-1}} + \sum_{j=1}^{T-2} [n^{-4(b_{j+1}-b_j)} + n^{-b_j}]\Big).$$

With the choice of $a = 1, b_{T-1} = \frac{2+T}{4+T}, b_{T-2} = \frac{1+T}{4+T}, \cdots, b_1 = \frac{4}{4+T}$, we have

$$\mathbb{E}[\|x_n - \Pi_{\mathcal{X}^*}(x_n)\|^2] \leq \mathcal{O}(n^{-4/(4+T)}).$$

Furthermore, if Assumption 3.2 also holds, Algorithm 2 achieves

$$\mathbb{E}[\|x_n - \Pi_{\mathcal{X}^*}(x_n)\|^2] \leq \mathcal{O}(n^{-4/(3+T)}),$$

with the stepsizes being $\alpha_k = \frac{1}{\lambda} k^{-a}$ and $\beta_{j,k} = 2k^{-b_j}$, where $a = 1$ and $b_j = \frac{3+j}{3+T}$ for $j = T-1, \cdots, 1$.

This result shows that our algorithm achieves a faster convergence for those problems of optimally strongly convexity in the objective functions. For the speical case $T = 1$ with a smooth strongly convex function, this result achieves a convergence rate of $\mathcal{O}(n^{-1})$, which meets the convergence rate of the single-level strongly convex stochastic optimization. Besides, for a special case $T = 2$ with a smooth first inner level function, this result achieves a convergence rate of $\mathcal{O}(n^{-4/5})$, which matches the convergence rate ASC-PG in [28] for optimally strongly convex problems.

# 4 Examples of Applications

In this section, we provide two motivating applications of the $T$-level stochastic compositional optimization problem (1.1). The first motivating application is the risk-averse stochastic optimization.

Risk-averse stochastic optimization finds wide applications in many fields such as risk management [4] and government planning [3]. Among different formulations of risk-averse stochastic optimization problems, one particular important problem is the mean-deviation risk-averse optimization problem that

$$\max_x \rho(U(x,w)) := \max_x \left\{ \mathbb{E}_\omega \left[U(x,\omega)\right] - \lambda \mathbb{E}\Big[ \left(\mathbb{E}\left[U(x,\omega)\right] - U(x,\omega)\right)_+^p \Big]^{1/p} \right\}. \tag{4.1}$$

Here the objective $\rho$ is the composition of three expected-value functions. It is also a law-invariant coherent risk measure. See [21, 1] for more detailed discussions.

This problem falls into the problem class (1.1) as a three-level stochastic compositional optimization problem. In particular, the problem is equivalent to

$$\min_x (f^{(1)} \circ f^{(2)} \circ f^{(3)})(x),$$

where

$$f^{(1)}\big((y_1, y_2)\big) = y_1 - y_2^{1/p}, f^{(2)}(z, x) = \big(z, \mathbb{E}_\omega\left[(z - U(x,\omega))_+^p\right]\big), \text{ and } f^{(3)}(x) = \left(\mathbb{E}_\omega\left[U(x,\omega)\right], x\right).$$

Another example is multi-level optimization problem. A $T$-level optimization problem takes the following general form (see textbook [24] for more details)

$$\min_{x_1} \mathbb{E}_{\omega_1}\Big[ \min_{x_2}(\mathbb{E}_{\omega_2|\omega_1}\Big[ \cdots \min_{x_T}(\mathbb{E}_{\omega_T|\omega_1,\omega_2,\cdots,\omega_{T-1}}\Big[U(x_1,\omega_1,x_2,\omega_2,\cdots,x_T,\omega_T)\Big])\Big])\Big],$$

where $x_1, \ldots, x_T$ are decision variables at levels from 1 to $T$, $\omega_1, \ldots, \omega_T$ are random variables which are revealed after each level, and $U(x_1,\omega_1,x_2,\omega_2,\cdots,x_T,\omega_T)$ is some utility function. In the case where the stochastic process $\omega_1, \ldots, \omega_T$ is generated by a random walk on a finite number of states, the problem becomes

$$\min_{x_1} \mathbb{E}_{\omega_1}\Big[ \min_{x_2}(\mathbb{E}_{\omega_2|\omega_1}\Big[ \cdots \min_{x_T}(\mathbb{E}_{\omega_T|\omega_{T-1}}\Big[U(x_1,\omega_2,x_2,\omega_3,\cdots,x_T,\omega_T)\Big])\Big])\Big].$$

It can be viewed as an extension of finite-horizon reinforcement learning, in which the overall objective is no longer additive with respect to levels and decisions are continuous. This problem takes a form similar to (1.1), especially when the state space and decision space are finite and discrete. Note that in this problem, $f^{(1)}(\cdot) = \min_{x_2}(\cdot)$, $f^{(2)}(\cdot) = \min_{x_3}(\cdot), \cdots$, and $f^{(T-1)}(\cdot) = \min_{x_T}(\cdot)$, where all of them are non-differentiable. Thus we may use a smooth approximation (e.g., a softmax operator) to replace these $\min_{x_2}, \cdots$, $\min_{x_T}$ operators, then apply the $T$-level SCGD methods to solve the approximate $T$-level optimization problem.

# 5 Numerical Experiments

In this section, we conduct numerical experiments. We consider the risk-averse stochastic optimization in a regression setting. In particular, consider a linear model $Y = X\beta^* + \epsilon$, where we assume all samples of $X$ and $\epsilon$ are independently and identically distributed. Our goal is to estimate $\beta^*$, and we consider a risk-averse formulation. Consider the risk-averse optimization problem (4.1). Denoting the $i$-th sample by $\omega_i = \{x_i, y_i\}$, we take

$$U(\beta, \omega_i) = -(y_i - x_i^T \beta)^2,$$

and we set $p = 2$. To the best of our knowledge, our algorithm is the first gradient-based method which can be adopted to solve this 3-level stochastic optimization problem. We point out that this approach of risk-averse regression tends to provide "stable" solutions. This defines a general notion of stability in statistics in [17, 25], where the stability is usually defined as variance, and we also penalize the "good" cases when the empirical error is smaller than its expectation. In comparison, in our approach, we do not penalize these "good" cases.

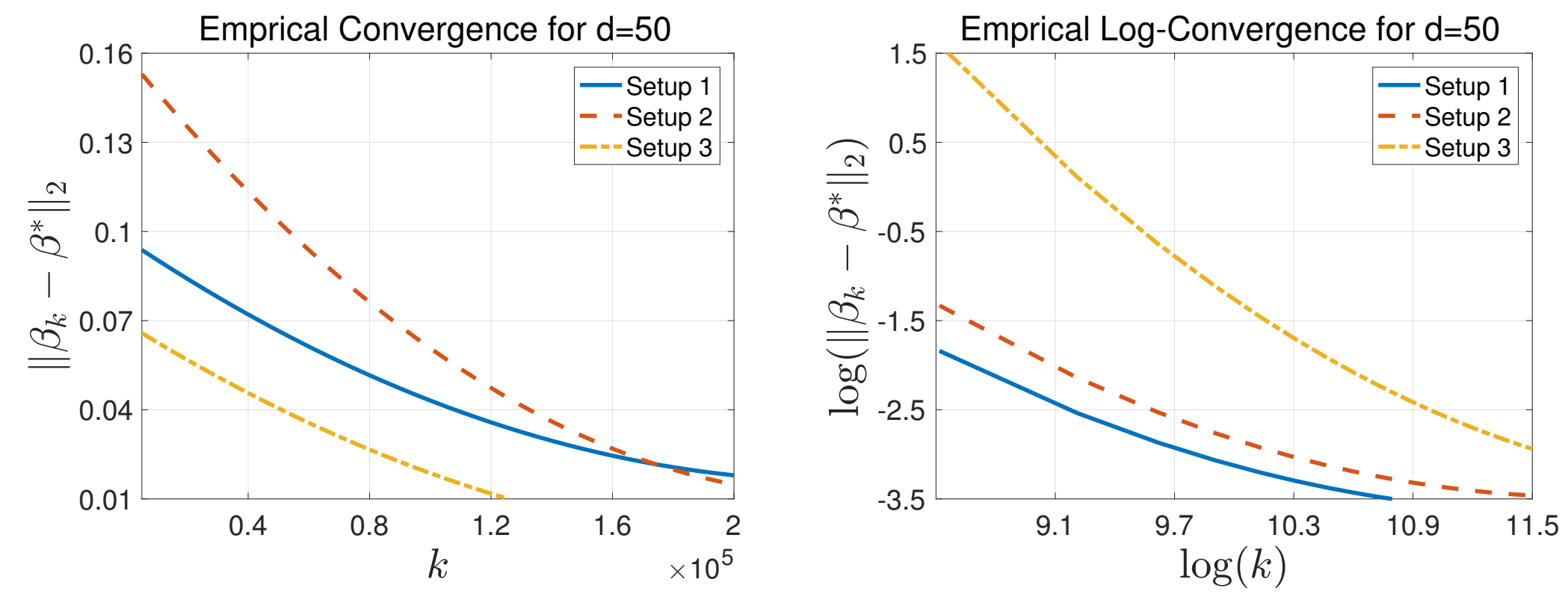

Figure 2: Averaged difference between generated solution and the optimal solution and empirical convergence rate when $d = 50$

Let the dimension of the covariate $x_i$ be $d$. We consider three setups to generate the data that

- Setup 1: $X \sim N(0, I_d)$.
- Setup 2: $X \sim N(0, \Sigma)$, where $\Sigma_{jj} = 1$ and $\Sigma_{jk} = 0.5$ for $j, k = 1, ..., d$ and $j \neq k$.
- Setup 3: $X \sim N(0, \Sigma)$, where $\Sigma_{jk} = 0.5e^{-\frac{|j-k|}{d}}$.

Since our problem is convex, by our theoretical analysis, the generated sequence of solutions converges to the optimal solution. As the true optimal solution is unknown (Note that $\beta^*$ is not necessarily the optimal solution), we take the solution after 500,000 iterations as the optimal solution. We run 100 replications and plot the averaged difference between the solution at the $k$-th iteration $\beta_k$ and the optimal solution $\widehat{\beta}$ in Figures 2, 3, 4 and 5.

Meanwhile, in all setups, we draw the error $\epsilon$ and generate each component of $\beta^* \in \mathbb{R}^d$ independently from a standard normal distribution. We also consider different $d \in \{50, 100, 150, 200\}$, which are specified in the Figures. In each iteration of the algorithm, we draw a new sample of $X$ and $Y$, and update the solution using Algorithm 2.

Besides, to further investigate empirical rates of convergence under all different settings, we plot the averaged $\log(k)$ vs. $\log(\|\beta_k - \widehat{\beta}\|)$ after 100 replications in the figures, where $\widehat{\beta}$ is the optimal solution. We find that for all cases, the slopes of the lines are close to $-2/5$, which matches our theoretical analysis that our algorithm converges at a rate of $\mathcal{O}(k^{-2/5})$ for three-level problems.

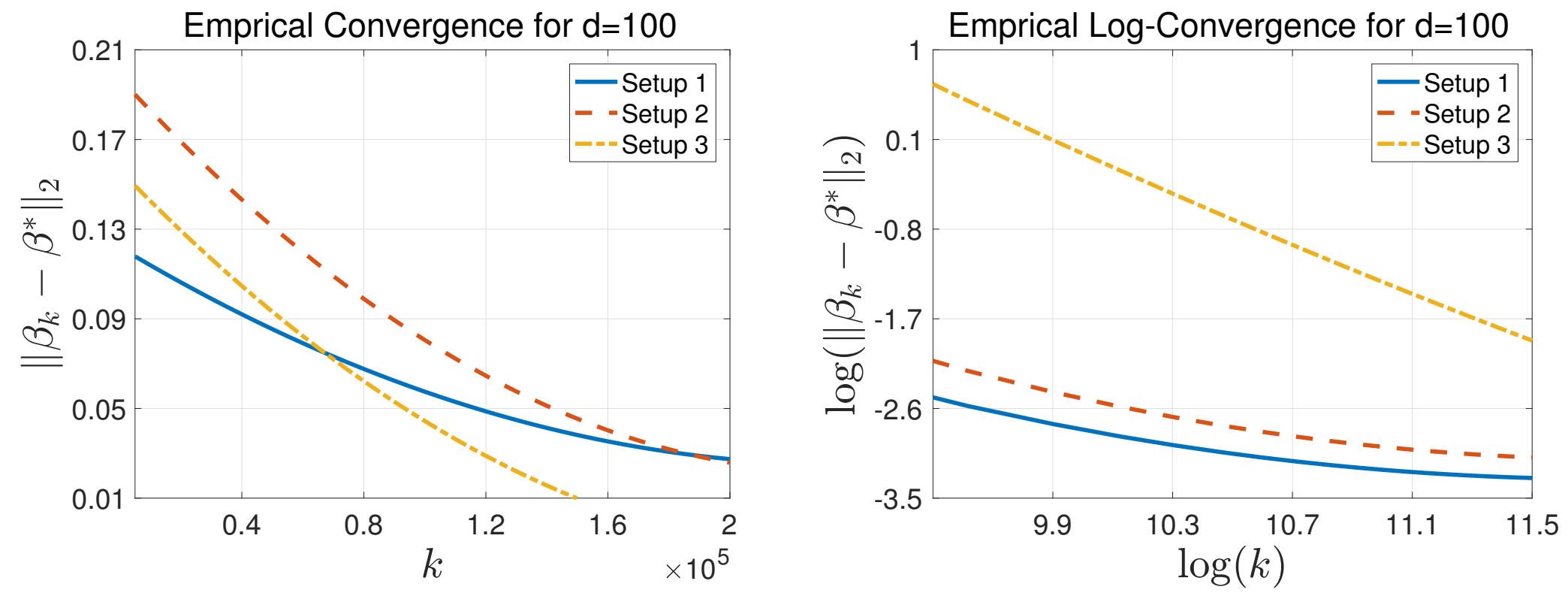

Figure 3: Averaged difference between generated solution and the optimal solution and empirical convergence rate when $d = 100$.

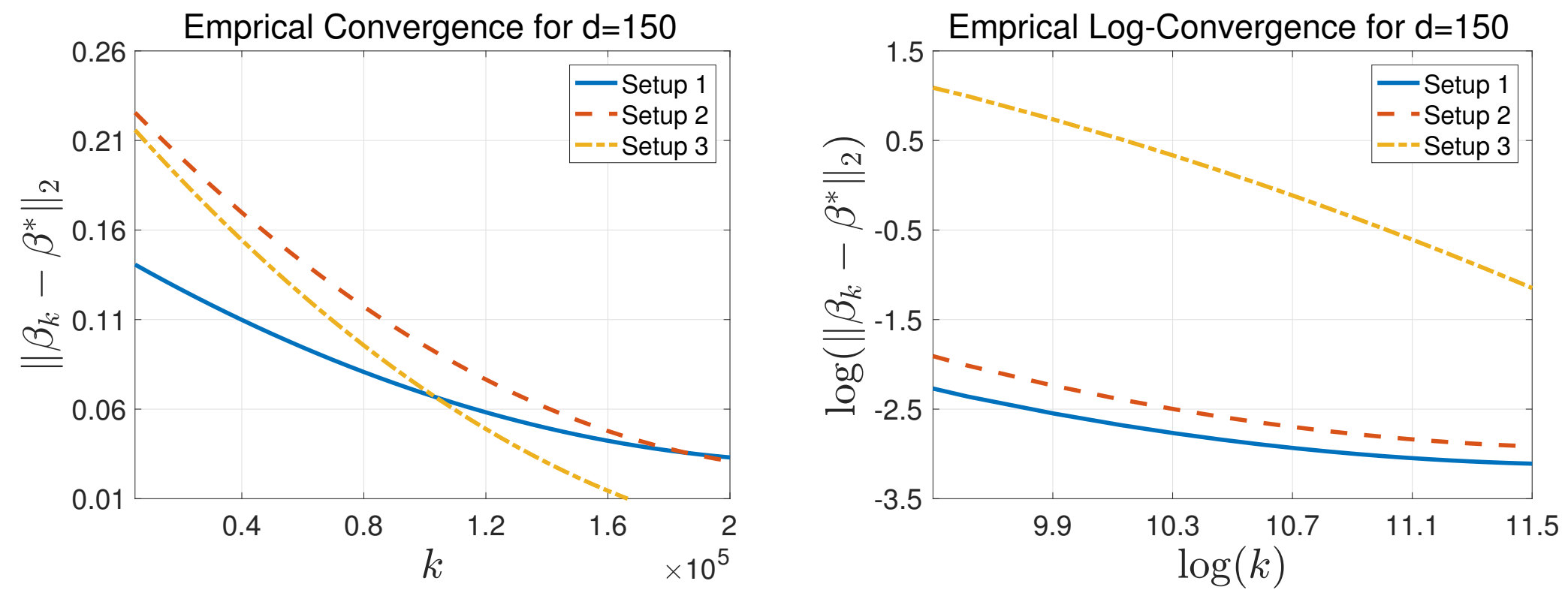

Figure 4: Averaged difference between generated solution and the optimal solution and empirical convergence rate when $d = 150$.

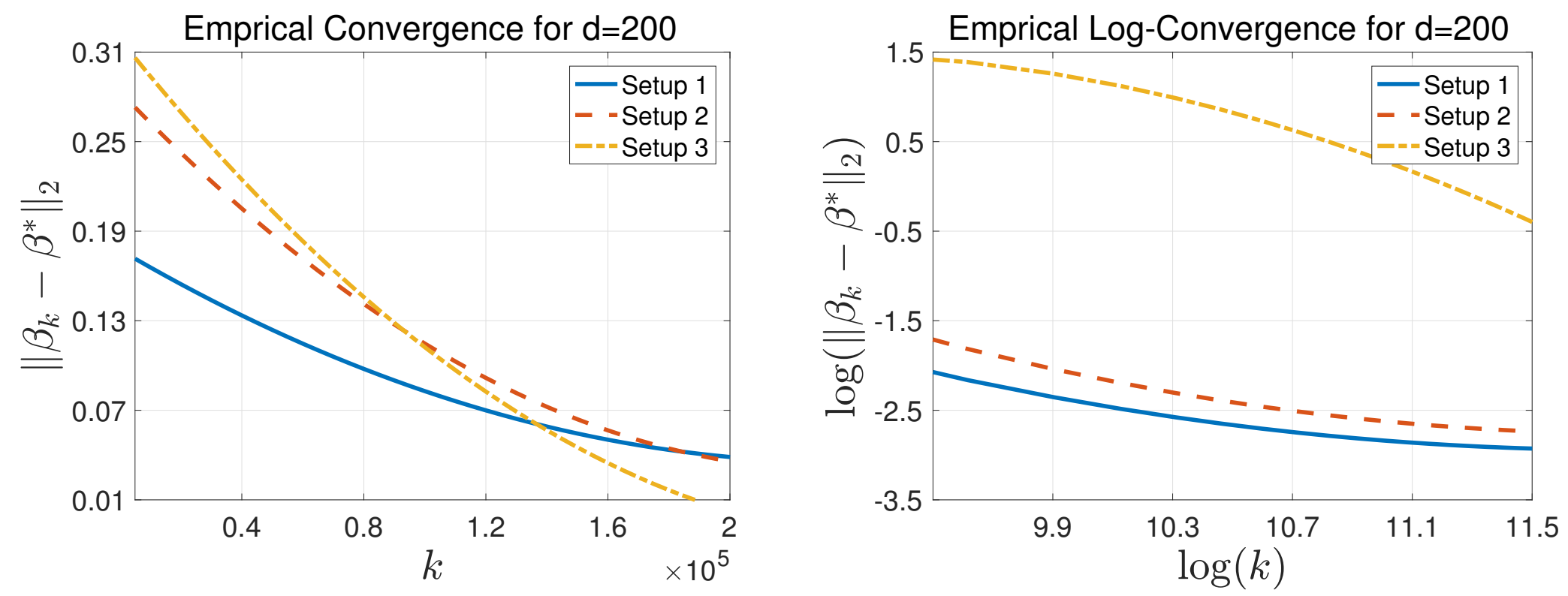

Figure 5: Averaged difference between generated solution and the optimal solution and empirical convergence rate when $d = 200$.

## 6 Conclusion

In this paper, we propose the first gradient-type algorithms for a class of multi-level stochastic compositional optimization problems. We provide strong theoretical guarantees for our algorithms. In particular, we prove almost sure convergence results that when the problem is convex, our algorithm converges to an optimal solution, and when the problem is nonconvex, every limiting point of the sequence of solutions is an stationary point. Under various assumptions, we further characterize the rates of convergence of our algorithms. In the case where $T = 2$, our convergence rate result matches and strictly generalizes the best known result by [28]. In the case where $T \geq 3$, our results provide the first few benchmarks on the sample complexity for solving multi-level stochastic optimization problems.

There are several interesting future research questions. First, our convergence rate result requires that the inner-level functions $f^{(2)}, \cdots, f^{(T)}$ be smooth. It is unclear how to achieve fast convergence when some of these functions are non-smooth. Second, it is not clear whether the convergence rate can be improved or not. We are not aware of any sample complexity lower bound for the multi-level stochastic optimization problem. Third, it is of practical interest to consider the special case where all expectations are finite sums. In this case, one may conjecture that variance reduction can be used to further improve the algorithms' efficiency.

Appendix to “Multi-Level Stochastic Compositonal Optimization”

Shuoguang Yang* Mengdi Wang† Ethan X. Fang‡

# A Proof of Theorem 2.1

In this section, we present the detailed proof for Theorem 2.1.

## A.1 Proof of Lemma 2.1

Before presenting the detail proof of Lemma 2.1, we present a lemma that is used in proving Lemma 2.1.

**Lemma A.1.** Suppose Assumption 2.1 holds, and let $\{(x_k, y_k^{(T-1)}, \cdots, y_k^{(1)})\}_{k=0}^{\infty}$ be the sequence generated by Algorithm 1. Denote by $\widetilde{\nabla} F_{\omega_k}(x_k) \equiv \widetilde{\nabla} f^{(T)}_{\omega_{T,k}}(x_k) \nabla f^{(T-1)}_{\omega_{T-1,k}}(f^{(T)}(x_k)) \cdots \nabla f^{(1)}_{\omega_{1,k}}(f^{(2)} \circ \cdots \circ f^{(T)}(x_k))$, and let $X_k \in \mathbb{F}_k$ be a vector of random variables, where $\mathbb{F}_k$ is the collection of random variables

$$\Big\{\{x_i\}_{i=0}^{k}, \{y_i^{(T-1)}\}_{i=0}^{k-1}, \cdots, \{y_i^{(1)}\}_{i=0}^{k-1}, \{\omega_{T,i}\}_{i=0}^{k-1}, \cdots, \{\omega_{1,i}\}_{i=0}^{k-1}\Big\}.$$

Then there exists a constant $C_0 > 0$ dependent only on the number of levels $T$ such that for all $k$, with probability 1,

$$\begin{aligned}
&X_k' \mathbb{E}\Big[\widetilde{\nabla} F_{\omega_k}(x_k) - \widetilde{\nabla} f^{(T)}_{\omega_{T,k}}(x_k) \nabla f^{(T-1)}_{\omega_{T-1,k}}(y_k^{(T-1)}) \cdots \nabla f^{(1)}_{\omega_{1,k}}(y_k^{(1)}) \Big| \mathbb{F}_k\Big] \\
\le& (T-1)\beta_{T-1,k}\mathbb{E}[\|y_k^{(T-1)} - f^{(T)}(x_k)\|^2 | \mathbb{F}_k] + (T-2)\beta_{T-2,k}\mathbb{E}[\|y_k^{(T-2)} - f^{(T-1)}(y_k^{(T-1)})\|^2 | \mathbb{F}_k] \\
&+ \cdots + \beta_{1,k}\mathbb{E}[\|y_k^{(1)} - f^{(2)}(y_k^{(2)})\|^2 | \mathbb{F}_k] + C_0\Big(\frac{1}{\beta_{T-1,k}} + \cdots + \frac{1}{\beta_{1,k}}\Big)\|X_k\|^2.
\end{aligned}$$

*Department of Industrial Engineering and Operations Research, Columbia University, New York City, NY, USA; e-mail: sy2614@columbia.edu

†Department of Operations Research and Financial Engineering, Princeton University, Princeton, NJ, USA; email: mengdiw@princeton.edu

‡Department of Statistics and Department of Industrial and Manufacturing Engineering, Pennsylvania State University, University Park, PA, USA; email: xxf13@psu.edu

PROOF: We begin our analysis by the chain rule as follows:

$$\begin{aligned}
&\widetilde{\nabla} F_{\omega_k}(x_k) - \widetilde{\nabla} f^{(T)}_{\omega_{T,k}}(x_k)\nabla f^{(T-1)}_{\omega_{T-1,k}}(y_k^{(T-1)})\cdots\nabla f^{(1)}_{\omega_{1,k}}(y_k^{(1)})\\
=&\widetilde{\nabla} f^{(T)}_{\omega_{T,k}}(x_k)\nabla f^{(T-1)}_{\omega_{T-1,k}}(f^{(T)}(x_k))\nabla f^{(T-2)}_{\omega_{T-2,k}}(f^{(T-1)}(f^{(T)}(x_k)))\cdots\nabla f^{(1)}_{\omega_{1,k}}(f^{(2)}\circ\cdots\circ f^{(T)}(x_k))\\
&-\widetilde{\nabla} f^{(T)}_{\omega_{T,k}}(x_k)\nabla f^{(T-1)}_{\omega_{T-1,k}}(y_k^{(T-1)})\nabla f^{(T-2)}_{\omega_{T-2,k}}(f^{(T-1)}(f^{(T)}(x_k)))\cdots\nabla f^{(1)}_{\omega_{1,k}}(f^{(2)}\circ\cdots\circ f^{(T)}(x_k))\\
&+\widetilde{\nabla} f^{(T)}_{\omega_{T,k}}(x_k)\nabla f^{(T-1)}_{\omega_{T-1,k}}(y_k^{(T-1)})\nabla f^{(T-2)}_{\omega_{T-2,k}}(f^{(T-1)}(f^{(T)}(x_k)))\cdots\widetilde{\nabla} f^{(1)}_{\omega_{1,k}}(f^{(2)}\circ\cdots\circ f^{(T)}(x_k))\\
&-\widetilde{\nabla} f^{(T)}_{\omega_{T,k}}(x_k)\nabla f^{(T-1)}_{\omega_{T-1,k}}(y_k^{(T-1)})\nabla f^{(T-2)}_{\omega_{T-2,k}}(y_k^{(T-2)})\nabla f^{(T-3)}_{\omega_{T-3,k}}(f^{(T-2)}(f^{(T-1)}(f^{(T)}(x_k))))\cdots\\
&\ \nabla f^{(1)}_{\omega_{1,k}}(f^{(2)}\circ\cdots\circ f^{(T)}(x_k))\\
&\vdots\\
&+\widetilde{\nabla} f^{(T)}_{\omega_{T,k}}(x_k)\nabla f^{(T-1)}_{\omega_{T-1,k}}(y_k^{(T-1)})\nabla f^{(T-2)}_{\omega_{T-2,k}}(y_k^{(T-2)})\cdots\nabla f^{(2)}_{\omega_{2,k}}(y_k^{(2)})\nabla f^{(1)}_{\omega_{1,k}}(f^{(2)}\circ\cdots\circ f^{(T)}(x_k))\\
&-\widetilde{\nabla} f^{(T)}_{\omega_{T,k}}(x_k)\nabla f^{(T-1)}_{\omega_{T-1,k}}(y_k^{(T-1)})\nabla f^{(T-2)}_{\omega_{T-2,k}}(y_k^{(T-2)})\cdots\nabla f^{(2)}_{\omega_{2,k}}(y_k^{(2)})\nabla f^{(1)}_{\omega_{1,k}}(y_k^{(1)}).
\end{aligned}$$

Denote by $S_m = \widetilde{\nabla} f^{(T)}_{\omega_{T,k}}(x_k^{(T)})\cdots\nabla f^{(m)}_{\omega_{m,k}}(y_k^{(m)})\nabla f^{(m-1)}_{\omega_{m-1,k}}(f^{(m)}\circ\cdots\circ f^{(T)}(x_k))\cdots\nabla f^{(1)}_{\omega_{1,k}}(f^{(2)}\circ\cdots\circ f^{(T)}(x_k))$. Clearly, $S_T = \widetilde{\nabla} F_{\omega_k}(x_k)$ and $S_1 = \widetilde{\nabla} f^{(T)}_{\omega_{T,k}}(x_k)\nabla f^{(T-1)}_{\omega_{T-1,k}}(y_k^{(T-1)})\cdots\nabla f^{(1)}_{\omega_{1,k}}(y_k^{(1)})$, and we have

$$\begin{aligned}
&\widetilde{\nabla} F_{\omega_k}(x_k) - \widetilde{\nabla} f^{(T)}_{\omega_{T,k}}(x_k)\nabla f^{(T-1)}_{\omega_{T-1,k}}(y_k^{(T-1)})\cdots\nabla f^{(1)}_{\omega_{1,k}}(y_k^{(1)})\\
=&(S_T - S_{T-1}) + (S_{T-1} - S_{T-2}) + \cdots + (S_2 - S_1).
\end{aligned} \tag{A.1}$$

Considering $S_m - S_{m-1}$, by definition, we obtain

$$\begin{aligned}
&\|S_m - S_{m-1}\|\\
=&\|\widetilde{\nabla} f^{(T)}_{\omega_{T,k}}(y_k^{(T)})\cdots\nabla f^{(m)}_{\omega_{m,k}}(y_k^{(m)})P_m\nabla f^{(m-2)}_{\omega_{m-2,k}}(f^{(m-1)}\circ\cdots\circ f^{(T)}(x_k))\cdots\nabla f^{(1)}_{\omega_{1,k}}(f^{(2)}\circ\cdots\circ f^{(T)}(x_k))\|\\
\le&M_m\|P_m\|,
\end{aligned} \tag{A.2}$$

where $P_m = \nabla f^{(m-1)}_{\omega_{m-1,k}}(f^{(m)}\circ\cdots\circ f^{(T)}(x_k)) - \nabla f^{(m-1)}_{\omega_{m-1,k}}(y_k^{(m-1)})$ and

$M_m = \|\widetilde{\nabla} f^{(T)}_{\omega_{T,k}}(y_k^{(T)})\|\cdots\|\nabla f^{(m)}_{\omega_{m,k}}(y_k^{(m)})\|\|\nabla f^{(m-2)}_{\omega_{m-2,k}}(f^{(m-1)}\circ\cdots\circ f^{(T)}(x_k))\|\cdots\|\nabla f^{(1)}_{\omega_{1,k}}(f^{(2)}\circ\cdots\circ f^{(T)}(x_k))\|$.

By Assumption 2.1 (iii)-(iv), we have

$$\mathbb{E}[M_m^2|\mathbb{F}_k] \le C_T C_{T-1}\cdots C_m C_{m-2}\cdots C_1. \tag{A.3}$$

Consider $P_m$,

$$\begin{aligned}
&\|P_m\| = \|\nabla f^{(m-1)}_{\omega_{m-1,k}}(f^{(m)} \circ \cdots \circ f^{(T)}(x_k))) - \nabla f^{(m-1)}_{\omega_{m-1,k}}(y_k^{(m-1)})\| \\
\leq & L_{m-1}\|f^{(m)} \circ \cdots \circ f^{(T)}(x_k) - y_k^{(m-1)}\| \\
\leq & L_{m-1}\|f^{(m)} \circ \cdots \circ f^{(T)}(x_k) - f^{(m)} \circ \cdots \circ f^{(T-1)}(y_k^{(T-1)})\| \\
& + L_{m-1}\|f^{(m)} \circ \cdots \circ f^{(T-1)}(y_k^{(T-1)}) - f^{(m)} \circ \cdots \circ f^{(T-2)}(y_k^{(T-2)})\| \\
& + \cdots + L_{m-1}\|f^{(m)}(y_k^{(m)}) - y_k^{(m-1)}\| \\
\leq & L_{m-1}\sqrt{C_m \cdots C_{T-1}}\|f^{(T)}(x_k) - y_k^{(T-1)}\| + L_{m-1}\sqrt{C_m \cdots C_{T-2}}\|f^{(T-1)}(y_k^{(T-1)}) - y_k^{(T-2)}\| \\
& + \cdots + L_{m-1}\|f^{(m)}(y_k^{(m)}) - y_k^{(m-1)}\|,
\end{aligned} \tag{A.4}$$

where the first inequality holds by Assumption 2.1 (iv) as $f^{(m-1)}_{\omega_{m-1,k}}$ has Lipschitz continuous gradient with parameter $L_j$, and the last inequality holds by Assumption 2.1 (iii)-(iv) that $f^{(j)}$ is Lipschitz continuous with parameter $C_j$ for all $j$'s.

Substituting Eq.(A.4) into Eq.(A.2) yields

$$\begin{aligned}
&\|X_k\|\|S_m - S_{m-1}\| \leq M_m\|X_k\|\|P_m\| \\
\leq & M_m\|X_k\|\Big(L_{m-1}\sqrt{C_m \cdots C_{T-1}}\|f^{(T)}(x_k) - y_k^{(T-1)}\| + L_{m-1}\sqrt{C_m \cdots C_{T-2}}\|f^{(T-1)}(y_k^{(T-1)}) - y_k^{(T-2)}\| \\
& + \cdots + L_{m-1}\|f^{(m)}(y_k^{(m)}) - y_k^{(m-1)}\|\Big) \\
\leq & \beta_{T-1,k}\|f^{(T)}(x_k) - y_k^{(T-1)}\|^2 + \beta_{T-2,k}\|f^{(T-1)}(y_k^{(T-1)}) - y_k^{(T-2)}\|^2 + \cdots + \beta_{m-1,k}\|f^{(m)}(y_k^{(m)}) - y_k^{(m-1)}\|^2 \\
& + M_k^2\|X_k\|^2\Big(\frac{L_{m-1}^2 C_m \cdots C_{T-1}}{4\beta_{T-1,k}} + \frac{L_{m-1}^2 C_m \cdots C_{T-2}}{4\beta_{T-2,k}} + \cdots + \frac{L_{m-1}^2}{4\beta_{m-1,k}}\Big),
\end{aligned}$$

where the last inequality holds by the fact that $2xy \leq ax^2 + \frac{1}{a}y^2$ for any $x, y \in \mathbb{R}$ and $a > 0$. Taking expectation on both sides of the previous inequality and combine it with Eq.(A.3) , since $X_k \in \mathbb{F}_k$, there exists a constant $R_m > 0$ such that almost surely

$$\begin{aligned}
&\mathbb{E}\Big[\|X_k\|\|S_m - S_{m-1}\|\Big|\mathbb{F}_k\Big] \\
\leq & \beta_{T-1,k}\mathbb{E}[\|f^{(T)}(x_k) - y_k^{(T-1)}\|^2|\mathbb{F}_k] + \beta_{T-2,k}\mathbb{E}[\|f^{(T-1)}(y_k^{(T-1)}) - y_k^{(T-2)}\|^2|\mathbb{F}_k] \\
& + \cdots + \beta_{m-1,k}\mathbb{E}[\|f^{(m)}(y_k^{(m)}) - y_k^{(m-1)}\|^2|\mathbb{F}_k] + R_m\Big(\frac{1}{\beta_{T-1,k}} + \cdots + \frac{1}{\beta_{m-1,k}}\Big)\|X_k\|^2.
\end{aligned} \tag{A.5}$$

Meanwhile, we have

$$\begin{aligned}
& X_k'\mathbb{E}\Big[\widetilde{\nabla} F_{\omega_k}(x_k) - \widetilde{\nabla} f^{(T)}_{\omega_{T,k}}(x_k)\nabla f^{(T-1)}_{\omega_{T-1,k}}(y_k^{(T-1)}) \cdots \nabla f^{(1)}_{\omega_{1,k}}(y_k^{(1)})\Big|\mathbb{F}_k\Big] \\
\leq & \|X_k\|\mathbb{E}\Big[\|S_T - S_{T-1}\| + \|S_{T-1} - S_{T-2}\| + \cdots + \|S_2 - S_1\|\Big|\mathbb{F}_k\Big].
\end{aligned} \tag{A.6}$$

Substituting Eq.(A.5) into Eq.(A.6) and sum up from $m = 2$ to $m = T$, with some algebraic manupulation, we conclude that there exists a constant $C_0 > 0$ dependent only on the number of levels $T$ such that with probability 1,

$$
\begin{aligned}
& X_k'\mathbb{E}\Big[\widetilde{\nabla} F_{\omega_k}(x_k) - \widetilde{\nabla} f^{(T)}_{\omega_{T,k}}(x_k)\nabla f^{(T-1)}_{\omega_{T-1,k}}(y_k^{(T-1)})\cdots\nabla f^{(1)}_{\omega_{1,k}}(y_k^{(1)})\Big|\mathbb{F}_k\Big] \\
\le & (T-1)\beta_{T-1,k}\mathbb{E}[\|y_k^{(T-1)} - f^{(T)}(x_k)\|^2|\mathbb{F}_k] + (T-2)\beta_{T-2,k}\mathbb{E}[\|y_k^{(T-2)} - f^{(T-1)}(y_k^{(T-1)})\|^2|\mathbb{F}_k] \\
& + \cdots + \beta_{1,k}\mathbb{E}[\|y_k^{(1)} - f^{(2)}(y_k^{(2)})\|^2|\mathbb{F}_k] + C_0\Big(\frac{1}{\beta_{T-1,k}} + \cdots + \frac{1}{\beta_{1,k}}\Big)\|X_k\|^2,
\end{aligned}
$$

which completes the proof. □

Next, we present the proof of Lemma 2.1.

Proof of lemma 2.1:

$$
\begin{aligned}
& \|x_{k+1} - x^*\|^2 \\
= & \|\Pi_{\mathcal{X}}\{x_k - \alpha_k \widetilde{\nabla} f^{(T)}_{\omega_{T,k}}(x_k)\nabla f^{(T-1)}_{\omega_{T-1,k}}(y_k^{(T-1)})\cdots\nabla f^{(1)}_{\omega_{1,k}}(y_k^{(1)})\} - x^*\|^2 \\
\le & \|x_k - x^* - \alpha_k \widetilde{\nabla} f^{(T)}_{\omega_{T,k}}(x_k)\nabla f^{(T-1)}_{\omega_{T-1,k}}(y_k^{(T-1)})\cdots\nabla f^{(1)}_{\omega_{1,k}}(y_k^{(1)})\|^2 \\
= & \|x_k - x^*\|^2 - 2\alpha_k(x_k - x^*)'\widetilde{\nabla} f^{(T)}_{\omega_{T,k}}(x_k)\nabla f^{(T-1)}_{\omega_{T-1,k}}(y_k^{(T-1)})\cdots\nabla f^{(1)}_{\omega_{1,k}}(y_k^{(1)}) \\
& + \alpha_k^2\|\widetilde{\nabla} f^{(T)}_{\omega_{T,k}}(x_k)\nabla f^{(T-1)}_{\omega_{T-1,k}}(y_k^{(T-1)})\cdots\widetilde{\nabla} f^{(1)}_{\omega_{1,k}}(y_k^{(1)})\|^2 \\
= & \|x_k - x^*\|^2 - 2\alpha_k(x_k - x^*)'\widetilde{\nabla} F_{\omega_k}(x_k) + u_k \\
& + \alpha_k^2\|\widetilde{\nabla} f^{(T)}_{\omega_{T,k}}(x_k)\nabla f^{(T-1)}_{\omega_{T-1,k}}(y_k^{(T-1)})\cdots\nabla f^{(1)}_{\omega_{1,k}}(y_k^{(1)})\|^2,
\end{aligned} \tag{A.7}
$$

where

$$
\widetilde{\nabla} F_{\omega_k}(x_k) = \widetilde{\nabla} f^{(T)}_{\omega_{T,k}}(x_k)\nabla f^{(T-1)}_{\omega_{T-1,k}}(f^{(T)}(x_k))\cdots\nabla f^{(1)}_{\omega_{1,k}}(f^{(2)}\circ\cdots\circ f^{(T)}(x_k)),
$$

as defined in the main text, and

$$
u_k = 2\alpha_k(x_k - x^*)'\Big[\widetilde{\nabla} F_{\omega_k}(x_k) - \widetilde{\nabla} f^{(T)}_{\omega_{T,k}}(x_k)\nabla f^{(T-1)}_{\omega_{T-1,k}}(y_k^{(T-1)})\cdots\nabla f^{(1)}_{\omega_{1,k}}(y_k^{(1)})\Big].
$$

By Assumption 2.1 (iii)-(iv), $\widetilde{\nabla} f^{(j)}_{\omega_{j,k}}$'s have bounded second-order moments, thus with probability 1,

$$
\mathbb{E}\Big[\|\widetilde{\nabla} f^{(T)}_{\omega_{T,k}}(x_k)\nabla f^{(T-1)}_{\omega_{T-1,k}}(y_k^{(T-1)})\cdots\nabla f^{(1)}_{\omega_{1,k}}(y_k^{(1)})\|^2\Big|\mathbb{F}_k\Big] \le C_1\cdots C_T. \tag{A.8}
$$

Taking expectation on both sides of Eq.(A.7), conditioning on $\mathbb{F}_k$, we have

$$
\begin{aligned}
& \mathbb{E}[\|x_{k+1} - x^*\|^2|\mathbb{F}_k] \\
\le & \|x_k - x^*\|^2 + \alpha_k^2 C_1C_2\cdots C_T + \mathbb{E}[u_k|\mathbb{F}_k] - 2\alpha_k(x_k - x^*)'\mathbb{E}\Big[\widetilde{\nabla} F_{\omega_k}(x_k)\Big|\mathbb{F}_k\Big].
\end{aligned}
$$

By the convexity of $F(x)$ and Assumption 2.1(ii), we obtain

$$
(x_k - x^*)'\mathbb{E}\Big[\widetilde{\nabla} F_{\omega_k}(x_k)\Big|\mathbb{F}_k\Big] \ge F(x_k) - F^*.
$$

Then we have

$$
\begin{aligned}
& \mathbb{E}[\|x_{k+1} - x^*\|^2|\mathbb{F}_k] \\
\le & \|x_k - x^*\|^2 + \alpha_k^2 C_1C_2\cdots C_T - 2\alpha_k(F(x_k) - F^*) + \mathbb{E}[u_k|\mathbb{F}_k].
\end{aligned} \tag{A.9}
$$

For the term $u_k$, there exists a constant $C_0 > 0$ such that with probability 1,

$$\begin{aligned}
&\mathbb{E}[u_k|\mathbb{F}_k] \\
=&2\alpha_k(x_k - x^*)'\mathbb{E}\Big[\widetilde{\nabla} F_{\omega_k}(x_k) - \widetilde{\nabla} f^{(T)}_{\omega_{T,k}}(x_k)\nabla f^{(T-1)}_{\omega_{T-1,k}}(y_k^{(T-1)})\cdots\nabla f^{(1)}_{\omega_{1,k}}(y_k^{(1)})\Big|\mathbb{F}_k\Big] \\
\le&(T-1)\beta_{T-1,k}\mathbb{E}[\|y_k^{(T-1)} - f^{(T)}(x_k)\|^2|\mathbb{F}_k] + (T-2)\beta_{T-2,k}\mathbb{E}[\|y_k^{(T-2)} - f^{(T-1)}(y_k^{(T-1)})\|^2|\mathbb{F}_k] \\
&+\cdots+\beta_{1,k}\mathbb{E}[\|y_k^{(1)} - f^{(2)}(y_k^{(2)})\|^2|\mathbb{F}_k] + C_0(\frac{\alpha_k^2}{\beta_{T-1,k}}+\cdots+\frac{\alpha_k^2}{\beta_{1,k}})\|x_k - x^*\|^2,
\end{aligned} \tag{A.10}$$

where the last inequality comes from Lemma A.1 by letting $X_k = 2\alpha_k(x_k - x^*) \in \mathbb{F}_k$. Substituting Eq.(A.10) into Eq.(A.9), we get

$$\begin{aligned}
&\mathbb{E}[\|x_{k+1} - x^*\|^2|\mathbb{F}_k] \\
\le&\Big(1 + C_0(\frac{\alpha_k^2}{\beta_{T-1,k}} + \frac{\alpha_k^2}{\beta_{T-2,k}} + \cdots + \frac{\alpha_k^2}{\beta_{1,k}})\Big)\|x_k - x^*\|^2 + \alpha_k^2 C_1 C_2\cdots C_T - 2\alpha_k(F(x_k) - F^*) \\
&+ (T-1)\beta_{T-1,k}\mathbb{E}[\|y_k^{(T-1)} - f^{(T)}(x_k)\|^2|\mathbb{F}_k] + (T-2)\beta_{T-2,k}\mathbb{E}[\|y_k^{(T-2)} - f^{(T-1)}(y_k^{(T-1)})\|^2|\mathbb{F}_k] \\
&+\cdots+\beta_{1,k}\mathbb{E}[\|y_k^{(1)} - f^{(2)}(y_k^{(2)})\|^2|\mathbb{F}_k],
\end{aligned}$$

which completes the proof. □

## A.2 Proof of Lemma 2.2

Proof of lemma 2.2: By the assumptions in part (b), $F$ has Lipschitz continuous gradient with parameter $L_F$ and $\mathcal{X} = \mathbb{R}^{d_T}$, we denote by $\nabla F(x)$ as the gradient of $F(x)$, and obtain

$$\begin{aligned}
&F(x_{k+1}) - F(x_k) \\
\le&\langle\nabla F(x_k), x_{k+1} - x_k\rangle + \frac{L_F}{2}\|x_{k+1} - x_k\|^2 \\
=&-\alpha_k\langle\nabla F(x_k), \widetilde{\nabla} f^{(T)}_{\omega_{T,k}}(x_k)\nabla f^{(T-1)}_{\omega_{T-1,k}}(y_k^{(T-1)})\cdots\nabla f^{(1)}_{\omega_{1,k}}(y_k^{(1)})\rangle + \frac{L_F}{2}\|x_{k+1} - x_k\|^2 \\
=&-\alpha_k\|\nabla F(x_k)\|^2 + \alpha_k\nabla F(x_k)'\Big[\nabla F(x_k) - \widetilde{\nabla} f^{(T)}_{\omega_{T,k}}(x_k)\nabla f^{(T-1)}_{\omega_{T-1,k}}(y_k^{(T-1)})\cdots\nabla f^{(1)}_{\omega_{1,k}}(y_k^{(1)})\Big] \\
&+\frac{L_F}{2}\|x_{k+1} - x_k\|^2.
\end{aligned} \tag{A.11}$$

As defined in the main text, $\widetilde{\nabla} F_{\omega_k}(x_k) = \widetilde{\nabla} f^{(T)}_{\omega_{T,k}}(x_k)\nabla f^{(T-1)}_{\omega_{T-1,k}}(f^{(T)}(x_k))\cdots\nabla f^{(1)}_{\omega_{1,k}}(f^{(2)}\circ\cdots\circ f^{(T)}(x_k))$, by Assumption 2.1 (ii), we have

$$\mathbb{E}[\widetilde{\nabla} F_{\omega_k}(x_k)|\mathbb{F}_k] = \nabla F(x_k).$$

We obtain that with probability 1, there exists a constant $C_0 > 0$ such that

$$\begin{aligned}
&\alpha_k\mathbb{E}\Big[\nabla F(x_k)'\Big(\nabla F(x_k) - \widetilde{\nabla} f^{(T)}_{\omega_{T,k}}(x_k)\nabla f^{(T-1)}_{\omega_{T-1,k}}(y_k^{(T-1)})\cdots\nabla f^{(1)}_{\omega_{1,k}}(y_k^{(1)})\Big)\Big|\mathbb{F}_k\Big] \\
=&\alpha_k\nabla F(x_k)'\mathbb{E}\Big[\widetilde{\nabla} F_{\omega_k}(x_k) - \widetilde{\nabla} f^{(T)}_{\omega_{T,k}}(x_k)\nabla f^{(T-1)}_{\omega_{T-1,k}}(y_k^{(T-1)})\cdots\nabla f^{(1)}_{\omega_{1,k}}(y_k^{(1)})\Big|\mathbb{F}_k\Big] \\
\le&(T-1)\beta_{T-1,k}\mathbb{E}[\|y_{k+1}^{(T-1)} - f^{(T)}(x_k)\|^2|\mathbb{F}_k] + (T-2)\beta_{T-2,k}\mathbb{E}[\|y_{k+1}^{(T-2)} - f^{(T-1)}(y_{k+1}^{(T-1)})\|^2|\mathbb{F}_k] \\
&+\cdots+\beta_{1,k}\mathbb{E}[\|y_{k+1}^{(1)} - f^{(2)}(y_{k+1}^{(2)})\|^2|\mathbb{F}_k] + C_0\Big(\frac{\alpha_k^2}{\beta_{T-1,k}}+\cdots+\frac{\alpha_k^2}{\beta_{1,k}}\Big)\|\nabla F(x_k)\|^2,
\end{aligned}$$

where the last inequality comes from Lemma A.1 by letting $X_k = \alpha_k \nabla F(x_k) \in \mathbb{F}_k$. Also note that

$$\mathbb{E}[\|x_{k+1} - x_k\|^2|\mathbb{F}_k] = \alpha_k^2 \mathbb{E}[\|\widetilde{\nabla} f^{(T)}_{\omega_{T,k}}(x_k)\nabla f^{(T-1)}_{\omega_{T-1,k}}(y_k^{(T-1)})\cdots\nabla f^{(1)}_{\omega_{1,k}}(y_k^{(1)})\|^2|\mathbb{F}_k] \leq \alpha_k^2 C_1\cdots C_T. \tag{A.12}$$

Combining the results above, we obtain

$$\begin{aligned}
&\mathbb{E}[F(x_{k+1}) - F^*|\mathbb{F}_k] \\
\leq& F(x_k) - F^* + \frac{1}{2}\alpha_k^2 L_F C_1 C_2 \cdots C_T - \alpha_k\Big(1 - (\frac{\alpha_k}{\beta_{T-1,k}} - \cdots - \frac{\alpha_k}{\beta_{1,k}})C_0\Big)\|\nabla F(x_k)\|^2 \\
&+ (T-1)\beta_{T-1,k}\mathbb{E}[\|y_k^{(T-1)} - f^{(T)}(x_k)\|^2|\mathbb{F}_k] + (T-2)\beta_{T-2,k}\mathbb{E}[\|y_k^{(T-2)} - f^{(T-1)}(y_k^{(T-1)})\|^2|\mathbb{F}_k] \\
&+ \cdots + \beta_{1,k}\mathbb{E}[\|y_k^{(1)} - f^{(2)}(y_k^{(2)})\|^2|\mathbb{F}_k].
\end{aligned}$$

Besides, we have $1/2 \leq 1 - (\frac{\alpha_k}{\beta_{T-1,k}} - \cdots - \frac{\alpha_k}{\beta_{T-1,k}})C_0$ for $k$ sufficiently large since $\alpha_k/\beta_{j,k} \to 0$ as $k \to \infty$ for all $j$'s. Finally, we conclude

$$\begin{aligned}
&\mathbb{E}[F(x_{k+1}) - F^*|\mathbb{F}_k] \\
\leq& F(x_k) - F^* - \frac{\alpha_k}{2}\|\nabla F(x_k)\|^2 + \frac{1}{2}\alpha_k^2 L_F C_1 C_2\cdots C_T + (T-1)\beta_{T-1,k}\mathbb{E}[\|y_k^{(T-1)} - f^{(T)}(x_k)\|^2|\mathbb{F}_k] \\
&+ \cdots + \beta_{1,k}\mathbb{E}[\|y_k^{(1)} - f^{(2)}(y_k^{(2)})\|^2|\mathbb{F}_k],
\end{aligned}$$

for $k$ sufficiently large, which completes our proof. □

### A.3 Proof of Lemma 2.3

Proof of Lemma 2.3: a) For ease of presentation, we denote $y_{k+1}^{(T-1)}$ by $z_{k+1}$, $\beta_{T-1,k}$ by $\gamma_k$, and $f^{(T)}_{\omega_{T,k+1}}(x_k)$ by $h_{u_{k+1}}(x_k)$. The corresponding update step can be written as

$$z_{k+1} = (1-\gamma_k)z_k + \gamma_k h_{u_{k+1}}(x_{k+1}).$$

Let $e_k = (1-\gamma_k)(h(x_{k+1}) - h(x_k))$. Together with the definition of $z_{k+1}$ above, we have

$$z_{k+1} - h(x_{k+1}) + e_k = (1-\gamma_k)\big(z_k - h(x_k)\big) + \gamma_k\big(h_{u_{k+1}}(x_{k+1}) - h(x_{k+1})\big). \tag{A.13}$$

In addition, by the Lipschitz continuity of $h$ in Assumption 2.1 (iii), we obtain

$$\|e_k\| \leq (1-\gamma_k)\sqrt{C_h}\|x_{k+1} - x_k\|. \tag{A.14}$$

Taking expectation of the squared norm on both sides of Eq.(A.13) conditioned on $\mathbb{F}_k$, and using Assumption 2.1 (ii)-(iii), we have

$$\begin{aligned}
&\mathbb{E}[\|z_{k+1} - h(x_{k+1}) + e_k\|^2|\mathbb{F}_k] \\
=&\mathbb{E}\Big[\|(1-\gamma_k)(z_k - h(x_k)) + \gamma_k(h_{u_{k+1}}(x_{k+1}) - h(x_{k+1}))\|^2\Big|\mathbb{F}_k\Big] \\
=&(1-\gamma_k)^2\mathbb{E}[\|z_k - h(x_k)\|^2|\mathbb{F}_k] + \gamma_k^2\mathbb{E}[\|h_{u_{k+1}}(x_{k+1}) - h(x_{k+1})\|^2|\mathbb{F}_k] \\
&+ 2(1-\gamma_k)\gamma_k\mathbb{E}[(z_k - h(x_k))^T(h_{u_{k+1}}(x_{k+1}) - h(x_{k+1}))|\mathbb{F}_k] \\
\leq&(1-\gamma_k)^2\mathbb{E}[\|z_k - h(x_k)\|^2|\mathbb{F}_k] + \gamma_k^2 V_h + 0.
\end{aligned} \tag{A.15}$$

Using the fact that $\|a+b\|^2 \le (1+\epsilon)\|a\|^2 + (1+1/\epsilon)\|b\|^2$ for any $\epsilon > 0$, we obtain

$$\|z_{k+1} - h(x_{k+1})\|^2 \le (1+\gamma_k)\|z_{k+1} - h(x_{k+1}) + e_k\|^2 + (1+1/\gamma_k)\|e_k\|^2. \tag{A.16}$$

Taking expectation of both sides conditioned on $\mathbb{F}_k$, and applying Eqs.(A.14)-(A.16), we conclude

$$\begin{aligned}
&\mathbb{E}[\|z_{k+1} - h(x_{k+1})\|^2|\mathbb{F}_k] \\
\le&(1+\gamma_k)(1-\gamma_k^2)\mathbb{E}[\|z_k - h(x_k)\|^2|\mathbb{F}_k] + (1+\gamma_k)\gamma_k^2 V_h + \frac{(1-\gamma_k^2)C_h}{\gamma_k}\mathbb{E}[\|x_{k+1} - x_k\|^2|\mathbb{F}_k] \\
\le&(1-\gamma_k)\mathbb{E}[\|z_k - h(x_k)\|^2|\mathbb{F}_k] + \frac{C_h}{\gamma_k}\mathbb{E}[\|x_{k+1} - x_k\|^2|\mathbb{F}_k] + 2V_h\gamma_k^2,
\end{aligned} \tag{A.17}$$

for all $k$, with probability 1.

b) Based on our assumption, we have $\mathbb{E}[\|x_{k+1} - x_k\|^2] \le \mathcal{O}(\alpha_k^2)$ and $\sum_{k=0}^{\infty} \alpha_k^2 \gamma_k^{-1} < \infty$, then we obtain

$$\sum_{k=0}^{\infty} \gamma_k^{-1}\mathbb{E}\Big[\mathbb{E}[\|x_{k+1} - x_k\|^2|\mathbb{F}_k]\Big] = \sum_{k=0}^{\infty} \gamma_k^{-1}\mathbb{E}[\|x_{k+1} - x_k\|^2] \le \sum_{k=0}^{\infty} \mathcal{O}(\frac{\alpha_k^2}{\gamma_k}) < \infty.$$

By the monotone convergence theorem, we conclude that $\sum_{k=1}^{n} \gamma_k^{-1} C_h \mathbb{E}[\|x_{k+1} - x_k\|^2|\mathbb{F}_k]$ converges almost surely to a random variable with finite expectation as $n \to \infty$. Therefore, the limit $\sum_{k=0}^{\infty} \gamma_k^{-1} C_h \mathbb{E}[\|x_{k+1} - x_k\|^2|\mathbb{F}_k]$ exists and is finite with probability 1.

c) By our assumption that $\mathbb{E}[\|x_{k+1} - x_k\|^2] \le \mathcal{O}(\alpha_k^2)$, there exists a constant $C_0 > 0$ such that $\mathbb{E}[\|x_{k+1} - x_k\|^2] \le \alpha_k^2 C_0$. Since $\alpha_k/\gamma_k \to 0$, there exists a constant $M > 0$ such that $\alpha_k \le M\gamma_k$ for all $k$. Let $D_z = \mathbb{E}[\|z_0 - h(x_0)]\|^2 + 2V_h + M^2 C_0 C_h$. Since $\alpha_k \le M\gamma_k$ and $\gamma_k \le 1$, we have $D_z \ge 2V_h\gamma_k + \gamma_k^{-2}\alpha_k^2 C_0 C_h$ for all $k$.

We prove by induction that $\mathbb{E}[\|z_{k+1} - h(x_k)\|^2] \le D_z$ for all $k$. Clearly, the claim holds for $k = 0$. Suppose the claim holds for $0, 1, ..., k-1$. By Eq.(A.17), we have that there exists a constant $C_0 > 0$ such that

$$\begin{aligned}
\mathbb{E}[\|z_{k+1} - h(x_{k+1})\|^2] \le&(1-\gamma_k)\mathbb{E}[\|z_k - h(x_k)\|^2] + 2V_h\gamma_k^2 + \gamma_k^{-1}\alpha_k^2 C_0 C_h \\
\le&(1-\gamma_k)D_z + 2V_h\gamma_k^2 + \gamma_k^{-1}\alpha_k^2 C_0 C_h \\
=&D_z - \gamma_k(D_z - 2V_h\gamma_k + \gamma_k^{-2}\alpha_k^2 C_0 C_h) \\
\le&D_z,
\end{aligned}$$

where the second inequality uses the fact $1 - \gamma_k \ge 0$ and $\mathbb{E}[\|z_k - h(x_k)\|^2] \le D_z$. Our claim holds as desired.

d) By the definition of $z_{k+1}$, we have $z_{k+1} = (1-\gamma_k)z_k + \gamma_k h_{u_{k+1}}(x_{k+1})$ and

$$\begin{aligned}
(1-\gamma_k)^2\|z_{k+1} - z_k\|^2 =& \gamma_k^2\|h_{u_{k+1}}(x_{k+1}) - z_{k+1}\|^2 \\
\le& 2\gamma_k^2\|h_{u_{k+1}}(x_{k+1}) - h(x_{k+1})\|^2 + 2\gamma_k^2\|z_{k+1} - h(x_{k+1})\|^2.
\end{aligned}$$

Then we obtain

$$\mathbb{E}[\|z_{k+1} - z_k\|^2] \le \frac{2\gamma_k^2}{(1-\gamma_k)^2}V_h + \frac{2\gamma_k^2}{(1-\gamma_k)^2}\mathbb{E}[\|z_{k+1} - h(x_{k+1})\|^2]. \tag{A.18}$$

From part (c), we have that there exists $D_z \ge 0$ such that $\mathbb{E}[\|z_{k+1} - h(x_{k+1})\|^2] \le D_z$. Plug this into Eq.(A.18), we conclude

$$\mathbb{E}[\|z_{k+1} - z_k\|^2] \le \mathcal{O}(\gamma_k^2).$$

□

### A.4 Proof of Lemma 2.4

This lemma proves the $T$-element super-martingale convergent lemma to establish convergence property of $\{x_{k+1} - x^*\}$.

Proof: Let $J_k$ be the random variable

$$J_k \equiv X_k + c_{T-1}Y_k^{(T-1)} + c_{T-2}Y_k^{(T-2)} + \cdots + c_1 Y_k^{(1)}.$$

We have

$$\begin{aligned}\mathbb{E}[J_{k+1}|\mathbb{F}_k] =& \mathbb{E}[X_{k+1}|\mathbb{F}_k] + c_{T-1}\mathbb{E}[Y_{k+1}^{(T-1)}|\mathbb{F}_k] + \cdots + c_1\mathbb{E}[Y_{k+1}^{(1)}|\mathbb{F}_k] \\ \leq& (1+\eta_k)X_k + c_{T-1}Y_k^{(T-1)} + \cdots + c_1 Y_k^{(1)} + \mu_k^{(T)} + c_{T-1}\mu_k^{(T-1)} + \cdots + c_1\mu_k^{(1)} \\ \leq& (1+\eta_k)J_k + \mu_k^{(T)} + c_{T-1}\mu_k^{(T-1)} + \cdots + c_1\mu_k^{(1)}.\end{aligned}$$

Since $\sum_{k=0}^{\infty} \mu_k^{(T)} + c_{T-1}\mu_k^{(T-1)} + \cdots + c_1\mu_k^{(1)} < \infty$ and $\sum_{k=0}^{\infty} \eta_k < \infty$, we obtain that $J_k$ converges almost surely to a random variable by Theorem 1 in [1], and $J_k$ is bounded by a constant with probability 1.

By the definition above, we have that $X_k \leq J_k$, $Y_k^{(2)} \leq \frac{1}{c_2}J_k, \cdots, Y_k^{(T)} \leq \frac{1}{c_T}J_k$. Then $X_k, Y_k^{(T-1)}, Y_k^{(T-2)}, \cdots, Y_k^{(1)}$ are also bounded with probability 1. Since

$$\sum_{k=0}^{\infty} \mu_k^{(j)} < \infty, \quad \text{for } j = 1, \cdots, T,$$

and $\theta_k^{(T-1)}, \cdots, \theta_k^{(1)}$ are nonnegative, we have

$$\mathbb{E}[Y_{k+1}^{(T-1)}|\mathbb{F}_k] \leq (1-\theta_k^{(j)})Y_k^{(j)} + \mu_k^{(j)}, \quad \text{for } j = T-1, \cdots, 1.$$

Again, by Theorem 1 in [1], we obtain that $Y_k^{(T-1)}, \cdots, Y_k^{(1)}$ converge almost surely to $T-1$ random variables, respectively, and

$$\sum_{k=1}^{\infty} \theta_k Y_k^{(j)} \leq \infty \;\; w.p.1 \text{ for } j = 1, \cdots, T-1.$$

Since $Y_k^{(1)}, \cdots, Y_k^{(T-1)}$ and $J_k = X_k + c_{T-1}Y_k^{(T-1)} + c_{T-2}Y_k^{(T-2)} + \cdots + c_1Y_k^{(1)}$ are almost surely convergent, then $X_k$ must converge almost surely to a random variable.

Together with the condition that

$$\sum_{k=0}^{\infty} \eta_k < \infty, \sum_{k=0}^{\infty} \mu_k^{(j)} < \infty, \quad \text{for } j = 1, \cdots, T,$$

we have

$$\sum_{j=1}^{T}\sum_{k=0}^{\infty} u_k^{(j)} < \infty \;\; w.p.1.$$

So we complete the proof. $\square$

### A.5 Proof of Lemma 2.5

Note that Lemma 2.5 and 2.6 are proved in Theorem 1 in [2], we present the details here for self-completeness.

Proof: Let $\Omega_{x^*}$ be the collection of sample paths that $\Omega_{x^*} = \{\omega : \lim_{k\to\infty} \|x_k(\omega) - x^*\|\}$ exists. It has been proved that $\mathbb{P}(\Omega_{x^*}) = 1$ for any $x^* \in \mathcal{X}^*$. We claim that $\cap_{x^*\in\mathcal{X}^*}\Omega_{x^*}$ is measurable and $\mathbb{P}(\cap_{x^*\in\mathcal{X}^*}\Omega_{x^*}) = 1$.

To see this, we consider a countable dense subset $\mathcal{X}_Q^*$ of $\mathcal{X}^*$. By the fact that $F$ is convex, we have that the set $\mathcal{X}^* \subset \mathbb{R}^n$ is separable, and such $\mathcal{X}_Q^*$ exists. So the probability of non-convergence for some $x^* \in \mathcal{X}_Q^*$ is the probability of a union of countably many sets, each having probability 0. As a result, we have

$$\mathbb{P}(\cap_{\mathcal{X}_Q^*}\Omega_{x^*}) = 1 - \mathbb{P}(\cup_{\mathcal{X}_Q^*}\Omega_{x^*}^c) \geq 1 - \sum_{x^*\in\mathcal{X}_Q^*} \mathbb{P}(\Omega_{x^*}^c) = 1.$$

Now we consider any arbitrary $\widetilde{x} \in \mathcal{X}^*$, which is the limit of a sequence of optimal solutions $\{\widetilde{x}_j\}_{j=1}^\infty \subset \mathcal{X}_Q^*$. By using a limit point argument, it is not hard to see that $\|x_k(\omega) - \widetilde{x}\|$ is convergent as long as $\|x_k(\omega) - \widetilde{x}_j\|$ is convergent for all $j$. Consider arbitrary $\omega \in \cap_{\mathcal{X}_Q^*}\Omega_{x^*}$. We have the triangle inequalities:

$$\|x_k(\omega) = \widetilde{x}_j\| - \|\widetilde{x}_j - \widetilde{x}\| \leq \|x_k(\omega) - \widetilde{x}_j\| + \|\widetilde{x}_j - \widetilde{x}\|.$$

By taking $k \to \infty$ and using the fact that $\lim_{k\to\infty} \|x_k(\omega) - \widetilde{x}_j\|$ exists, we obtain

$$\begin{aligned}
&\lim_{k\to\infty} \|x_k(\omega) - \widetilde{x}_j\| - \|\widetilde{x}_j - \widetilde{x}\| \leq \liminf_{k\to\infty} \|x_k(\omega) - \widetilde{x}\| \\
&\leq \limsup_{k\to\infty} \|x_k(\omega) - \widetilde{x}\| \leq \lim_{k\to\infty} \|x_k(\omega) - \widetilde{x}_j\| + \|\widetilde{x}_j - \widetilde{x}\|.
\end{aligned}$$

So we have

$$\limsup_{k\to\infty} \|x_k(\omega) - \widetilde{x}\| - \liminf_{k\to\infty} \|x_k(\omega) - \widetilde{x}_j\| \leq 2\|\widetilde{x}_j - \widetilde{x}\|.$$

By taking $j \to \infty$, we have $\|\widetilde{x}_j - \widetilde{x}\| \to 0$ and

$$\limsup_{k\to\infty} \|x_k(\omega) - \widetilde{x}\| = \liminf_{k\to\infty} \|x_k(\omega) - \widetilde{x}_j\|$$

It follows that the limit of $\|x_j(\omega) - \widetilde{x}\|$ exists, therefore $\omega \in \Omega_{\widetilde{x}}$. Now we have proved that $\cap_{\mathcal{X}_Q^*}\Omega_{x^*} \subset \Omega_{\widetilde{x}}$ for all $\widetilde{x} \in \mathcal{X}^*$. As a result, we have $\cap_{\mathcal{X}_Q^*}\Omega_{x^*} \subset \cap_{\mathcal{X}^*}\Omega_{x^*}$. Note that $\mathbb{P}(\cap_{\mathcal{X}_Q^*}\Omega_{x^*}) = 1$ and that the considered probability measure is complete. We have that $(\cap_{\mathcal{X}^*}\Omega_{x^*})^c$ is the subset of a null set $(\cap_{\mathcal{X}_Q^*}\Omega_{x^*})^c$, therefore it is measurable and satisfies $\mathbb{P}(\cap_{\mathcal{X}^*}\Omega_{x^*}) \leq \mathbb{P}(\cap_{\mathcal{X}_Q^*}\Omega_{x^*}) = 0$. It follows that $\cap_{\mathcal{X}^*}\Omega_{x^*}$ is measurable and satisfied $\mathbb{P}(\cap_{\mathcal{X}^*}\Omega_{x^*}) = 1$. In other words, $\|x_k - \widetilde{x}\|$ is convergent for all $\widetilde{x} \in \mathcal{X}^*$, with probability 1.

Consider an arbitrary sample trajectory $\{x_k(\omega)\}$ such that $\omega \in \cap_{\mathcal{X}^*}\Omega_{x^*}$ and $\liminf_{k\to\infty} F(x_k(\omega)) = F^*$. Take arbitrary $x^* \in \mathcal{X}^*$,since $\{\|x_k(\omega) - x^*\|\}$ converges, the sequence is bounded. By the continuity of $F$, the sequence $\{x_k(\omega)\}$ must have a limit point $\bar{x}$ being an optimal solution, i.e., $F(\bar{x}) = F^*$ and $\bar{x} \in \mathcal{X}^*$.

Since $\omega \in \cap_{\mathcal{X}^*}\Omega_{x^*} \subset \Omega_{\bar{x}}$, we obtain that $\{\|x_k(\omega) - \bar{x}\|\}$ is also a convergent sequence. Since it is convergent while having a limit point 0, we have $\|x_k(\omega) - \bar{x}\| \to 0$.Then $x_k(\omega) \to \bar{x}$ on this sample

trajectory. Note that $\bar{x}$ depends on the sample path $\omega$, so it is a random variable. Also note that the set of all such sample paths has a probability measure equaling to 1. Therefore $x_k$ converges almost surely to a random point in the set of optimal solutions to problem (1.1). □

### A.6 Proof of Lemma 2.6

Proof: We focus on a single sample trajectory $x_k(\omega)$ such that the preceding inequalities hold. For simplicity, we omit the notation $(\omega)$ in the rest of the proof.

Let $\epsilon > 0$ be arbitrary. We note that $\|\nabla F(x_k)\| \leq \epsilon$ holds for infinitely many $k$. Otherwise we would have for some $\bar{k} > 0$ that $\sum_{k=0}^{\infty} \alpha_k \|\nabla F(x_k)\|^2 \geq \sum_{k=\bar{k}}^{\infty} \alpha_k \epsilon^2 = \infty$, yielding a contradiction. Consequently, there exists a closed set $\bar{N}$ (e.g., closed union of neighborhoods of all $\epsilon$- stationary $x_k$'s) such that $\{x_k\}$ visits infinitely often, and

$$\|\nabla F(x)\| = \begin{cases} \leq \epsilon & \text{if } x \in \bar{N}, \\ > \epsilon & \text{if } x \notin \bar{N}, x \in \{x_k\}. \end{cases}$$

We assume to the contrary that there exists a limit point $\widetilde{x}$ such that $\|\nabla F(\widetilde{x})\| > 2\epsilon$. Then there exists a closed set $\widetilde{N}$ (e.g., union of neighborhoods of all $x_k$?s such that $\|\nabla F(x_k)\| > 2\epsilon$) such that $\{x_k\}$ visits infinitely often, and

$$\|\nabla F(x)\| = \begin{cases} \leq 2\epsilon & \text{if } x \in \widetilde{N}, \\ > 2\epsilon & \text{if } x \notin \widetilde{N}, x \in \{x_k\}. \end{cases}$$

By the continuity of $\nabla F$ and $\epsilon > 0$, we obtain that the sets $\widetilde{N}$ and $\bar{N}$ are disjoint, i.e., $\text{dist}(\bar{N}, \widetilde{N}) > 0$. Since the sequence $x_k$ enters both sets $\bar{N}$ and $\widetilde{N}$ infinitely often, there exists a subsequence

$$\{x_k\}_{k \in \mathcal{K}} = \left\{ \{x_k\}_{k=s_i}^{t_i - 1} \right\}_{i=1}^{\infty}$$

that crosses the two sets infinitely often, with $x_{s_i} \in \widetilde{N}$ , $x_{t_i} \in \widetilde{N}$ for all $i$. In other words, we have for every $i$ that

$$\|\nabla F(x_{s_i})\| \geq 2\epsilon > \|\nabla F(x_k)\| > \epsilon \geq \|\nabla F(x_{t_i})\|, \quad \forall k = s_i + 1, \cdots, t_i - 1.$$

On one hand, by using the triangle inequality, we have

$$\begin{aligned} \sum_{k \in \mathcal{K}} \|x_{k+1} - x_k\| &= \sum_{k=1}^{\infty} \sum_{k=s_i}^{t_i - 1} \|x_{k+1} - x_k\| \geq \sum_{i=1}^{\infty} \|x_{t_i} - x_{s_i}\| \\ &\geq \sum_{i=1}^{\infty} \text{dist}(\bar{N}, \widetilde{N}) = \infty. \end{aligned}$$

On the other hand, we have

$$\infty > \sum_{k=0}^{\infty} \alpha_k \|\nabla F(x_k)\|^2 \geq \sum_{k \in \mathcal{K}} \alpha_k \|\nabla F(x_k)\|^2 > \epsilon \sum_{k \in \mathcal{K}} \alpha_k.$$

By using the uniform boundedness of random variables generated from the $SO$, we have for some $M > 0$ that $\|x_{k+1} - x_k\| \leq \|\widetilde{\nabla} f^{(1)}_{\omega_{1,k}}(x_k) \nabla f^{(2)}_{\omega_{2,k}}(y_k^{(T-1)}) \cdots \nabla f_{\omega_{T,k}}(y_k^{(1)})\| \leq M\alpha_k$ for all $k$. It follows that

$$\sum_{k\in\mathcal{K}} \|x_{k+1} - x_k\| \leq M \sum_{k\in\mathcal{K}} \alpha_k < \infty.$$

This yields a contradiction with the result $\sum_{k\in\mathcal{K}} \|x_{k+1} - x_k\| = \infty$ we just obtained. It follows that there does not exist a limit point $\widetilde{x}$ such that $\|\nabla F(\widetilde{x})\| > 2\epsilon$. Since $\epsilon$ can be made arbitrarily small, there does not exist any limit point that is nonstationary. Finally, we note that the set of such sample paths (to which the preceding analysis applies) has a probability measure 1. In other words, any limit point of $x_k$ is a stationary point of $F(x)$ with probability 1. □

### A.7 Proof of Theorem 2.1

Proof: a) Let $x^*$ be an arbitrary optimal solution to problem (1.1), and let $F^* = F(x^*)$. By Lemma 2.1 , we have with probability 1,

$$\begin{aligned}
&\mathbb{E}[\|x_{k+1} - x^*\|^2 | \mathbb{F}_k] \\
\leq &\Big(1 + C_0\big(\frac{\alpha_k^2}{\beta_{T-1,k}} + \frac{\alpha_k^2}{\beta_{T-2,k}} + \cdots + \frac{\alpha_k^2}{\beta_{1,k}}\big)\Big)\|x_k - x^*\|^2 + \alpha_k^2 C_1 C_2 \cdots C_T - 2\alpha_k\big(F(x_k) - F^*\big) \\
&+ (T-1)\beta_{T-1,k}\mathbb{E}[\|y_k^{(T-1)} - f^{(T)}(x_k)\|^2|\mathbb{F}_k] + (T-2)\beta_{T-2,k}\mathbb{E}[\|y_k^{(T-2)} - f^{(T-1)}(y_k^{(T-1)})\|^2|\mathbb{F}_k] \\
&+ \cdots + \beta_{1,k}\mathbb{E}[\|y_k^{(1)} - f^{(2)}(y_k^{(2)})\|^2|\mathbb{F}_k].
\end{aligned} \tag{A.19}$$

By Lemma 2.3, for any basic update step $j = T-1, \cdots, 1$, if $\mathbb{E}[\|y_{k+1}^{(j+1)} - y_k^{(j+1)}\|^2] \leq \mathcal{O}(\beta^2_{j+1,k})$, then

$$\begin{aligned}
&\mathbb{E}\Big[\mathbb{E}[\|y_{k+1}^{(j)} - f^{(j+1)}(y_{k+1}^{(j+1)})\|^2|\mathbb{F}_{k+1}]\Big|\mathbb{F}_k\Big] = \mathbb{E}[\|y_{k+1}^{(j)} - f^{(j+1)}(y_{k+1}^{(j+1)})\|^2|\mathbb{F}_k] \\
\leq &(1-\beta_{j,k})\mathbb{E}[\|y_k^{(j)} - f^{(j+1)}(y_k^{(j+1)})\|^2|\mathbb{F}_k] + \beta_{j,k}^{-1}C_{j+1}\mathbb{E}[\|y_{k+1}^{(j+1)} - y_k^{(j+1)}\|^2|\mathbb{F}_k] + 2V_{j+1}\beta^2_{j,k},
\end{aligned} \tag{A.20}$$

and $\mathbb{E}[\|y_{k+1}^{(j)} - y_k^{(j)}\|^2] \leq \mathcal{O}(\beta^2_{j,k})$. By Eq.(A.12), we have $\mathbb{E}[\|x_{k+1} - x_k\|^2] \leq \mathcal{O}(\alpha_k^2)$, which serves as the sufficient condition for level $T-1$ for Eq.(A.20) to be true so that $\mathbb{E}[\|y_{k+1}^{(T-1)} - y_k^{(T-1)}\|^2] \leq \mathcal{O}(\beta^2_{T-1,k})$, a sufficient condition for level $T-2$ for Eq.(A.20). Conducting this chain to all the levels, we have that Eq.(A.20) holds for $j = T-1, \cdots, 1$.

Besides, by Lemma 2.3 (b), under the condition $\sum_{k=1}^{\infty} \beta^2_{j+1,k}/\beta_{j,k} < \infty$, we have

$$\beta_{j,k}^{-1}C_{j+1}\mathbb{E}[\|y_{k+1}^{(j+1)} - y_k^{(j+1)}\|^2|\mathbb{F}_k] + 2V_{j+1}\beta^2_{j,k} < \infty,$$

with probability 1.

Now we apply the $T$-element super-martingale convergent lemma 2.4 to Eqs.(A.19)-(A.20). By

letting

$$X_k = \|x_k - x^*\|^2, Y_k^{(T-1)} = \mathbb{E}[\|y_k^{(T-1)} - f^{(T)}(x_k)\|^2|\mathbb{F}_k],$$
$$Y_k^{(T-2)} = \mathbb{E}[\|y_k^{(T-2)} - f^{(T-1)}(y_k^{(T-1)})\|^2|\mathbb{F}_k], \cdots, Y_k^{(1)} = \mathbb{E}[\|y_k^{(1)} - f^{(2)}(y_k^{(2)})\|^2|\mathbb{F}_k],$$
$$\eta_k = [\frac{\alpha_k^2}{\beta_{T-1,k}} + \cdots + \frac{\alpha_k^2}{\beta_{1,k}}]C_0, u_k^{(T)} = 2\alpha_k(F(x_k) - F^*),$$
$$u_k^{(1)} = u_k^{(2)} = \cdots = u_k^{(T-1)} = 0, c_1 = 1, \cdots, c_{T-1} = T-1,$$
$$\mu_k^{(1)} = 2\beta_{1,k}^2 V_1 + \mathcal{O}(\frac{\mathbb{E}[\|y_{k+1}^{(2)} - y_k^{(2)}\|^2|\mathbb{F}_k]}{\beta_{1,k}}), \cdots,$$
$$\mu_k^{(T-2)} = 2\beta_{T-2,k}^2 V_{T-1} + \mathcal{O}(\frac{\mathbb{E}[\|y_{k+1}^{(T-1)} - y_k^{(T-1)}\|^2|\mathbb{F}_k]}{\beta_{T-2,k}}),$$
$$\mu_k^{(T-1)} = 2\beta_{T-1,k}^2 V_T + \mathcal{O}(\frac{\mathbb{E}[\|x_{k+1} - x_k\|^2|\mathbb{F}_k]}{\beta_{T-1,k}}),$$
$$\mu_k^{(T)} = \alpha_k^2 C_1 C_2 \cdots C_T, \theta_j^{(1)} = \beta_{1,k}, \cdots, \theta_j^{(T-1)} = \beta_{T-1,k},$$

we obtain the sequence $\{\|x_k - x^*\|\}$ converges almost surely to a random variable, and

$$\sum_{k=0}^{\infty} \alpha_k (F(x_k) - F^*) < \infty,$$

which further implies

$$\liminf_{k\to\infty} F(x_k) = F^*, \qquad w.p.1.$$

Using Lemma 2.5, we conclude that the sequence $\{x_k\}$ converges almost surely to a random point in the set of optimal solutions to problem (1.1).

b) Since problem (1.1) has at least one optimal solution, the function $F$ is bounded from below, and denote by $F^*$ the minimal value of $F(x)$ over $\mathcal{X}$. As a result, we can treat $\{F(x_k) - F^*\}$ as a nonnegative random variable. By Lemma 2.2, we have

$$\begin{aligned} &\mathbb{E}[F(x_{k+1}) - F^*|\mathbb{F}_k] \\ \le &F(x_k) - F^* - \frac{\alpha_k}{2}\|\nabla F(x_k)\|^2 + \frac{1}{2}\alpha_k^2 L_F C_1 C_2 \cdots C_T + (T-1)\beta_{T-1,k}\mathbb{E}[\|y_{k+1}^{(T-1)} - f^{(T)}(x_k)\|^2|\mathbb{F}_k] \\ &+ \cdots + \beta_{1,k}\mathbb{E}[\|y_{k+1}^{(1)} - f^{(2)}(y_{k+1}^{(2)})\|^2|\mathbb{F}_k], \end{aligned} \tag{A.21}$$

for sufficiently large $k$. We apply the $T$-element super-martingale convergent lemma to Eqs.(A.20)-

(A.21). By letting

$$
\begin{aligned}
&X_k = F(x_{k+1}) - F^*, Y_k^{(T-1)} = \mathbb{E}[\|y_k^{(T-1)} - f^{(T)}(x_k)\||\mathbb{F}_k],\\
&Y_k^{(T-2)} = \mathbb{E}[\|y_k^{(T-2)} - f^{(T-1)}(y_k^{(T-1)})\|^2|\mathbb{F}_k], \cdots\\
&Y_k^{(1)} = \mathbb{E}[\|y_k^{(1)} - f^{(2)}(y_k^{(2)})\|^2|\mathbb{F}_k],\\
&\eta_k = 0, u_k^{(T)} = \frac{1}{2}\alpha_k\|\nabla F(x_k)\|^2,\\
&u_k^{(1)} = u_k^{(2)} = \cdots = u_k^{(T-1)} = 0, c_1 = 1, \cdots, c_{T-1} = T-1,\\
&\mu_k^{(1)} = 2\beta_{1,k}^2 V_2 + \mathcal{O}(\frac{\mathbb{E}[\|y_{k+1}^{(2)} - y_k^{(2)}\|^2|\mathbb{F}_k]}{\beta_{1,k}}), \cdots,\\
&\mu_k^{(T-2)} = 2\beta_{T-2,k}^2 V_{T-1} + \mathcal{O}(\frac{\mathbb{E}[\|y_{k+1}^{(T-1)} - y_k^{(T-1)}\|^2|\mathbb{F}_k]}{\beta_{T-2,k}}),\\
&\mu_k^{(T-1)} = 2\beta_{T-1,k}^2 V_T + \mathcal{O}(\frac{\mathbb{E}[\|x_{k+1} - x_k\|^2|\mathbb{F}_k]}{\beta_{T-1,k}}),\\
&\mu_k^{(T)} = \frac{1}{2}\alpha_k^2 L_F C_1 C_2 \cdots C_T, \theta_j^{(1)} = \beta_{1,k}, \cdots, \theta_j^{(T-1)} = \beta_{T-1,k},
\end{aligned}
$$

we obtain that $\{F(x_k) - F^*\}$ converges almost surely to a random variable, and

$$\sum_{k=0}^{\infty} \alpha_k \|\nabla F(x_k)\|^2 < \infty, \quad w.p.1.$$

Using Lemma 2.6, we conclude that any limiting point of the sequence $\{x_k\}$ is a stationary point with probability 1, which completes the proof. □

# B Proof of Theorem 2.2

Note that we let the step-sizes be $\alpha_k = k^{-a}, \beta_{T-1,k} = k^{-b_{T-1}}, \cdots, \beta_{1,k} = k^{-b_1}$. We slightly modify Lemma 5 in [3] to help us derive the convergence rates.

**Lemma B.1.** Given a sequence of positive real numbers $\{w_k\}_{k=1}^{\infty}$ satisfying

$$w_{k+1} \le (1 - C_1\beta_k)w_k + C_2 k^{-a},$$

where $C_1 \ge 0$, $C_2 \ge 0$, and $a \ge 0$. Choosing $c = a - b$ and $\beta_k$ to be $\beta_k = C_3 k^{-b}$ where $b \in (0, 1]$ and $C_3 > c/C_1$, the sequence can be bounded by $w_k \le Ck^{-c}$, where $C$ is defined as

$$C := w_0 + \frac{C_2}{C_1 C_3 - c}.$$

In other words, we have

$$w_k \le \mathcal{O}(k^{-a+b}).$$

PROOF: We prove it by induction. Clearly, the claim holds for $k = 0$. Next, suppose the claim holds for "$k$", we prove it also true for "$k+1$". That is, given $w_k \le Ck^{-c}$, we need to prove $w_{k+1} \le C(k+1)^{-c}$.

$$\begin{aligned} w_{k+1} \le& (1 - C_1\beta_k)w_k + Ck^{-a} \\ \le& (1 - C_1C_3k^{-b})Ck^{-c} + C_2k^{-a} \\ =& Ck^{-c} - CC_1C_3k^{-b-c} + C_2k^{-a}. \end{aligned} \tag{B.1}$$

To prove B.1 is bounded by $C(k+1)^{-c}$, it suffices to show that

$$\Delta := (k+1)^{-c} - k^{-c} + C_1C_3k^{-b-c} > 0 \text{ and } C \ge \frac{C_2k^{-a}}{\Delta}.$$

From the convexity of function $h(t) = t^{-c}$, we have the inequality $(k+1)^{-c} - k^{-c} \ge -ck^{-c-1}$. Therefore we obtain

$$\Delta \ge -ck^{-c-1} + C_3k^{-b-c} \ge (C_1C_3 - c)k^{-b-c} \ge 0.$$

To verify the second one, we have

$$\frac{C_2k^{-a}}{\Delta} \le \frac{C_2}{C_1C_3 - 2}k^{-a+b+c} = \frac{C_2}{C_1C_3 - c} \le C,$$

where the equality comes from the definition that $c = a - b$ and the last inequality holds by the definition of $C$. It completes the proof. □

## B.1 Proof of Lemma 2.7

PROOF: By Lemma 2.3, for any $j = T-1, \cdots, 1$, $\mathbb{E}[\|y_{k+1}^{(j+1)} - y_k^{(j+1)}\|^2] \le \mathcal{O}(\beta_{j+1,k}^2)$, and

$$\begin{aligned} &\mathbb{E}[\|y_{k+1}^{(j)} - f^{(j+1)}(y_{k+1}^{(j+1)})\|^2|\mathbb{F}_k] \\ \le& (1 - \beta_{j,k})\mathbb{E}[\|y_k^{(j)} - f^{(j+1)}(y_k^{(j+1)})\|^2|\mathbb{F}_k] + \beta_{j,k}^{-1}C_{j+1}\mathbb{E}[\|y_{k+1}^{(j+1)} - y_k^{(j)}\|^2|\mathbb{F}_k] + 2V_{j+1}\beta_{j,k}^2, \end{aligned}$$

with probability 1. Thus,

$$\begin{aligned} &\mathbb{E}[\|y_{k+1}^{(j)} - f^{(j+1)}(y_{k+1}^{(j+1)})\|^2] \\ \le& (1 - \beta_{j,k})\mathbb{E}[\|y_k^{(j)} - f^{(j+1)}(y_k^{(j+1)})\|^2] + \beta_{j,k}^{-1}C_{j+1}\mathbb{E}[\|y_{k+1}^{(j+1)} - y_k^{(j+1)}\|^2] + 2V_{j+1}\beta_{j,k}^2 \\ \le& (1 - \beta_{j,k})\mathbb{E}[\|y_k^{(j)} - f^{(j+1)}(y_k^{(j+1)})\|^2] + \beta_{j,k}^{-1}\mathcal{O}(\beta_{j+1,k}^2) + 2V_{j+1}\beta_{j,k}^2. \end{aligned} \tag{B.2}$$

Substitute $\alpha_k = k^{-a}$, $\beta_{j,k} = 2k^{-b_j}$ into Eq.(B.2), we get

$$\begin{aligned} &\mathbb{E}[\|y_k^{(j)} - f^{(j+1)}(y_k^{(j+1)})\|^2] \\ \le& (1 - k^{-b_j})\mathbb{E}[\|y_k^{(j)} - f^{(j+1)}(y_k^{(j+1)})\|^2] + \mathcal{O}(k^{-2b_{j+1}+b_j}) + 2V_{j+1}k^{-2b_j}. \end{aligned}$$

By Lemma B.1, we obtain

$$\mathbb{E}[\|y_k^{(j)} - f^{(j+1)}(y_k^{(j+1)})\|^2] \le \mathcal{O}(k^{-2b_{j+1}+2b_j}) + \mathcal{O}(k^{-b_j}).$$

□

### B.2 Proof of Theorem 2.2

PROOF: Define the random variable

$$J_k = \|x_k - x^*\|^2 + (T-1)\mathbb{E}[\|y_k^{(T-1)} - f^{(T)}(x_k)\|^2|\mathbb{F}_k] + \cdots + \mathbb{E}[\|y_k^{(1)} - f^{(2)}(y_k^{(2)})\|^2|\mathbb{F}_k],$$

so we have $\mathbb{E}[J_k] \leq D_x + (T-1)D_{T-1} + (T-2)D_{T-2} + \cdots D_1 \equiv D_J$. In this basic algorithm, all steps are updated by the basic update rule. We multiply Eq.(A.20) by $j \times (1+\beta_{j,k})$ for every $j = T-1, ..., 1$ and take their sums with Eq.(A.19). We obtain

$$\begin{aligned}
&\mathbb{E}[J_{k+1}|\mathbb{F}_k] \leq \Big(1 + (\frac{\alpha_k^2}{\beta_{T-1,k}} + \frac{\alpha_k^2}{\beta_{T-2,k}} + \cdots + \frac{\alpha_k^2}{\beta_{1,k}})\Big) J_k - 2\alpha_k(F(x_k - F^*)) \\
&+ C_1C_2\cdots C_T\alpha_k^2 + 2(T-1)V_T\beta_{T-1,k}^2(1+\beta_{T-1,k}) + \frac{(T-1)C_T(1+\beta_{T-1,k})}{\beta_{T-1,k}}\mathbb{E}[\|x_{k+1} - x_k\|^2|\mathbb{F}_k] \\
&+ \sum_{j=1}^{T-2} j \times \Big\{2V_{j+1}\beta_{j,k}^2(1+\beta_{j,k}) + \frac{C_{j+1}(1+\beta_{j,k})}{\beta_{j,k}}\mathbb{E}[\|y_{k+1}^{(j+1)} - y_k^{(j+1)}\|^2|\mathbb{F}_k]\Big\}.
\end{aligned}$$

Taking expectation on both sides using the fact $1+\beta_{j,k} \leq 2$ for all $j = T-1, \cdots, 1$ and Eq.(A.12), we obtain

$$\begin{aligned}
\mathbb{E}[J_{k+1}] \leq &\Big(1 + (\frac{\alpha_k^2}{\beta_{T-1,k}} + \frac{\alpha_k^2}{\beta_{T-2,k}} + \cdots + \frac{\alpha_k^2}{\beta_{1,k}})\Big)\mathbb{E}[J_k] - 2\alpha_k(F(x_k - F^*)) \\
&+ C_1C_2\cdots C_T\alpha_k^2 + 4(T-1)V_T\beta_{T-1,k}^2 + \frac{2(T-1)C_1C_2\cdots C_{T-1}C_T^2\alpha_k^2}{\beta_{T-1,k}} \\
&+ \sum_{j=1}^{T-2} j\Big\{4V_{j+1}\beta_{j,k}^2 + \frac{2C_{j+1}}{\beta_{j,k}}\mathbb{E}[\|y_{k+1}^{(j+1)} - y_k^{(j+1)}\|^2]\Big\}.
\end{aligned}$$

By Lemma 2.4, we have $\mathbb{E}[\|y_{k+1}^{(j)} - y_k^{(j)}\|^2] \leq \mathcal{O}(\beta_{j,k}^2)$ for $j = T-1, \cdots, 1$. Substitute this into the previous inequality, we obtain

$$\begin{aligned}
\mathbb{E}[J_{k+1}] \leq &\Big(1 + (\frac{\alpha_k^2}{\beta_{T-1,k}} + \frac{\alpha_k^2}{\beta_{T-2,k}} + \cdots + \frac{\alpha_k^2}{\beta_{1,k}})C_0\Big)\mathbb{E}[J_k] - 2\alpha_k(F(x_k - F^*)) \\
&+ C_1C_2\cdots C_T\alpha_k^2 + 4(T-1)V_T\beta_{T-1,k}^2 + \frac{2(T-1)C_1C_2\cdots C_{T-1}C_T^2\alpha_k^2}{\beta_{T-1,k}} \\
&+ \sum_{j=1}^{T-2} j \times \Big\{4V_{j+1}\beta_{j,k}^2 + \frac{2C_{j+1}}{\beta_{j,k}}\mathcal{O}(\beta_{j+1,k}^2)\Big\}.
\end{aligned}$$

Let $N > 0$, by reordering the terms in the preceding relation and taking its sum over $k-N, \cdots, k$, we have

$$\begin{aligned}
& 2\sum_{t=k-N}^{k} \mathbb{E}[F(x_t) - F^*] \\
\leq & \sum_{t=k-N}^{k} \frac{1}{\alpha_t}\Big\{\Big(1 + (\frac{\alpha_k}{\beta_{T-1,k}} + \cdots + \frac{\alpha_k}{\beta_{1,k}})C_0\Big)\mathbb{E}[J_t] - \mathbb{E}[J_{t+1}]\Big\} \\
& + \sum_{t=k-N}^{k} \Big(\mathcal{O}(\alpha_k) + \sum_{j=1}^{T-1}\mathcal{O}(\frac{\beta_{j,k}^2}{\alpha_k}) + \mathcal{O}(\frac{\alpha_k}{\beta_{T-1,k}}) + \sum_{j=1}^{T-2}\mathcal{O}(\frac{\beta_{j+1,k}^2}{\alpha_k\beta_{j,k}})\Big) \\
= & \sum_{t=k-N}^{k} (\frac{1}{\alpha_t} - \frac{1}{\alpha_{t-1}})\mathbb{E}[J_t] - \frac{1}{\alpha_k}\mathbb{E}[J_{k+1}] + \frac{1}{\alpha_{k-N-1}}\mathbb{E}[J_{k-N}] + \sum_{t=k-N}^{k}(\frac{\alpha_k}{\beta_{T-1,k}} + \cdots + \frac{\alpha_k}{\beta_{1,k}})C_0\mathbb{E}[J_k] \\
& + \sum_{t=k-N}^{k} \Big(\mathcal{O}(\alpha_k) + \sum_{j=1}^{T-1}\mathcal{O}(\frac{\beta_{j,k}^2}{\alpha_k}) + \mathcal{O}(\frac{\alpha_k}{\beta_{T-1,k}}) + \sum_{j=1}^{T-2}\mathcal{O}(\frac{\beta_{j+1,k}^2}{\alpha_k\beta_{j,k}})\Big) \\
\leq & \sum_{t=k-N}^{k} (\frac{1}{\alpha_t} - \frac{1}{\alpha_{t-1}})D_J + \frac{1}{\alpha_{k-N-1}}D_J + \sum_{t=k-N}^{k}(\frac{\alpha_k}{\beta_{T-1,k}} + \cdots + \frac{\alpha_k}{\beta_{1,k}})C_0 D_J \\
& + \sum_{t=k-N}^{k}\mathcal{O}(\alpha_k) + \sum_{j=1}^{T-1}\sum_{t=k-N}^{k}\mathcal{O}(\frac{\beta_{j,k}^2}{\alpha_k}) + \sum_{t=k-N}^{k}\mathcal{O}(\frac{\alpha_k}{\beta_{T-1,k}}) + \sum_{j=1}^{T-2}\sum_{t=k-N}^{k}\mathcal{O}(\frac{\beta_{j+1,k}^2}{\alpha_k\beta_{j,k}}).
\end{aligned}$$

Since $\sum_{t=k-N}^{k}(\frac{1}{\alpha_t} - \frac{1}{\alpha_{t-1}})D_J + \frac{1}{\alpha_{k-N-1}}D_J = \frac{1}{\alpha_k}D_J$, let $\alpha_k = k^{-a}$, $\beta_{T-1,k} = k^{-b_{T-1}}$, $\cdots$, and $\beta_{1,k} = k^{-b_1}$, we obtain

$$\begin{aligned}
& 2\sum_{t=k-N}^{k} \mathbb{E}[F(x_t) - F^*] \\
\leq & k^a D + \sum_{t=k-N}^{k}(k^{-a+b_{T-1}} + \cdots + k^{-a+b_1})C_0 D_J \\
& + \sum_{t=k-N}^{k}\mathcal{O}(k^{-a}) + \sum_{j=1}^{T-1}\sum_{t=k-N}^{k}\mathcal{O}(k^{-2b_j+a}) + \sum_{t=k-N}^{k}\mathcal{O}(k^{-a+b_{T-1}}) + \sum_{j=1}^{T-2}\sum_{t=k-N}^{k}\mathcal{O}(k^{-2b_{j+1}+a+b_j}).
\end{aligned}$$

Choosing $a = 1 - 1/2^T$, $b_{T-1} = 1 - 1/2^{T-1}$, $\cdots$, $b_1 = 1 - 1/2 = 1/2$, it is derived that

$$\begin{aligned}
& 2\sum_{t=k-N}^{k} \mathbb{E}[F(x_t) - F^*] \\
\leq & k^{1-1/2^T} D_J + \sum_{t=k-N}^{k}(k^{-1/2^T} + \cdots + k^{1/2^T - 1/2})C_0 D_J \\
& + \sum_{t=k-N}^{k}\mathcal{O}(k^{-1+1/2^T}) + \sum_{j=1}^{T-1}\sum_{t=k-N}^{k}\mathcal{O}(k^{-1+1/2^{j-1}-1/2^T}) + \sum_{t=k-N}^{k}\mathcal{O}(k^{-1/2^T}) + \sum_{j=1}^{T-2}\sum_{t=k-N}^{k}\mathcal{O}(k^{-1/2^T}).
\end{aligned}$$

Since for any $p > 0$, $\sum_{t=k-N}^{k} \mathcal{O}(k^p) \le \mathcal{O}(k^{1+p} - (k-N)^{1+p})$, using the convexity of $F$ and taking $N = N_k = \frac{k}{2}$, we obtain

$$\mathbb{E}[F(\widehat{x}_k) - F^*] \le \frac{1}{N_k} \sum_{t=k-N}^{k} \mathbb{E}[F(x_t) - F^*] \le \mathcal{O}(k^{-1/2^T}),$$

which completes the proof. □

# C Proof of Theorem 3.1

## C.1 Proof of Lemma 3.1

Note that here we use the same notation as in the proof of Lemma 2.3 . Consider a basic update step

$$y_{k+1}^{(j)} = (1-\beta_{j,k}) y_k^{(j)} + \beta_{j,k} f_{w_{j+1,k+1}}^{(j+1)}(y_{k+1}^{(j+1)}),$$

and denote it by

$$z_{k+1} = (1-\gamma_k) z_k + \gamma_k h_{u_{k+1}}(x_{k+1}).$$

Now we show the detailed proof of Lemma 3.1.

Proof: (a) Under the assumption $\mathbb{E}[\|x_{k+1} - x_k\|^4] \le \mathcal{O}(\alpha_k^4)$, there exists a constant $C_0 > 0$ such that $\mathbb{E}[\|x_{k+1} - x_k\|^4] \le C_0 \alpha_k^4$. By Lemma 2.3, there exists a constant $D_z > 0$ such that $\mathbb{E}[\|z_{k+1} - h(x_{k+1})\|^2] \le D_z$, and we also have $\mathbb{E}[\|z_{k+1} - h(x_{k+1})\|] \le \sqrt{\mathbb{E}[\|z_{k+1} - h(x_{k+1})\|^2]} \le \sqrt{D_z}$.

Let $e_{k+1} = (1-\gamma_k)\big(h(x_{k+1}) - h(x_k)\big)$. Together with the definition of $z_{k+1}$, we get

$$z_{k+1} - h(x_{k+1}) + e_{k+1} - e_{k+1} = (1-\gamma_k)\big(z_k - h(x_k)\big) + \gamma_k\big(h_{u_{k+1}}(x_{k+1}) - h(x_{k+1})\big) - e_{k+1}. \quad \text{(C.1)}$$

By the Lipschitz continuity of $h$, we obtain

$$\|e_{k+1}\| \le (1-\gamma_k)\sqrt{C_h}\|x_{k+1} - x_k\|.$$

Meanwhile, we have

$$\|z_{k+1} - h(x_{k+1})\| \le \|(1-\gamma_k)(z_k - h(x_k)) + \gamma_k(h_{u_{k+1}}(x_{k+1}) - h(x_{k+1}))\| + \|e_{k+1}\|.$$

Let $P_k = \|(1-\gamma_k)(z_k - h(x_k)) + \gamma_k(h_{u_{k+1}}(x_{k+1}) - h(x_{k+1}))\|$, considering the fourth moment, we get

$$\begin{aligned}
P_k^4 =& \|(1-\gamma_k)(z_k - h(x_k)) + \gamma_k(h_{u_{k+1}}(x_{k+1}) - h(x_{k+1}))\|^4 \\
=& (1-\gamma_k)^4\|z_k - h(x_k)\|^4 + 4(1-\gamma_k)^3\gamma_k(z_k - h(x_k))^3(h_{u_{k+1}}(x_{k+1}) - h(x_{k+1})) \\
&+ 6(1-\gamma_k)^2\gamma_k^2\|z_k - h(x_k)\|^2\|h_{u_{k+1}}(x_{k+1}) - h(x_{k+1})\|^2 \\
&+ 4(1-\gamma_k)\gamma_k^3(z_k - h(x_k))(h_{u_{k+1}}(x_{k+1}) - h(x_{k+1}))^3 \\
&+ \gamma_k^4\|h_{u_{k+1}}(x_{k+1}) - h(x_{k+1})\|^4 \\
\le& (1-\gamma_k)^4\|z_k - h(x_k)\|^4 + 4(1-\gamma_k)^3\gamma_k(z_k - h(x_k))^3(h_{u_{k+1}}(x_{k+1}) - h(x_{k+1})) \\
&+ 6(1-\gamma_k)^2\gamma_k^2\|z_k - h(x_k)\|^2\|h_{u_{k+1}}(x_{k+1}) - h(x_{k+1})\|^2 \\
&+ 4(1-\gamma_k)\gamma_k^3\|z_k - h(x_k)\|\|h_{u_{k+1}}(x_{k+1}) - h(x_{k+1})\|^3 \\
&+ \gamma_k^4\|h_{u_{k+1}}(x_{k+1}) - h(x_{k+1})\|^4.
\end{aligned}$$

So we obtain

$$
\begin{aligned}
\mathbb{E}[P_k^4] \leq & (1-\gamma_k)^4 \mathbb{E}[\|z_k - h(x_k)\|^4] + 0 \\
& + 6(1-\gamma_k)^2\gamma_k^2 \mathbb{E}\Big[\|z_k - h(x_k)\|^2 \mathbb{E}[\|h_{u_{k+1}}(x_{k+1}) - h(x_{k+1})\|^2 | \mathbb{F}_{k+1}]\Big] \\
& + 4(1-\gamma_k)\gamma_k^3 \mathbb{E}\Big[\|z_k - h(x_k)\| \mathbb{E}[\|h_{u_{k+1}}(x_{k+1}) - h(x_{k+1})\|^3 | \mathbb{F}_{k+1}]\Big] \\
& + \gamma_k^4 \mathbb{E}[\|h_{u_{k+1}}(x_{k+1}) - h(x_{k+1})\|^4] \\
\leq & (1-\gamma_k)^4 \mathbb{E}[\|z_k - h(x_k)\|^4] + 6(1-\gamma_k)^2\gamma_k^2 D V_h + 4(1-\gamma_k)\gamma_k^3 \sqrt{D} V_h^{3/2} + \gamma_k^4 V_h^2.
\end{aligned}
$$

Using the fact that $\|a+b\|^2 \leq (1+\epsilon)\|a\|^2 + (1+1/\epsilon)\|b\|^2$ and $\|a+b\|^4 \leq (1+\epsilon)^3\|a\|^4 + (1+1/\epsilon)^3\|b\|^4$ for any $\epsilon > 0$, we have

$$
\|z_{k+1} - h(x_{k+1})\|^4 \leq (1+\gamma_k)^3 P_k^4 + (1+1/\gamma_k)^3 \|e_{k+1}\|^4, \tag{C.2}
$$

and

$$
\begin{aligned}
& \mathbb{E}[\|z_{k+1} - h(x_{k+1})\|^4] \\
\leq & (1+\gamma_k)^3 \Big\{ (1-\gamma_k)^4 \mathbb{E}[\|z_k - h(x_k)\|^4] + 6(1-\gamma_k)^2\gamma_k^2 D_z V_h + 4(1-\gamma_k)\gamma_k^3 \sqrt{D_z} V_h^{3/2} + \gamma_k^4 V_h^2 \Big\} \\
& + (1+1/\gamma_k)^3 (1-\gamma_k)^4 C_h^2 \mathbb{E}[\|x_{k+1} - x_k\|^4] \\
\leq & (1-\gamma_k)\mathbb{E}[\|z_k - h(x_k)\|^4] + 12\gamma_k^2 D_z V_h + 16\gamma_k^3 \sqrt{D_z} V_h^{3/2} + 8\gamma_k^4 V_h^2 + \frac{C_h^2}{\gamma_k^3} \mathbb{E}[\|x_{k+1} - x_k\|^4] \\
\leq & (1-\gamma_k)\mathbb{E}[\|z_k - h(x_k)\|^4] + 12\gamma_k^2 D_z V_h + 16\gamma_k^3 \sqrt{D_z} V_h^{3/2} + 8\gamma_k^4 V_h^2 + \frac{\alpha_k^4}{\gamma_k^3} C_h^2 C_0.
\end{aligned} \tag{C.3}
$$

Finally, we complete the proof by induction. Since $\alpha_k/\gamma_k \to 0$, there exists a constant $M > 0$ such that $\alpha_k/\gamma_k \leq M$ for all $k$. Let $S_z = \|z_0 - h(x_0)\|^4 + 12 D_z V_h + 16\sqrt{D_z} V_h^{3/2} + 8M^4 V_h^2 + C_h^2 C_0$, then $\|z_0 - h(x_0)\|^4 \leq S_z$, and $S_z - 12 D_z \gamma_k V_h - 16\gamma_k^2 \sqrt{D_z} V_h^{3/2} - 8\gamma_k^4 V_h^2 - \frac{\alpha_k^4}{\gamma_k^4} C_h^2 C_0 \geq 0$ for all $k$. Suppose the claim is true for $0, 1, \cdots, k$, then

$$
\begin{aligned}
& \mathbb{E}[\|z_{k+1} - h(x_{k+1})\|^4] \\
\leq & (1-\gamma_k) S_z + 12\gamma_k^2 D_z V_h + 16\gamma_k^3 \sqrt{D_z} V_h^{3/2} + 8\gamma_k^4 V_h^2 + \frac{\alpha_k^4}{\gamma_k^3} C_h^2 C_0 \\
= & S_z - \gamma_k \Big( S_z - 12\gamma_k D_z V_h - 16\gamma_k^2 \sqrt{D_z} V_h^{3/2} - 8\gamma_k^4 V_h^2 - \frac{\alpha_k^4}{\gamma_k^4} C_h^2 C_0 \Big) \\
\leq & S_z,
\end{aligned}
$$

which completes the proof.

b) By the definition of $z_{k+1}$, we have $z_{k+1} = (1-\gamma_k) z_k + \gamma_k h_{u_{k+1}}(x_{k+1})$, and

$$
\begin{aligned}
(1-\gamma_k)^4 \|z_{k+1} - z_k\|^4 = & \gamma_k^4 \|h_{u_{k+1}}(x_{k+1}) - z_{k+1}\|^4 \\
\leq & 8\gamma_k^4 \|h_{u_{k+1}}(x_{k+1}) - h(x_{k+1})\|^4 + 8\gamma_k^4 \|z_{k+1} - h(x_{k+1})\|^4,
\end{aligned}
$$

where the inequality comes from $(a+b)^4 \le (2a^2+2b^2)^2 \le 8a^4+8b^4$. Thus

$$\mathbb{E}[\|z_{k+1}-z_k\|^4] \le \frac{8\gamma_k^4}{(1-\gamma_k)^4}V_h^2 + \frac{8\gamma_k^4}{(1-\gamma_k)^4}\mathbb{E}[\|z_{k+1}-h(x_{k+1})\|^4].$$

By part (a) that there exists a constant $S_z \ge 0$ such that $\mathbb{E}[\|z_{k+1}-h(x_{k+1})\|^4] \le S_z$, we obtain

$$\mathbb{E}[\|z_{k+1}-z_k\|^4] \le \mathcal{O}(\gamma_k^4),$$

which completes the proof. □

## C.2 Proof of Lemma 3.2

Let Assumption 2.1 and 3.1 hold, we apply the basic update rule to the first inner level and accelerated update rule to the remaining levels. We consider the analysis for the second inner level updated as:

$$\begin{aligned}
\widehat{y}_{k+1}^{(T-2)} &= (1-1/\beta_{T-2,k})y_k^{(T-1)} + y_{k+1}^{(T-1)}/\beta_{T-2,k},\\
y_{k+1}^{(T-2)} &= (1-\beta_{T-2,k})y_k^{(T-2)} + \beta_{T-2,k}\cdot f_{\omega_{T-1,k+1}}^{(T-1)}(\widehat{y}_{k+1}^{(T-2)}).
\end{aligned}$$

For ease of notation, we denote $y_k^{(T-1)}$ by $z_k$, $\beta_{T-2,k}$ by $\beta_k$, $\widehat{y}_k^{(T-2)}$ by $\widehat{y}_{k+1}$, $y_k^{(T-2)}$ by $y_k$ and $f_{\omega_{T-1,k+1}}^{(T-1)}(\cdot)$ by $g_{w_{k+1}}(\cdot)$. The corresponding update step can be written as

$$\begin{aligned}
\widehat{y}_{k+1} &= (1-\frac{1}{\beta_k})z_k + \frac{1}{\beta_k}z_{k+1},\\
y_{k+1} &= (1-\beta_k)y_k + \beta_k g_{w_{k+1}}(\widehat{y}_{k+1}).
\end{aligned}$$

To construct a super-martingale involving $\{y_{k+1}-g(z_{k+1})\}$ while utilizing the smoothness of $g(\cdot)$, we write $y_k$ as the weighted average of samples at the extrapolated points $\{g_{w_t}(\widehat{y}_t)\}_{t=0}^k$ with weights $\zeta_t^{(k)}$ defined as

$$\zeta_t^{(k)} = \begin{cases} \beta_t \Pi_{i=t+1}^k (1-\beta_i) & \text{if } \ k > t \ge 0,\\ \beta_k & \text{if } \ k = t \ge 0. \end{cases}$$

By the definitions of $\zeta_t^{(k)}$, $z_k$, and $y_k$, we have the following identities

$$\zeta_t^{(k+1)} = (1-\beta_{k+1})\zeta_t^{(k)}, \sum_{t=0}^k \zeta_t^{(k)} = 1,$$

and

$$z_{k+1} = \sum_{t=0}^k \zeta_t^{(k)}\widehat{y}_{t+1}, y_{k+1} = \sum_{t=0}^k \zeta_t^{(k)} g_{w_{t+1}}(\widehat{y}_{t+1}).$$

In what follows, we prove Lemma 3.2 through a series of lemmas. These lemmas are also the building blocks for the subsequent rate of convergence analysis. It can be seen that $\widehat{y}_k$ plays a role as a link between $y_k$ and $z_k$, and it is natural to consider decomposing $\{y_k - g(z_k)\}$ into some terms containing $\widehat{y}_k$ to start our analysis, which is presented in the following lemma.

**Lemma C.1** (Decomposition of $y_k - g(z_k)$)**.** Suppose that Assumption 2.1 holds, then we have

$$\|y_{k+1} - g(z_{k+1})\| \leq L_g \sum_{t=0}^{k} \zeta_t^{(k)} \|z_{k+1} - \widehat{y}_{t+1}\|^2 + \|\sum_{t=0}^{k} \zeta_t^{(k)} (g_{w_{t+1}}(\widehat{y}_{t+1}) - g(\widehat{y}_{t+1}))\|.$$

Proof: By definition, we have

$$\begin{aligned}
y_{k+1} &= \sum_{t=0}^{k} \zeta_t^{(k)} g_{w_{t+1}}(\widehat{y}_{t+1}) \\
&= \sum_{t=0}^{k} \zeta_t^{(k)} g(\widehat{y}_{t+1}) + \sum_{t=0}^{k} \zeta_t^{(k)} [g_{w_{t+1}}(\widehat{y}_{t+1}) - g(\widehat{y}_{t+1})] \\
&= \sum_{t=0}^{k} \zeta_t^{(k)} [g(z_{k+1}) + \nabla g(z_{k+1})^T (\widehat{y}_{t+1} - z_{k+1}) + \kappa(z_{k+1}, \widehat{y}_{t+1}))] \\
&\quad + \sum_{t=0}^{k} \zeta_t^{(k)} (g_{w_{t+1}}(\widehat{y}_{t+1}) - g(\widehat{y}_{t+1})) \\
&= (\sum_{t=0}^{k} \zeta_t^{(k)}) g(z_{k+1}) + \nabla g(z_{k+1})^T \sum_{t=0}^{k} \zeta_t^{(k)} (\widehat{y}_{t+1} - z_{k+1}) \\
&\quad + \sum_{t=0}^{k} \zeta_t^{(k)} \kappa(z_{k+1}, \widehat{y}_{t+1}) + \sum_{t=0}^{k} \zeta_t^{(k)} (g_{w_{t+1}}(\widehat{y}_{t+1}) - g(\widehat{y}_{t+1})),
\end{aligned} \tag{C.4}$$

where the second equality holds by the Taylor expansion of $g(\cdot)$ at $z_{k+1}$ that

$$g(\widehat{y}_{t+1}) = g(z_{k+1}) + \nabla g(z_{k+1})'(\widehat{y}_{t+1} - z_{k+1}) + \kappa(z_{k+1}, \widehat{y}_{t+1}),$$

with $\kappa(z_{k+1}, \widehat{y}_{t+1})$ summarizing the second and higher order terms. An important observation here is that the first order term in Eq.(C.4) cancels out that

$$\nabla g(z_{k+1})' \sum_{t=0}^{k} \zeta_t^{(k)} (\widehat{y}_{t+1} - z_{k+1}) = \nabla g(z_{k+1})' \Big( \sum_{t=0}^{k} \zeta_t^{(k)} \widehat{y}_{t+1} - (\sum_{t=0}^{k} \zeta_t^{(k)}) z_{k+1} \Big) = 0,$$

where we use the identities $z_{k+1} = \sum_{t=0}^{k} \zeta_t^{(k)} \widehat{y}_{t+1}$ and $\sum_{t=0}^{k} \zeta_t^{(k)} = 1$. Next, it follows from (C.4) that

$$y_{k+1} = g_{z_{k+1}} + 0 + \sum_{t=0}^{k} \zeta_t^{(k)} \kappa(z_{k+1}, \widehat{y}_{t+1}) + \sum_{t=0}^{k} \zeta_t^{(k)} [g_{w_{t+1}}(\widehat{y}_{t+1}) - g(\widehat{y}_{t+1})].$$

Using triangle inequality, we have

$$\begin{aligned}
\|y_{k+1} - g(z_{k+1})\| &\leq \sum_{t=0}^{k} \zeta_t^{(k)} \|\kappa(z_{k+1}, \widehat{y}_{t+1})\| + \|\sum_{t=0}^{k} \zeta_t^{(k)} (g_{w_{t+1}}(\widehat{y}_{t+1}) - g(\widehat{y}_{t+1}))\| \\
&\leq \frac{L_g}{2} \sum_{t=0}^{k} \zeta_t^{(k)} \|z_{k+1} - \widehat{y}_{t+1}\|^2 + \|\sum_{t=0}^{k} \zeta_t^{(k)} (g_{w_{t+1}}(\widehat{y}_{t+1}) - g(\widehat{y}_{t+1}))\|,
\end{aligned}$$

where the second inequality holds by Assumption 2.1 that $\|\kappa(z_{k+1}, \widehat{y}_{t+1})\| \leq \frac{L_g}{2}\|z_{k+1} - \widehat{y}_{t+1}\|^2$. □

In Lemma C.1, we bound $\{y_k - g(z_k)\}$ by the sum of two iterative sums $\sum_{t=0}^{k} \zeta_t^{(k)} \|z_{k+1} - \widehat{y}_{t+1}\|^2$ and $\sum_{t=0}^{k} \zeta_t^{(k)} (g_{w_{t+1}}(\widehat{y}_{t+1}) - g(\widehat{y}_{t+1}))$. In the next lemma, we consider convergence properties of the two iterative sums.

**Lemma C.2.** Let Assumption 2.1 hold, and let

$$q_{k+1} = \sum_{t=0}^{k} \zeta_t^{(k)} \|z_{k+1} - \widehat{y}_{t+1}\|, m_{k+1} = \sum_{t=0}^{k} \zeta_t^{(k)} \|z_{k+1} - \widehat{y}_{t+1}\|^2,$$

and

$$n_{k+1} = \sum_{t=0}^{k} \zeta_t^{(k)} [g_{w_{t+1}}(\widehat{y}_{t+1}) - g(\widehat{y}_{t+1})].$$

We have , with probability 1,

(a) $q_{k+1}^2 \leq (1 - \beta_k) q_k^2 + \frac{4}{\beta_k} \|z_{k+1} - z_k\|^2$.

(b) $m_{k+1} \leq (1 - \beta_k) m_k + \beta_k q_k^2 + \frac{2}{\beta_k} \|z_{k+1} - z_k\|^2$.

(c) $\mathbb{E}[\|n_{k+1}\|^2 | \mathbb{F}_{k+1}] \leq (1 - \beta_k)^2 \|n_k\|^2 + \beta_k V_g$.

PROOF: a) By definition, we have

$$\begin{aligned}
q_{k+1} &= \sum_{t=0}^{k} \zeta_t^{(k)} \|z_{k+1} - \widehat{y}_{t+1}\| \\
&= (1 - \beta_k) \sum_{t=0}^{k-1} \zeta_t^{(k-1)} \|z_{k+1} - \widehat{y}_{t+1}\| + \beta_k \|z_{k+1} - \widehat{y}_{k+1}\| \\
&\leq (1 - \beta_k) \sum_{t=0}^{k-1} \zeta_t^{(k-1)} \Big( \|z_k - \widehat{y}_{t+1}\| + \|z_{k+1} - z_k\| \Big) + \beta_k \|z_{k+1} - \widehat{y}_{k+1}\| \\
&\leq (1 - \beta_k) q_k + (1 - \beta_k) \|z_{k+1} - z_k\| + \beta_k \|z_{k+1} - \widehat{y}_{k+1}\| \\
&= (1 - \beta_k) q_k + (1 - \beta_k) \|z_{k+1} - z_k\| + \beta_k \|z_{k+1} - \widehat{y}_{k+1}\| \\
&= (1 - \beta_k) q_k + 2(1 - \beta_k) \|z_{k+1} - z_k\|,
\end{aligned}$$

where the inequality holds by the triangle inequality, the third equality holds by the definition of $q_k$ and the fact that $\sum_{t=0}^{k} \zeta_t^{(k)} = 1$, and the last equality holds by the definition of $\widehat{y}_{k+1}$. Taking squares of both sides and using the inequality $(a + b)^2 \leq (1 + \beta)a^2 + (1 + 1/\beta)b^2$ for any $\beta > 0$, we obtain

$$\begin{aligned}
q_{k+1}^2 &\leq (1 + \beta_k)(1 - \beta_k)^2 q_k^2 + 4(1 + \beta_k^{-1})(1 - \beta_k)^2 \|z_{k+1} - z_k\|^2 \\
&\leq (1 - \beta_k) q_k^2 + 4\beta_k^{-1} \|z_{k+1} - z_k\|^2.
\end{aligned}$$

b) By the definition of $m_k$, we have

$$
\begin{aligned}
m_{k+1} &= \sum_{t=0}^{k} \zeta_t^{(k)} \|z_{k+1} - \widehat{y}_{t+1}\|^2 \\
&= (1-\beta_k) \sum_{t=0}^{k-1} \zeta_t^{(k-1)} \|z_{k+1} - \widehat{y}_{t+1}\|^2 + \beta_k \|z_{k+1} - \widehat{y}_{k+1}\|^2 \\
&= (1-\beta_k) m_k + (1-\beta_k) \sum_{t=0}^{k-1} \xi_t^{(k-1)} \Big( \|z_{k+1} - \widehat{y}_{t+1}\|^2 - \|z_k - \widehat{y}_{t+1}\|^2 \Big) + \beta_k \|z_{k+1} - \widehat{y}_{k+1}\|^2.
\end{aligned}
$$

Using the triangle inequality, we have

$$
\begin{aligned}
\|z_{k+1} - \widehat{y}_{t+1}\|^2 - \|z_k - \widehat{y}_{t+1}\|^2 \leq& \Big( \|z_{k+1} - \widehat{y}_{t+1})\| - \|z_k - \widehat{y}_{t+1}\| \Big) \Big( \|z_{k+1} - \widehat{y}_{t+1}\| + \|z_k - \widehat{y}_{t+1}\| \Big) \\
\leq& \|z_{k+1} - x_k\| \Big( 2\|z_k - \widehat{y}_{t+1}\| + \|z_{k+1} - z_k\| \Big).
\end{aligned}
$$

It follows that

$$
\begin{aligned}
m_{k+1} \leq& (1-\beta_k) m_k + (1-\beta_k) \Big( 2\|z_{k+1} - z_k\| \sum_{t=0}^{k-1} \zeta_t^{(k)} \|z_k - \widehat{y}_{t+1}\| \\
&+ \sum_{t=0}^{k-1} \zeta_t^{(k-1)} \|z_{k+1} - z_k\|^2 \Big) + \beta_k \|z_{k+1} - \widehat{y}_{k+1}\|^2 \\
=& (1-\beta_k) m_k + 2(1-\beta_k) \|z_{k+1} - z_k\| q_k + (1-\beta_k) \|z_{k+1} - z_k\|^2 \\
&+ (1-\beta_k)^2 / \beta_k \|z_{k+1} - z_k\|^2 \\
\leq& (1-\beta_k) m_k + (1-\beta_k) \Big( \beta_k^{-1} \|z_{k+1} - z_k\|^2 + \beta_k q_k^2 \Big) \\
&+ (1/\beta_k - 1) \|z_{k+1} - z_k\|^2 \\
\leq& (1-\beta_k) m_k + \beta_k q_k^2 + 2\beta_k^{-1} \|z_{k+1} - z_k\|^2,
\end{aligned}
$$

where the first equality holds by the definition of $q_k$, $\widehat{y}_{k+1}$ and the inequality $\sum_{t=0}^{k-1} \zeta_t^{(k-1)} = 1$.

c) By the definition of $n_k$, we have

$$
n_{k+1} = (1-\beta_k) n_k + \beta_k \Big( g_{w_{k+1}}(\widehat{y}_{k+1}) - g(\widehat{y}_{k+1}) \Big).
$$

Taking conditional expectations on both sides and using the fact $n_k \in \mathbb{F}_{k+1}$, we further obtain

$$
\begin{aligned}
\mathbb{E}[\|n_{k+1}\|^2 | \mathbb{F}_{k+1}] =& (1-\beta_k)^2 \|n_k\|^2 + 2\beta_k (1-\beta_k) n_k^T \mathbb{E}[g_{w_{k+1}}(\widehat{y}_{k+1}) - g(\widehat{y}_{k+1}) | \mathbb{F}_{k+1}] \\
&+ \beta_k^2 \mathbb{E}[\|g_{w_{k+1}}(\widehat{y}_{k+1}) - g(\widehat{y}_{k+1})\|^2 | \mathbb{F}_{k+1}] \\
\leq& (1-\beta_k)^2 \|n_k\|^2 + \beta_k^2 V_g,
\end{aligned}
$$

where we use the unbiasedness and moment boundedness of $g_w$ by Assumption 2.1. □

Next, we prove Lemma 3.2 by using Lemma C.1 and Lemma C.2 to construct the supermartingale with respect to an upper bound of $\{y_k - g(z_k)\}$.

Proof of Lemma 3.2: a) By Lemma C.1 and Lemma C.2, we have

$$\|y_k - g(z_k)\|^2 \le (L_g m_k + \|n_k\|)^2 \le 2L_g^2 m_k^2 + 2\|n_k\|^2.$$

By the iterative inequalities for $q_k$ and $m_k$ derived in Lemma C.2, we obtain that

$$m_{k+1} + 4q_{k+1}^2 \le (1-\beta_k/2)(m_k + 4q_k^2) + \mathcal{O}(\beta_k^{-1}\|z_{k+1} - z_k\|^2).$$

Taking squares on both sides of the above inequality and using the fact that $(a+b)^2 \le (1+\beta/2)a^2 + (1+2/\beta)b^2$ for $\beta > 0$, we have

$$(m_{k+1} + 4q_{k+1}^2)^2 \le (1-\beta_k/2)(m_k + 4q_k^2)^2 + \mathcal{O}(\beta_k^{-3}\|z_{k+1} - z_k\|^4). \tag{C.5}$$

Let

$$e_k^2 = 2L_g^2(m_k + 4q_k^2) + 2\|n_k\|^2.$$

Clearly, $\|y_k - g(z_k)\| \le e_k$ for all $k$ and $e_k \in \mathbb{F}_{k+1}$. Taking the sum of Eq.(C.5) and the iterative inequality for $n_k$ in Lemma C.2, we have

$$\mathbb{E}[e_{k+1}^2|\mathbb{F}_k] \le (1-\beta_k/2)\mathbb{E}[e_k^2|\mathbb{F}_k] + 2\beta_k^2 V_g + \mathcal{O}(\beta_k^{-3}\mathbb{E}[\|z_{k+1} - z_k\|^4|\mathbb{F}_k]).$$

b) Under the condition $\sum_{k=1}^{\infty}\gamma_k^4\beta_k^{-3} < \infty$, we have

$$\sum_{k=1}^{\infty}\beta_{j,k}^{-3}\{\mathbb{E}[\|z_{k+1} - z_k\|^4|\mathbb{F}_k]\} = \sum_{k=1}^{\infty}\beta_{j,k}^{-3}\mathbb{E}[\|z_{k+1} - z_k\|]^4 \le \sum_{k=1}^{\infty}\mathcal{O}\Big(\frac{\gamma_k^4}{\beta_k^3}\Big) < \infty,$$

where the inequality holds by the assumption $\mathbb{E}[\|z_{k+1} - z_k\|^4] \le \mathcal{O}(\gamma_k^4)$. By the monotone convergence theorem, we obtain that $\sum_{k=1}^{n}\beta_k^{-3}\mathcal{O}(\mathbb{E}[\|z_{k+1}-z_k\|^4|\mathbb{F}_k])$ converges almost surely to some random variable with finite expectation as $n \to \infty$. Therefore, the limit $\sum_{k=1}^{\infty}\beta_k^{-3}\mathcal{O}(\mathbb{E}[\|z_k - z_{k-1}\|^4|\mathbb{F}_k])$ exists and is finite with probability 1.

c) By part (a), we have that there exists a constant $C \ge 0$ such that

$$\mathbb{E}[e_{k+1}^2] \le (1-\beta_k/2)\mathbb{E}[e_k^2] + 2\beta_k^2 V_g + C\frac{\gamma_k^4}{\beta_k^3}. \tag{C.6}$$

Since $\gamma_k/\beta_k \to 0$, there exists an $M > 0$ such that $\gamma_k \le M\beta_k$ for all $k$. Letting $D_y = \mathbb{E}[e_1^2] + 4V_g + 2M^4 C$, by $\gamma_k \le M\beta_k$ and $\beta_k \le 1$, we have $D_y \ge 4V_g\beta_k + 2\beta_k^{-4}\gamma_k^4 C$ for all $k$.

We prove by induction that $\mathbb{E}[e_k^2] \le D_y$ for all $k$. Clearly, the claim holds for $k = 1$. Suppose the claim holds for $1, 2, \cdots, k$. We have by Eq.(C.6),

$$\begin{aligned}
\mathbb{E}[e_{k+1}^2] &\le (1-\frac{\beta_k}{2})\mathbb{E}[e_k^2] + 2\beta_k^2 V_g + C\frac{\gamma_k^4}{\beta_k^3} \\
&\le (1-\frac{\beta_k}{2})D_y + 2\beta_k^2 V_g + C\frac{\gamma_k^4}{\beta_k^3} \\
&\le D_y - \frac{\beta_k}{2}(D_y - 4\beta_k V_g - 2\beta_k^{-4}\gamma_k^4 C) \\
&\le D_y,
\end{aligned}$$

where the last inequality uses the fact that $D_y - 4\beta_k V_g - 2\beta_k^{-4}\gamma_k^4 C \geq 0$ for all $k$. The claim thus holds as desired.

d) By the assumption that $\mathbb{E}[\|z_{k+1} - z_k\|^4] \leq \mathcal{O}(\gamma_k^4)$, there exists a constant $C_0 > 0$ such that $\mathbb{E}[\|z_{k+1} - z_k\|^4] \leq \gamma_k^4 C_0$. By part (c), we have that there exists a constant $D_y > 0$ such that $\mathbb{E}[\|y_{k+1} - g(z_{k+1})\|^2] \leq \mathbb{E}[e_{k+1}^2] \leq D_y$ for all $k$. Thus, $\mathbb{E}[\|y_{k+1} - g(z_{k+1})\|] \leq \sqrt{\mathbb{E}[\|y_{k+1} - g(z_{k+1})\|^2]} \leq \sqrt{D_y}$. Now we begin our analysis to show the finiteness of the fourth moment $\mathbb{E}[\|y_{k+1} - g(z_{k+1})\|^4]$.

Let $e_{k+1} = (1-\beta_k)(g(z_{k+1}) - g(z_k))$, by the Lipschitz continuity of $g$, we have

$$\|e_{k+1}\| \leq (1-\beta_k)\sqrt{C_g}\|z_{k+1} - z_k\|.$$

By the definition of $y_{k+1}$, we have

$$\begin{aligned}
&y_{k+1} - g(z_{k+1}) + e_{k+1} - e_{k+1} \\
&= (1-\beta_k)(y_k - g(z_k)) + \beta_k(g_{w_{k+1}}(\widehat{y}_{k+1}) - g(z_{k+1})) - e_{k+1} \\
&= (1-\beta_k)(y_k - g(z_k)) + \beta_k(g_{w_{k+1}}(\widehat{y}_{k+1}) - g(\widehat{y}_{k+1})) + \beta_k(g(\widehat{y}_{k+1}) - g(z_{k+1})) - e_{k+1}.
\end{aligned}$$

Again, by the Lipschitz continuity of $g$, we get

$$\begin{aligned}
\beta_k\|g(\widehat{y}_{k+1}) - g(z_{k+1})\| &\leq \beta_k\sqrt{C_g}\|\widehat{y}_{k+1} - z_{k+1}\| \\
&= \beta_k\sqrt{C_g}\|\frac{1-\beta_k}{\beta_k}(z_{k+1} - z_k)\| \\
&= (1-\beta_k)\sqrt{C_g}\|z_{k+1} - z_k\|.
\end{aligned}$$

So we obtain

$$\|y_{k+1} - g(z_{k+1})\| \leq \|(1-\beta_k)(y_k - g(z_k)) + \beta_k(g_{w_{k+1}}(\widehat{y}_{k+1}) - g(\widehat{y}_{k+1}))\| + 2(1-\beta_k)\sqrt{C_g}\|z_{k+1} - z_k\|.$$

Let $P_k = \|(1-\beta_k)(y_k - g(z_k)) + \beta_k(g_{w_{k+1}}(\widehat{y}_{k+1}) - g(\widehat{y}_{k+1}))\|$, then we have

$$\begin{aligned}
P_k^4 &= \|(1-\beta_k)(y_k - g(z_k)) + \beta_k(g_{w_{k+1}}(\widehat{y}_{k+1}) - g(\widehat{y}_{k+1}))\|^4 \\
&\leq (1-\beta_k)^4\|y_k - g(z_k)\|^4 + 4(1-\beta_k)^3\beta_k(y_k - g(z_k))^3(g_{w_{k+1}}(\widehat{y}_{k+1}) - g(\widehat{y}_{k+1})) \\
&\quad + 6(1-\beta_k)^2\beta_k^2\|y_k - g(z_k)\|^2\|g_{w_{k+1}}(\widehat{y}_{k+1}) - g(\widehat{y}_{k+1})\|^2 \\
&\quad + 4(1-\beta_k)\beta_k^3\|(y_k - g(z_k))\|\|(g_{w_{k+1}}(\widehat{y}_{k+1}) - g(\widehat{y}_{k+1}))\|^3 \\
&\quad + \beta_k^4\|g_{w_{k+1}}(\widehat{y}_{k+1}) - g(\widehat{y}_{k+1})\|^4,
\end{aligned}$$

which implies

$$\begin{aligned}
\mathbb{E}[P_k^4] &\leq (1-\beta_k)^4\mathbb{E}[\|y_k - g(z_k)\|^4] + 0 \\
&\quad + 6(1-\beta_k)^2\beta_k^2\mathbb{E}\Big[\|y_k - g(z_k)\|^2\mathbb{E}[\|g_{w_{k+1}}(\widehat{y}_{k+1}) - g(\widehat{y}_{k+1})\|^2|\mathbb{F}_{k+1}]\Big] \\
&\quad + 4(1-\beta_k)\beta_k^3\mathbb{E}\Big[\|y_k - g(z_k)\|\mathbb{E}[\|(g_{w_{k+1}}(\widehat{y}_{k+1}) - g(\widehat{y}_{k+1}))\|^3|\mathbb{F}_{k+1}]\Big] \\
&\quad + \beta_k^4\mathbb{E}[\|g_{w_{k+1}}(\widehat{y}_{k+1}) - g(\widehat{y}_{k+1})\|^4] \\
&\leq (1-\beta_k)^4\mathbb{E}[\|y_k - g(z_k)\|^4] + 6(1-\beta_k)^2\beta_k^2 D_y V_g + 4\beta_k^3\sqrt{D_y}V_g^{3/2} + \beta_k^4 V_g^2.
\end{aligned}$$

Using the fact that $\|a+b\|^2 \leq (1+\epsilon)\|a\|^2 + (1+1/\epsilon)\|b\|^2$ and $\|a+b\|^4 \leq (1+\epsilon)^3\|a\|^4 + (1+1/\epsilon)^3\|b\|^4$ for $\epsilon > 0$, we obtain

$$\|y_{k+1} - g(z_{k+1})\|^4 \leq (1+\beta_k)^3 P_k^4 + 16(1+1/\beta_k)^3(1-\beta_k)^4 C_g^2 \|z_{k+1} - z_k\|^4,$$

which implies

$$\mathbb{E}[\|y_{k+1} - g(z_{k+1})\|^4] \leq (1-\beta_k)\mathbb{E}[\|y_k - g(z_k)\|^4] + 12\beta_k^2 D_y V_g + 16\beta_k^3 \sqrt{D_y} V_g^{3/2} + 8\beta_k^4 V_g^2 + \frac{16\gamma_k^4}{\beta_k^3} C_g^2 C_0.$$

Since $\gamma_k/\beta_k \to 0$, there exists a constant $M > 0$ such that $\gamma_k \leq \beta_k M$ for all $k$. We set $S_y = \|y_0 - g(z_0)\|^4 + 12 D_y V_g + 16\sqrt{D_y} V_g^{3/2} + 8V_g^2 + 16M^4 C_g^2 C_0$, and prove the claim by induction on $k$. The rest of the analysis follows the same argument as in the proof of Lemma 3.1 (a). We omit the details to avoid repetition. We conclude that there exists a constant $S_y > 0$ such that $\mathbb{E}[\|y_{k+1} - g(z_{k+1})\|^4] \leq S_y$ for all $k$.

e) By the definition of $y_k, \widehat{y}_k, z_k$, we have

$$\begin{aligned}
||y_{k+1} - y_k|| &= ||(1-\beta_{k+1})y_k + \beta_k g_{w_{k+1}}(\widehat{y}_{k+1}) - y_k|| \\
&= ||\beta_k(g_{w_{k+1}}(\widehat{y}_{k+1}) - y_k)|| = \beta_k ||g_{w_{k+1}}(\widehat{y}_{k+1}) - y_k|| \\
&= \beta_k ||g_{w_{k+1}}(\widehat{y}_{k+1}) - g(\widehat{y}_{k+1}) + g(\widehat{y}_{k+1}) - g(z_k) + g(z_k) - y_k|| \\
&\leq \beta_k \Big[ ||g_{w_{k+1}}(\widehat{y}_{k+1}) - g(\widehat{y}_{k+1})|| + ||g(\widehat{y}_{k+1}) - g(z_k)|| + ||g(z_k) - y_k|| \Big] \\
&\leq \beta_k \Big[ ||g_{w_{k+1}}(\widehat{y}_{k+1}) - g(\widehat{y}_{k+1})|| + L_g ||\widehat{y}_{k+1} - z_k|| + ||g(z_k) - y_k|| \Big] \\
&= \beta_k \Big[ ||g_{w_{k+1}}(\widehat{y}_{k+1}) - g(\widehat{y}_{k+1})|| + L_g || \frac{1}{\beta_k}(z_{k+1} - z_k)|| + ||g(z_k) - y_k|| \Big] \\
&= \beta_k ||g_{w_{k+1}}(\widehat{y}_{k+1}) - g(\widehat{y}_{k+1})|| + L_g ||z_{k+1} - z_k|| + \beta_k ||g(z_k) - y_k||.
\end{aligned}$$

Using the fact $(a+b+c)^4 \leq [2(a+b)^2 + 2c^2]^2 \leq 8(a+b)^4 + 8c^4 \leq 64a^4 + 64b^4 + 8c^4$, it is easy to see

$$||y_{k+1} - y_k||^4 \leq 64\beta_k^4 ||g_{w_{k+1}}(\widehat{y}_{k+1}) - g(\widehat{y}_{k+1})||^4 + 64L_g^4 ||z_{k+1} - z_k||^4 + 8\beta_k^4 ||g(z_k) - y_k||^4.$$

Then we get,

$$\mathbb{E}[||y_{k+1} - y_k||^2] \leq 64\beta_k^4 V_g^2 + 64L_g^2 \mathbb{E}[||z_{k+1} - z_k||^4] + 8\beta_k^4 \mathbb{E}[||g(z_k) - y_k||^4], \tag{C.7}$$

By part (d), we have that $\mathbb{E}[\|y_k - g(z_k)\|^4] \leq S_y$. Since $\mathbb{E}[\|z_{k+1} - z_k\|^4] \leq \mathcal{O}(\gamma_k^4)$ and $\gamma_k/\beta_k \to 0$, we obtain

$$\mathbb{E}[\|y_{k+1} - y_k\|^4] \leq \mathcal{O}(\beta_k^4) + \mathcal{O}(\gamma_k^4) \leq \mathcal{O}(\beta_k^4),$$

which concludes the proof.

□

### C.3 Proof of Theorem 3.1

PROOF OF THEOREM 3.1: (a) Let $x^*$ be an arbitrary optimal solution to problem (1.1), and let $F^* = F(x^*)$. By the same argument as in the proof of Lemma 2.1, we obtain there exists $C_0 > 0$ such that

$$\begin{aligned}
&\mathbb{E}[\|x_{k+1} - x^*\|^2 | \mathbb{F}_k] \\
\le&\Big(1 + C_0(\frac{\alpha_k^2}{\beta_{T-1,k}} + \frac{\alpha_k^2}{\beta_{T-2,k}} + \cdots + \frac{\alpha_k^2}{\beta_{1,k}})\Big)\|x_k - x^*\|^2 + \alpha_k^2 C_1 C_2 \cdots C_T - 2\alpha_k(F(x_k) - F^*) \\
&+ (T-1)\beta_{T-1,k}\mathbb{E}[\|y_k^{(T-1)} - f^{(T)}(x_k)\| | \mathbb{F}_k] + (T-2)\beta_{T-2,k}\mathbb{E}[\|y_k^{(T-2)} - f^{(T-1)}(y_k^{(T-1)})\|^2 | \mathbb{F}_k] \\
&+ \cdots + \beta_{1,k}\mathbb{E}[\|y_k^{(1)} - f^{(2)}(y_k^{(2)})\|^2 | \mathbb{F}_k].
\end{aligned} \tag{C.8}$$

First, we consider the case when the first inner level function $f^{(T)}$ does not have Lipschitz continuous gradients. In this case, Algorithm 2 runs with basic update step for the first inner level and accelerated update steps for all other levels. By Assumption 3.1, we have

$$\mathbb{E}[\|x_{k+1} - x_k\|^4] \le \alpha_k^4 C_1^2 C_2^2 \cdots C_T^2,$$

which is the sufficient condition for Lemma 3.1 to be true. Apply Lemmas 2.3 and 3.1 to the first update step, we have

$$\begin{aligned}
&\mathbb{E}\Big[\mathbb{E}[\|y_{k+1}^{(T-1)} - f^{(T)}(x_{k+1})\|^2 | \mathbb{F}_{k+1}]\Big|\mathbb{F}_k\Big] \\
\le&(1 - \beta_{T-1,k})\mathbb{E}[\|y_{k-1}^{(T-1)} - f^{(T)}(x_{k-1})\|^2 | \mathbb{F}_k] + \beta_{T-1}^{-1} C_T \mathbb{E}[\|x_k - x_{k-1}\|^2 | \mathbb{F}_k] + 2V_T\beta_{T-1,k}^2,
\end{aligned} \tag{C.9}$$

and $\mathbb{E}[\|y_{k+1}^{(T-1)} - y_k^{(T-1)}\|^4] \le \mathcal{O}(\beta_{T-1,k}^4)$, which serves as the sufficient condition for level $(T-1)$ in Lemma 3.2 to be true so that $\mathbb{E}[\|y_{k+1}^{(T-2)} - y_k^{(T-2)}\|^4] \le \mathcal{O}(\beta_{T-2,k}^4)$.

Apply Lemma 3.2 recursively for the accelerated update from $j = T-2$ to 1, we have that for all $j$'s and all $k$, there exists an $e_k^{(j)} \in \mathbb{F}_{k+1}$, such that almost surely

$$\mathbb{E}[\|y_k^{(j)} - f^{(j+1)}(y_k^{(j+1)})\|^2 | \mathbb{F}_k] \le \mathbb{E}\Big[[e_k^{(j)}]^2 \Big| \mathbb{F}_k\Big],$$

$$\mathbb{E}\Big[\mathbb{E}[[e_{k+1}^{(j)}]^2 | \mathbb{F}_{k+1}] | \mathbb{F}_k\Big] \le (1 - \beta_{j,k}/2)\mathbb{E}[[e_k^{(j)}]^2 | \mathbb{F}_k] + 2\beta_{j,k}^2 V_{j+1} + \mathcal{O}\Big(\frac{\mathbb{E}[\|y_{k+1}^{(j+1)} - y_k^{(j+1)}\|^4 | \mathbb{F}_k]}{\beta_{j,k}^3}\Big), \tag{C.10}$$

and $\mathbb{E}[\|y_{k+1}^{(j)} - y_k^{(j)}\|^4] \le \mathcal{O}(\beta_{j,k}^4)$.

By Lemma 3.2 part (b), under the condition that $\sum_{k=1}^{\infty} \frac{\beta_{j+1,k}^4}{\beta_{j,k}^3} < \infty$, we have

$$\sum_{k=1}^{\infty} \frac{\mathbb{E}[\|y_{k+1}^{(j+1)} - y_k^{(j+1)}\|^4 | \mathbb{F}_k]}{\beta_{j,k}^3} < \infty,$$

with probability 1. Together with the condition $\sum_{k=1}^{\infty} \beta_{j,k}^2 < \infty$, the sum of tail part of this super-martingale, $2\beta_{j,k}^2 V_{j+1} + \mathcal{O}(\frac{\mathbb{E}[\|y_{k+1}^{(j+1)} - y_k^{(j+1)}\|^4 | \mathbb{F}_k]}{\beta_{j,k}^3})$, converges almost surely.

Similarly, by Lemma 2.3 part (b), under the condition $\sum_{k=1}^{\infty} \frac{\alpha_k^2}{\beta_{T-1,k}} < \infty$, we have that with probability 1,

$$\sum_{k=1}^{\infty} \frac{\mathbb{E}[\|x_{k+1} - x_k\|^2 | \mathbb{F}_k]}{\beta_{T-1,k}} < \infty.$$

Together with the condition $\sum_{k=1}^{\infty} \beta_{T-1,k}^2 < \infty$, the sum of the tail part of the super-martingale for Eq.(C.9), $2V_T\beta_k^2 + \frac{C_T \mathbb{E}[\|x_{k+1}-x_k\|^2|\mathbb{F}_k]}{\beta_{T-1,k}}$, converges almost surely.

Now we apply the $T$-element super-martingale convergent lemma to Eqs.(C.9), (C.8), and (C.10). By letting

$$\begin{aligned}
&X_k = \|x_k - x^*\|^2, Y_k^{(T-1)} = \mathbb{E}[\|y_k^{(T-1)} - f^{(T)}(x_k)\|^2 | \mathbb{F}_k],\\
&Y_k^{(T-2)} = \mathbb{E}[[e_k^{(T-2)}]^2 | \mathbb{F}_k], \cdots, Y_k^{(1)} = \mathbb{E}[[e_k^{(1)}]^2 | \mathbb{F}_k],\\
&\eta_k = [\frac{\alpha_k^2}{\beta_{T-1,k}} + \cdots + \frac{\alpha_k^2}{\beta_{1,k}}]C_0, u_k^{(T)} = 2\alpha_k(F(x_k) - F^*),\\
&u_k^{(1)} = u_k^{(2)} = \cdots = u_k^{(T-1)} = 0, c_1 = 2, \cdots, c_{T-2} = 2(T-2), c_{T-1} = T-1,\\
&\mu_k^{(1)} = 2\beta_{1,k}^2 V_1 + \mathcal{O}(\frac{\mathbb{E}[\|y_{k+1}^{(2)} - y_k^{(2)}\|^4 | \mathbb{F}_k]}{\beta_{1,k}^3}), \cdots,\\
&\mu_k^{(T-2)} = 2\beta_{T-2,k}^2 V_{T-1} + \mathcal{O}(\frac{\mathbb{E}[\|y_{k+1}^{(T-1)} - y_k^{(T-1)}\|^4 | \mathbb{F}_k]}{\beta_{T-2,k}^3}),\\
&\mu_k^{(T-1)} = C_T \beta_{T-1,k}^{-1} \mathbb{E}[\|x_{k+1} - x_k\|^2 | \mathbb{F}_k] + 2V_T \beta_{T-1,k}^2,\\
&\mu_k^{(T)} = \alpha_k^2 C_1 C_2 \cdots C_T,\\
&\theta_k^{(1)} = \beta_{1,k}/2, \cdots, \theta_k^{(T-2)} = \beta_{T-2,k}/2, \theta_k^{(T-1)} = \beta_{T-1,k},
\end{aligned}$$

we obtain that $\|x_{k+1} - x^*\|$ converges almost surely to a nonnegative random variable, and

$$\sum_{k=0}^{\infty} \alpha_k (F(x_k) - F^*) < \infty,$$

which further implies

$$\liminf_{k \to \infty} F(x_k) = F^*, \qquad w.p.1.$$

Using Lemma 2.5, we conclude that the sequence $\{x_k\}$ converges almost surely to a random point in the set of optimal solutions to problem (1.1).

Next, we study the case when Assumption 3.2 also holds, i.e.,$f^{(T)}$ has Lipschitz continuous gradient. In this case, Algorithm 2 runs with accelerating update steps for all levels. The only difference is that we use the accelerated update rule for the first inner level instead of the basic one, so Lemma 2.1 and 3.2 still hold. Consider the first inner level, since it is also updated by the accelerated update rule, we apply similar analysis as in Lemma 3.2 for this level. We have that there exists $e_k^{(T-1)} \in \mathbb{F}_{k+1}$ such that with probability 1,

$$\mathbb{E}[\|y_k^{(T-1)} - f^{(T)}(x_k)\|^2 | \mathbb{F}_k] \leq \mathbb{E}\Big[[e_k^{(T-1)}]^2 \Big| \mathbb{F}_k\Big],$$

and

$$\mathbb{E}\Big[\mathbb{E}[[e_{k+1}^{(T-1)}]^2|\mathbb{F}_{k+1}]\Big|\mathbb{F}_k\Big] \le (1-\frac{\beta_{j,k}}{2})\mathbb{E}[[e_k^{(T-1)}]^2|\mathbb{F}_k] + 2\beta_{j,k}^2 V_{j+1} + \mathcal{O}(\frac{\mathbb{E}[\|x_{k+1}-x_k\|^4|\mathbb{F}_k]}{\beta_{j,k}^3}). \quad (C.11)$$

By similar argument as in Lemma 3.2 part (b), under the condition $\sum_{k=1}^{\infty}\frac{\alpha_k^4}{\beta_{T-1,k}^3} < \infty$, we have with probability 1,

$$\sum_{k=1}^{\infty}\frac{\mathbb{E}[\|x_{k+1}-x_k\|^4|\mathbb{F}_k]}{\beta_{j,k}^3} < \infty.$$

Now we apply the $T$-element super-martingale convergent lemma to Eqs.(C.8),(C.10) and (C.11). The remaining part follows the same line as in the case where $f^{(T)}$ does not have Lipschitz continuous gradient. We conclude $\{x_k\}$ converges almost surely to a random optimal solution.

(b) Firstly, consider the case when $f^{(T)}$ does not have Lipschitz continuous gradient. Since problem (1.1) has at least one optimal solution, the function $F$ is bounded from below, and denote $F^*$ as the optimal value of $F(x)$ over $\mathcal{X}$. As a result, we can treat $\{F(x_k)-F^*\}$ as nonnegative random variables. The $x_k$ update steps of Algorithm 2 and Algorithm 1 are same. Thus, Lemma 2.2 holds for Algorithm 2 as well. It follows that for sufficiently large $k$,

$$\begin{aligned}&\mathbb{E}[F(x_{k+1})-F^*|\mathbb{F}_k]\\ \le& F(x_k)-F^*-\frac{\alpha_k}{2}\|\nabla F(x_k)\|^2+\frac{1}{2}\alpha_k^2 L_F C_1C_2\cdots C_T+(T-1)\beta_{T-1,k}\mathbb{E}[\|y_k^{(T-1)}-f^{(T)}(x_k)\|^2|\mathbb{F}_k]\\ &+\cdots+\beta_{1,k}\mathbb{E}[\|y_k^{(1)}-f^{(2)}(y_k^{(2)})\|^2|\mathbb{F}_k].\end{aligned} \quad (C.12)$$

Using similar argument as in part (a), we apply the $T$-element super-martingale convergence lemma to Eqs.(C.10) , (C.11) and (C.12). By letting

$$\begin{aligned}&X_k = F(x_k)-F^*, Y_k^{(T-1)} = \mathbb{E}[\|y_k^{(T-1)}-f^{(T)}(x_k)\|^2|\mathbb{F}_k],\\ &Y_k^{(T-2)} = \mathbb{E}[[e_k^{(T-2)}]^2|\mathbb{F}_k],\cdots, Y_k^{(1)} = \mathbb{E}[[e_k^{(1)}]^2|\mathbb{F}_k],\\ &\eta_k = 0, u_k^{(T)} = \frac{1}{2}\alpha_k\|\nabla F(x)\|^2,\\ &u_k^{(1)} = u_k^{(2)} = \cdots = u_k^{(T-1)} = 0, c_1 = 2,\cdots, c_{T-2} = 2(T-2), c_{T-1} = T-1,\\ &\mu_k^{(T-1)} = C_T\beta_{T-1,k}^{-1}\mathbb{E}[\|x_{k+1}-x_k\|^2|\mathbb{F}_k]+2V_T\beta_{T-1,k}^2,\\ &\mu_k^{(T-2)} = 2\beta_{T-2,k}^2V_{T-1}+\mathcal{O}(\frac{\mathbb{E}[\|y_{k+1}^{(T-1)}-y_k^{(T-1)}\|^4|\mathbb{F}_k]}{\beta_{T-2,k}^3}),\cdots,\\ &\mu_k^{(1)} = 2\beta_{1,k}^2V_1+\mathcal{O}(\frac{\mathbb{E}[\|y_{k+1}^{(2)}-y_k^{(2)}\|^4|\mathbb{F}_k]}{\beta_{1,k}^3}),\\ &\mu_k^{(T)} = \frac{1}{2}\alpha_k^2L_FC_1C_2\cdots C_T,\\ &\theta_k^{(1)} = \beta_{1,k}/2,\cdots,\theta_k^{(T-2)} = \beta_{T-2,k}/2, \theta_k^{(T-1)} = \beta_{T-1,k},\end{aligned}$$

we obtain $\{F(x_k) - F^*\}$ converges almost surely to a nonnegative random variable, and

$$\sum_{k=0}^{\infty} \alpha_k \|\nabla F(x_k)\|^2 < \infty, \ \ w.p.1.$$

By Lemma 2.6, we conclude any limit point of the sequence $\{x_k\}$ is a stationary point with probability 1, which completes the proof.

When $f^{(T)}$ has Lipschitz continuous gradient, we apply the $T$-element super-martingale convergent lemma to Eqs.(C.10), (C.11), and (C.12). The rest of the proof follows the same line as in the case where $f^{(T)}$ is non-smooth. □

# D Proof of Theorem 3.2

## D.1 Proof of Lemma 3.3

Proof: By Lemma 3.2, we have that for $j = T-2, \cdots, 1$, there exists a random variables $e_k^{(j)} \in \mathbb{F}_{k+1}$ for all $k$ satisfying $\|y_k^{(j)} - f^{(j+1)}(y_k^{(j+1)})\| \le e_k^{(j)}$ such that

$$\mathbb{E}[[e_{k+1}^{(j)}]^2|\mathbb{F}_k] \le (1-\beta_{j,k}/2)[e_k^{(j)}]^2 + 2\beta_{j,k}^2 V_{j+1} + \mathcal{O}\Big(\frac{\mathbb{E}[\|y_{k+1}^{(j+1)} - y_k^{(j+1)}\|^4|\mathbb{F}_k]}{\beta_{j,k}^3}\Big),$$

almost surely and

$$\mathbb{E}[\|y_{k+1}^{(j)} - y_k^{(j)}\|^4] \le \mathcal{O}(\beta_{j,k}^4).$$

So we have

$$\mathbb{E}[[e_{k+1}^{(j)}]^2] \le (1-\beta_{j,k}/2)\mathbb{E}[[e_k^{(j)}]^2] + 2\beta_{j,k}^2 V_{j+1} + \mathcal{O}\Big(\frac{\beta_{j+1,k}^4}{\beta_{j,k}^3}\Big). \tag{D.1}$$

Substituting $\beta_{j,k} = 2k^{-b_j}$ into Eq.(D.1) and applying Lemma B.1, we obtain

$$\mathbb{E}[\|y_k^{(j)} - f^{(j+1)}(y_k^{(j+1)})\|^2] \le \mathbb{E}[[e_k^{(j)}]^2] \le \mathcal{O}(k^{-4b_{j+1}+4b_j}) + \mathcal{O}(k^{-b_j}),$$

which completes the proof. □

## D.2 Proof of Lemma 3.4

Proof: By the Lipschitz continuous gradient condition in Assumption 2.2, we have

$$\begin{aligned}
&F(x_{k+1}) - F(x_k)\\
\le&\langle \nabla F(x_k), x_{k+1} - x_k\rangle + \frac{L_F}{2}\|x_{k+1} - x_k\|^2\\
\le& -\alpha_k \langle \nabla F(x_k), \widetilde{\nabla} f_{\omega_{T,k}}^{(T)}(x_k) \nabla f_{\omega_{T-1,k}}^{(T-1)}(y_k^{(T-1)}) \cdots \nabla f_{\omega_{1,k}}^{(1)}(y_k^{(1)})\rangle + \mathcal{O}(\alpha_k^2)\\
\le& -\alpha_k \|\nabla F(x_k)\|^2 + \alpha_k \langle \nabla F(x_k), \nabla F(x_k) - \widetilde{\nabla} f_{\omega_{T,k}}^{(T)}(x_k) \nabla f_{\omega_{T-1,k}}^{(T-1)}(y_k^{(T-1)}) \cdots \nabla f_{\omega_{1,k}}^{(1)}(y_k^{(1)})\rangle + \mathcal{O}(\alpha_k^2).
\end{aligned} \tag{D.2}$$

Let $Q$ be

$$Q = \Big\langle \nabla F(x_k), \nabla F(x_k) - \widetilde{\nabla} f_{\omega_{T,k}}^{(T)}(x_k) \nabla f_{\omega_{T-1,k}}^{(T-1)}(y_k^{(T-1)}) \cdots \nabla f_{\omega_{1,k}}^{(1)}(y_k^{(1)})\Big\rangle.$$

We have

$$\begin{aligned}\mathbb{E}[Q] =&\mathbb{E}\left[\langle \nabla F(x_k), \nabla F(x_k) - \widetilde{\nabla} f^{(T)}_{\omega_{T,k}}(x_k) \nabla f^{(T-1)}_{\omega_{T-1,k}}(y_k^{(T-1)}) \cdots \nabla f^{(1)}_{\omega_{1,k}}(y_k^{(1)})\rangle\right]\\ =&\mathbb{E}\left[\langle \nabla F(x_k), \widetilde{\nabla} F_{\omega_k}(x_k) - \widetilde{\nabla} f^{(T)}_{\omega_{T,k}}(x_k) \nabla f^{(T-1)}_{\omega_{T-1,k}}(y_k^{(T-1)}) \cdots \nabla f^{(1)}_{\omega_{1,k}}(y_k^{(1)})\rangle\right],\end{aligned}$$

where $\widetilde{\nabla} F_{\omega_k}(x_k) = \widetilde{\nabla} f^{(T)}_{\omega_{T,k}}(x_k) \nabla f^{(T-1)}_{\omega_{T-1,k}}(f^{(T)}(x_k)) \cdots \nabla f^{(1)}_{\omega_{1,k}}(f^{(2)} \circ \cdots \circ (f^{(T)}(x_k)))$ and the equality comes from Assumption 2.1 (ii). Based on the fact $2ab \le a^2 + b^2$ for all $a, b$, we obtain

$$\begin{aligned}&\mathbb{E}\Big[\langle \nabla F(x_k), \widetilde{\nabla} F_{\omega_k}(x_k) - \widetilde{\nabla} f^{(T)}_{\omega_{T,k}}(x_k) \nabla f^{(T-1)}_{\omega_{T-1,k}}(y_k^{(T-1)}) \cdots \nabla f^{(1)}_{\omega_{1,k}}(y_k^{(1)})\rangle\Big]\\ \le&\frac{1}{2}\mathbb{E}[\|\nabla F(x_k)\|^2] + \frac{1}{2}\mathbb{E}[\|\widetilde{\nabla} F_{\omega_k}(x_k) - \widetilde{\nabla} f^{(T)}_{\omega_{T,k}}(x_k) \nabla f^{(T-1)}_{\omega_{T-1,k}}(y_k^{(T-1)}) \cdots \nabla f^{(1)}_{\omega_{1,k}}(y_k^{(1)})\|^2].\end{aligned}$$

Applying the inequality $\|a+b\|^2 \le (\|a\| + \|b\|)^2 \le 2\|a\|^2 + 2\|b\|^2$ to Eq.(A.1)-(A.4) in Lemma A.1, we have

$$\begin{aligned}&\frac{1}{2}\mathbb{E}[\|\widetilde{\nabla} F_{\omega_k}(x_k) - \widetilde{\nabla} f^{(T)}_{\omega_{T,k}}(x_k) \nabla f^{(T-1)}_{\omega_{T-1,k}}(y_k^{(T-1)}) \cdots \nabla f^{(1)}_{\omega_{1,k}}(y_k^{(1)})\|^2]\\ \le&\mathcal{O}(\mathbb{E}[\|y_k^{(T-1)} - f^{(T)}(x_k)\|^2]) + \cdots + \mathcal{O}(\mathbb{E}[\|y_k^{(1)} - f^{(2)}(y_k^{(2)})\|^2]).\end{aligned}$$

Thus, we obtain

$$\begin{aligned}\mathbb{E}[Q] =&\mathbb{E}\big[\langle \nabla F(x_k), \widetilde{\nabla} F_{\omega_k}(x_k) - \widetilde{\nabla} f^{(T)}_{\omega_{T,k}}(x_k) \nabla f^{(T-1)}_{\omega_{T-1,k}}(y_k^{(T-1)}) \cdots \nabla f^{(1)}_{\omega_{1,k}}(y_k^{(1)})\rangle\big]\\ \le&\frac{1}{2}\mathbb{E}[\|\nabla F(x_k)\|^2] + \mathcal{O}(\mathbb{E}[\|y_k^{(T-1)} - f^{(T)}(x_k)\|^2]) + \mathcal{O}(\mathbb{E}[\|y_k^{(T-2)} - f^{(T-1)}(y_k^{(T-1)})\|^2])\\ &+ \cdots + \mathcal{O}(\mathbb{E}[\|y_k^{(1)} - f^{(2)}(y_k^{(2)})\|^2]).\end{aligned}$$

Taking expectations on both sides of Eq.(D.2) and substituting $\mathbb{E}[Q]$ by its upper bound derived above, we have

$$\begin{aligned}&\frac{\alpha_k}{2}\|\nabla F(x_k)\|^2\\ \le&\mathbb{E}[F(x_k)] - \mathbb{E}[F(x_{k+1})] + \mathcal{O}(\alpha_k \mathbb{E}[\|y_k^{(T-1)} - f^{(T)}(x_k)\|^2]) + \mathcal{O}(\alpha_k \mathbb{E}[\|y_k^{(T-2)} - f^{(T-1)}(y_k^{(T-1)})\|^2])\\ &+ \cdots + \mathcal{O}(\alpha_k \mathbb{E}[\|y_k^{(1)} - f^{(2)}(y_k^{(2)})\|^2]) + \mathcal{O}(\alpha_k^2).\end{aligned}$$

This implies that

$$\begin{aligned}&\mathbb{E}[\|\nabla F(x_k)\|^2]\\ \le&\ 2\alpha_k^{-1}\mathbb{E}[F(x_k)] - 2\alpha_k^{-1}\mathbb{E}[F(x_{k+1})] + \mathcal{O}(\mathbb{E}[\|y_k^{(T-1)} - f^{(T)}(x_k)\|^2])\\ &+ \mathcal{O}(\mathbb{E}[\|y_k^{(T-2)} - f^{(T-1)}(y_k^{(T-1)})\|^2]) + \cdots + \mathcal{O}(\mathbb{E}[\|y_k^{(1)} - f^{(2)}(y_k^{(2)})\|^2]) + \mathcal{O}(\alpha_k),\end{aligned} \tag{D.3}$$

which completes the proof. □

## D.3 Proof of Theorem 3.2

PROOF OF THEOREM 3.2: Firstly, we consider the case where $f^{(T)}$ does not have Lipschitz continuous gradient. Since Assumption 2.1 and 2.2 hold, we apply Lemma 3.4 and sum up Eq.(D.3) from $k = 1$ to $n$ , then we obtain

$$\begin{aligned}
&\frac{\sum_{k=1}^{n}\mathbb{E}(\|\nabla(x_k)\|^2)}{n} \\
&\leq 2n^{-1}\alpha_1^{-1}F(x_0) + n^{-1}\sum_{k=1}^{n}((k+1)^a - k^a)\mathbb{E}[F(x_k)] + n^{-1}\sum_{k=1}^{n}\mathcal{O}(\mathbb{E}[\|y_k^{(T-1)} - f^{(T)}(x_k)\|^2]) + \\
&\quad \cdots + K^{-1}\sum_{k=1}^{n}\mathcal{O}(\mathbb{E}[\|y_k^{(1)} - f^{(2)}(y_k^{(2)})\|^2]) + n^{-1}\sum_{k=1}^{n}\mathcal{O}(\alpha_k) \\
&\leq 2n^{-1}F(x_0) + n^{-1}\sum_{k=1}^{n}ak^{a-1}\mathbb{E}[F(x_k)] + n^{-1}\sum_{k=1}^{n}\mathcal{O}(\mathbb{E}[\|y_k^{(T-1)} - f^{(T)}(x_k)\|^2]) + \cdots \\
&\quad + n^{-1}\sum_{k=1}^{n}\mathcal{O}(\mathbb{E}[\|y_k^{(1)} - f^{(2)}(y_k^{(2)})\|^2]) + n^{-1}\sum_{k=1}^{n}\mathcal{O}(n^{-a}),
\end{aligned} \tag{D.4}$$

where the second inequality holds by the fact $(k+1)^a \leq k^a + ak^{a-1}$ since $h(t) = t^a$ is a concave function for $0 < a < 1$.

Meanwhile, by Lemma 2.7 and Lemma 3.3, with the choice of $a = \frac{4+T}{8+T}$ and $b_j = \frac{3+j}{8+T}$ for $j = T-1, T-2, \cdots, 1$, we have $\mathbb{E}[\|y_{k+1}^{(T-1)} - f^{(T)}(x_k)\|^2] \leq \mathcal{O}(k^{-4/(8+T)})$ and $\mathbb{E}[\|y_{k+1}^{(j)} - f^{(j+1)}(y_{k+1}^{(j+1)})\|^2] \leq \mathcal{O}(k^{-4/(8+T)})$ for $j = T-2, \cdots, 1$. Plug them into Eq.(D.4), we have

$$\begin{aligned}
&\frac{\sum_{k=1}^{n}\mathbb{E}(\|\nabla F(x_k)\|^2)}{n} \\
\leq &\mathcal{O}\Big(n^{a-1} + n^{2(b_{T-1}-a)}\mathbb{I}^{\log n}_{2(a-b_{T-1})=1} + n^{-b_{T-1}} + \sum_{j=1}^{T-2}[n^{4(b_j-b_{j+1})}\mathbb{I}^{\log n}_{4(b_{j+1}-b_j)=1} + n^{-b_j}] + n^{-a}\Big) \\
\leq &\mathcal{O}(n^{-4/(8+T)}),
\end{aligned} \tag{D.5}$$

which completes the proof.

Next, when $f^{(T)}$ has Lipschitz continuous gradient, the first inner level is also updated by the accelerating rule. By similar analysis as in Lemma 3.3, we have

$$\mathbb{E}[\|y_k^{(T-1)} - f^{(T)}(x_k)\|^2] \leq \mathcal{O}(k^{-4a+4b_{T-1}}) + \mathcal{O}(k^{-b_{T-1}}).$$

Plug this convergent rate into Eq.(D.5), together with Lemma 3.3, we obtain

$$\begin{aligned}
&\frac{\sum_{k=1}^{n}\mathbb{E}(\|\nabla F(x_k)\|^2)}{n} \\
\leq &\mathcal{O}\Big(n^{a-1} + n^{-4a+4b_{T-1}}\mathbb{I}^{\log n}_{4(a-b_{T-1})=1} + n^{-b_{T-1}} + \sum_{j=1}^{T-2}[n^{4(b_j-b_{j+1})}\mathbb{I}^{\log n}_{4(b_{j+1}-b_j)} + n^{-b_j}] + n^{-a}\Big) \\
\leq &\mathcal{O}(n^{-4/(7+T)}),
\end{aligned}$$

by choosing $a = \frac{3+T}{7+T}$ and $b_j = \frac{3+j}{7+T}$ for $j = T-1, T-2, \cdots, 1$, which completes the proof. □

## E Proof of Theorem 3.3

PROOF OF THEOREM 3.3: Denote $F^* = \min_{x\in\mathcal{X}} F(x)$, note that $F^* = F(\Pi_{\mathcal{X}^*}(x))$ for all $x \in \mathcal{X}$. When $\mathcal{X} = \mathbb{R}^{d_T}$, based on the definition of $x_k$, we have

$$x_{k+1} - x_k = -\alpha_k \widetilde{\nabla} f^{(T)}_{\omega_{T,k}}(x_k) \nabla f^{(T-1)}_{\omega_{T-1,k}}(y_k^{(T-1)}) \cdots \nabla f_{\omega_{1,k}}(y_k^{(1)}).$$

Then, for the term $||x_{k+1} - \Pi_{\mathcal{X}^*}(x_{k+1})||^2$, we have

$$\begin{aligned}
&||x_{k+1} - \Pi_{\mathcal{X}^*}(x_{k+1})||^2 \\
\leq&||x_{k+1} - \Pi_{\mathcal{X}^*}(x_k)||^2 \\
\leq&||x_{k+1} - x_k + x_k - \Pi_{\mathcal{X}^*}(x_k)||^2 \\
=&||x_k - \Pi_{\mathcal{X}^*}(x_k)||^2 - ||x_{k+1} - x_k||^2 + 2\langle x_{k+1} - x_k, x_{k+1} - \Pi_{\mathcal{X}^*}(x_k)\rangle \\
=&||x_k - \Pi_{\mathcal{X}^*}(x_k)||^2 - ||x_{k+1} - x_k||^2 \\
&+ 2\alpha_k \langle \widetilde{\nabla} f^{(T)}_{\omega_{T,k}}(x_k) \nabla f^{(T-1)}_{\omega_{T-1,k}}(y_k^{(T-1)}) \cdots \nabla f^{(1)}_{\omega_{1,k}}(y_k^{(1)}), \Pi_{\mathcal{X}^*}(x_k) - x_{k+1}\rangle \\
\leq&||x_k - \Pi_{\mathcal{X}^*}(x_k)||^2 - ||x_{k+1} - x_k||^2 + 2\alpha_k \langle \nabla F(x_k), \Pi_{\mathcal{X}^*}(x_k) - x_{k+1}\rangle \\
&+ 2\alpha_k \langle \widetilde{\nabla} f^{(T)}_{\omega_{T,k}}(x_k) \nabla f^{(T-1)}_{\omega_{T-1,k}}(y_k^{(T-1)}) \cdots \nabla f^{(1)}_{\omega_{1,k}}(y_k^{(1)}) - \nabla F(x_k), \Pi_{\mathcal{X}^*}(x_k) - x_{k+1}\rangle,
\end{aligned} \tag{E.1}$$

where the second equality comes from $\|a+b\|^2 = \|b\|^2 - \|a\|^2 + 2\langle a, a+b\rangle$ with $a = x_{k+1} - x_k$ and $b = x_k - \Pi_{\mathcal{X}^*}(x_k)$. Define

$$T_1 = \langle \nabla F(x_k), \Pi_{\mathcal{X}^*}(x_k) - x_{k+1}\rangle,$$

and

$$T_2 = \langle \widetilde{\nabla} f^{(T)}_{\omega_{T,k}}(x_k) \nabla f^{(T-1)}_{\omega_{T-1,k}}(y_k^{(T-1)}) \cdots \nabla f^{(1)}_{\omega_{1,k}}(y_k^{(1)}) - \nabla F(x_k), \Pi_{\mathcal{X}^*}(x_k) - x_{k+1}\rangle.$$

For the term $T_1$, we have

$$\begin{aligned}
T_1 =&\langle \nabla F(x_k), x_k - x_{k+1}\rangle + \langle \nabla F(x_k), \Pi_{\mathcal{X}^*}(x_k) - x_k\rangle \\
\leq& F(x_k) - F(x_{k+1}) + \frac{L_F}{2}\|x_{k+1} - x_k\|^2 + \langle \nabla F(x_k), \Pi_{\mathcal{X}^*}(x_k) - x_k\rangle \\
\leq& F(x_k) - F(x_{k+1}) + \frac{L_F}{2}\|x_{k+1} - x_k\|^2 + F(\Pi_{\mathcal{X}^*}(x_k)) - F(x_k) \\
=& F(\Pi_{\mathcal{X}^*}(x_k)) - F(x_{k+1}) + \frac{L_F}{2}\|x_{k+1} - x_k\|^2 \\
=& F^* - F(x_{k+1}) + \frac{L_F}{2}\|x_{k+1} - x_k\|^2 \\
\leq& F^* - F(x_{k+1}) + \mathcal{O}(\alpha_k^2)
\end{aligned}$$

where the first inequality is due to the Lipschitz gradient of $F$, the second inequality comes from the convexity of $F$, and the last inequality holds by Eq.(A.8).

For the term $T_2$, we have

$$\begin{aligned}\mathbb{E}[T_2] =&\mathbb{E}[\langle \widetilde{\nabla} F_{\omega_k}(x_k) - \widetilde{\nabla} f^{(T)}_{\omega_{T,k}}(x_k)\nabla f^{(T-1)}_{\omega_{T-1,k}}(y_k^{(T-1)})\cdots\nabla f^{(1)}_{\omega_{1,k}}(y_k^{(1)}), x_k - \Pi_{\mathcal{X}^*}(x_k)\rangle] \\ &+ \mathbb{E}[\langle \widetilde{\nabla} F_{\omega_k}(x_k) - \widetilde{\nabla} f^{(T)}_{\omega_{T,k}}(x_k)\nabla f^{(T-1)}_{\omega_{T-1,k}}(y_k^{(T-1)})\cdots\nabla f^{(1)}_{\omega_{1,k}}(y_k^{(1)}), x_{k+1} - x_k\rangle] \\ \le& \underbrace{\mathbb{E}[\langle \widetilde{\nabla} F_{\omega_k}(x_k) - \widetilde{\nabla} f^{(T)}_{\omega_{T,k}}(x_k)\nabla f^{(T-1)}_{\omega_{T-1,k}}(y_k^{(T-1)})\cdots\nabla f^{(1)}_{\omega_{1,k}}(y_k^{(1)}), x_k - \Pi_{\mathcal{X}^*}(x_k)\rangle]}_{T_{2,1}} \\ &+ \frac{\alpha_k}{2}\underbrace{\mathbb{E}[\|\widetilde{\nabla} F_{\omega_k}(x_k) - \widetilde{\nabla} f^{(T)}_{\omega_{T,k}}(x_k)\nabla f^{(T-1)}_{\omega_{T-1,k}}(y_k^{(T-1)})\cdots\nabla f^{(1)}_{\omega_{1,k}}(y_k^{(1)})\|^2]}_{T_{2,2}} \\ &+ \frac{1}{2\alpha_k}\|x_k - x_{k+1}\|^2,\end{aligned}$$

where the inequality comes from the fact that $\langle a, b\rangle \le \frac{1}{2\alpha_k}\|a\|^2 + \frac{\alpha_k}{2}\|b\|^2$. For $T_{2,1}$, we have

$$\begin{aligned}T_{2,1} \le& \frac{\alpha_k}{2\phi_k}\mathbb{E}\Big[\|\widetilde{\nabla} F_{\omega_k}(x_k) - \widetilde{\nabla} f^{(T)}_{\omega_{T,k}}(x_k)\nabla f^{(T-1)}_{\omega_{T-1,k}}(y_k^{(T-1)})\cdots\nabla f^{(1)}_{\omega_{1,k}}(y_k^{(1)})\|^2\Big] + \mathbb{E}[\frac{\phi_k}{2\alpha_k}\|x_k - \Pi_{\mathcal{X}^*}(x_k)\|^2] \\ \le& \mathcal{O}(\frac{\alpha_k}{\phi_k})\Big[\mathbb{E}[\|y_k^{(T-1)} - f^{(T)}(x_k)\|^2] + \cdots + \mathbb{E}[\|y_k^{(1)} - f^{(2)}(y_k^{(2)})\|^2]\Big] + \frac{\phi_k}{2\alpha_k}\mathbb{E}[\|x_k - \Pi_{\mathcal{X}^*}(x_k)\|^2],\end{aligned}$$

where $\phi_k$ is a scalar and will be specified later. By the fact $\|a-b\|^2 \le 2a^2 + 2b^2$ and Assumption 2.1 (iii)-(iv), we have

$$T_{2,2} \le 2\mathbb{E}[\|\widetilde{\nabla} F_{\omega_k}(x_k)\|^2] + 2\mathbb{E}[\|\widetilde{\nabla} f^{(T)}_{\omega_{T,k}}(x_k)\nabla f^{(T-1)}_{\omega_{T-1,k}}(y_k^{(T-1)})\cdots\nabla f^{(1)}_{\omega_{1,k}}(y_k^{(1)})\|^2] \le \mathcal{O}(1).$$

Taking expectations on both sides of Eq.(E.1) and plugging in the upper bounds of $T_1$ and $T_2$ derived above, we obtain

$$\begin{aligned}&2\alpha_k(\mathbb{E}[F(x_{k+1})] - F^*) + \mathbb{E}[\|x_{k+1} - \Pi_{\mathcal{X}^*}(x_{k+1})\|^2] \\ \le& (1+\phi_k)\mathbb{E}[\|x_k - \Pi_{\mathcal{X}^*}(x_k)\|^2] + \mathcal{O}(\alpha_k^3) + \mathcal{O}(\alpha_k^2/\phi_k)\Big[\mathbb{E}[\|y_k^{(T-1)} - f^{(T)}(x_k)\|^2] \\ &+ \cdots + \mathbb{E}[\|y_k^{(1)} - f^{(2)}(y_k^{(2)})\|^2]\Big] + \mathcal{O}(\alpha_k^2).\end{aligned}$$

By the definition of optimally strong convexity in (3.1), we have

$$F(x_{k+1}) - F^* \ge \lambda\|x_{k+1} - \Pi_{\mathcal{X}^*}(x_{k+1})\|^2,$$

which futher implies

$$\begin{aligned}&(1+2\lambda\alpha_k)\mathbb{E}[\|x_{k+1} - \Pi_{\mathcal{X}^*}(x_{k+1})\|^2] \\ \le& (1+\phi_k)\mathbb{E}[\|x_k - \Pi_{\mathcal{X}^*}(x_k)\|^2] + \mathcal{O}(\alpha_k^2/\phi_k)\Big(\mathbb{E}[\|y_k^{(T-1)} - f^{(T)}(x_k)\|^2] + \cdots + \mathbb{E}[\|y_k^{(1)} - f^{(2)}(y_k^{(2)})\|^2]\Big) \\ &+ \mathcal{O}(\alpha_k^3) + \mathcal{O}(\alpha_k^2).\end{aligned}$$

It follows by dividing $(1+2\lambda\alpha_k) > 0$ on both sides,

$$\begin{aligned}\mathbb{E}[\|x_{k+1} - \Pi_{\mathcal{X}^*}(x_{k+1})\|^2] \le& \frac{1+\phi_k}{1+2\lambda\alpha_k}\mathbb{E}[\|x_k - \Pi_{\mathcal{X}^*}(x_k)\|^2] + +\mathcal{O}(\alpha_k^2/\phi_k)\Big(\mathbb{E}[\|y_k^{(T-1)} - f^{(T)}(x_k)\|^2] \\ &+ \cdots + \mathbb{E}[\|y_k^{(1)} - f^{(2)}(y_k^{(2)})\|^2]\Big) + \mathcal{O}(\alpha_k^3) + \mathcal{O}(\alpha_k^2).\end{aligned}$$

Choosing $\phi_k = \lambda\alpha_k - 2\lambda^2\alpha_k^2$ yields that

$$\begin{aligned}
&\mathbb{E}[\|x_{k+1} - \Pi_{\mathcal{X}^*}(x_{k+1})\|^2] \\
&\leq (1 - \frac{\lambda\alpha_k + 2\lambda^2\alpha_k^2}{1 + 2\lambda\alpha_k})\mathbb{E}[\|x_k - \Pi_{\mathcal{X}^*}(x_k)\|^2] + \mathcal{O}(\alpha_k^2) + \mathcal{O}(\frac{\alpha_k}{\lambda})\Big(\mathbb{E}[\|y_k^{(T-1)} - f^{(T)}(x_k)\|^2] \\
&\quad + \cdots + \mathbb{E}[\|y_k^{(1)} - f^{(2)}(y_k^{(2)})\|^2]\Big) \\
&= (1 - \lambda\alpha_k)\mathbb{E}[\|x_k - \Pi_{\mathcal{X}^*}(x_k)\|^2] \\
&\quad + \mathcal{O}(\alpha_k^2) + \mathcal{O}(\frac{\alpha_k}{\lambda})\Big(\mathbb{E}[\|y_k^{(T-1)} - f^{(T)}(x_k)\|^2] + \cdots + \mathbb{E}[\|y_k^{(1)} - f^{(2)}(y_k^{(2)})\|^2]\Big).
\end{aligned} \tag{E.2}$$

When $f^{(T)}$ has Lipschitz continuous gradient, applying Lemma 2.7 and Lemma 3.3, we have

$$\begin{aligned}
&\mathbb{E}[\|x_{k+1} - \Pi_{\mathcal{X}^*}(x_{k+1})\|^2] \\
\leq &(1 - \lambda\alpha_k)\mathbb{E}[\|x_k - \Pi_{\mathcal{X}^*}(x_k)\|^2] + \mathcal{O}(k^{-2a}) + \mathcal{O}(\frac{\alpha_k}{\lambda})\Big(\mathcal{O}(k^{-2a+2b_{T-1}}) \\
&+ \mathcal{O}(k^{-b_{T-1}}) + \sum_{j=1}^{T-2}\mathcal{O}(k^{-4b_j+4b_{j+1}}) + \mathcal{O}(k^{-b_j})\Big).
\end{aligned}$$

Apply Lemma B.1 to the previous inequality, we have

$$\mathbb{E}[\|x_k - \Pi_{\mathcal{X}^*}(x_k)\|^2] \leq \mathcal{O}(k^{-a}) + \mathcal{O}(k^{-2a+2b_{T-1}}) + \mathcal{O}(k^{-b_{T-1}}) + \sum_{j=1}^{T-2}\Big[\mathcal{O}(k^{-4b_j+4b_{j+1}}) + \mathcal{O}(k^{-b_j})\Big].$$

Letting $a = 1$, $b_{T-1} = \frac{2+T}{4+T}$, $b_{T-2} = \frac{1+T}{4+T}$, $\cdots$, $b_1 = \frac{4}{4+T}$, we have

$$\mathbb{E}[\|x_k - \Pi_{\mathcal{X}^*}(x_k)\|^2] \leq \mathcal{O}(k^{-4/(4+T)}),$$

which provides the convergence rate result for the optimally strongly convex $T$-level accelerated SCGD.

Next, when $f^{(T)}$ has Lipschitz continuous gradient, the first inner function is also updated by the accelerated update rule. By similar analysis as in Lemma 3.3, we have

$$\mathbb{E}[\|y_k^{(T-1)} - f^{(T)}(x_k)\|^2] \leq \mathcal{O}(k^{-4a+4b_{T-1}}) + \mathcal{O}(k^{-b_{T-1}}).$$

Plug this convergent rate into Eq.(E.2), we have

$$\begin{aligned}
&\mathbb{E}[\|x_k - \Pi_{\mathcal{X}^*}(x_k)\|^2] \\
\leq &\mathcal{O}(k^{-a}) + \mathcal{O}(k^{-4a+4b_{T-1}}) + \mathcal{O}(k^{-b_{T-1}}) + \sum_{j=1}^{T-2}\Big[\mathcal{O}(k^{-4b_j+4b_{j+1}}) + \mathcal{O}(k^{-b_j})\Big] \\
\leq &\mathcal{O}(k^{-4/(3+T)}),
\end{aligned}$$

by choosing $a = 1$,$b_{T-1} = \frac{2+T}{3+T}$, $b_{T-2} = \frac{1+T}{3+T}$, $\cdots$, $b_1 = \frac{4}{3+T}$, which completes the proof. □